\documentclass[11pt,letterpaper]{amsart}
\usepackage[dvipsnames]{xcolor}
\usepackage[unicode]{hyperref}
\hypersetup{
	linktoc=all,
	colorlinks = true,
	citecolor = Brown,
	anchorcolor = green,
	linkcolor = darkgray}
\usepackage[top=1.2in, bottom=1.2in, left=1.3in, right=1.25in]{geometry}
\numberwithin{equation}{section}
\usepackage{amsmath,amssymb,bm}

\newtheorem{theorem}{Theorem}[section]
\newtheorem{lemma}[theorem]{Lemma}
\newtheorem{proposition}[theorem]{Proposition}
\newtheorem{corollary}[theorem]{Corollary}

\theoremstyle{definition}

\usepackage{indentfirst}
\usepackage{fancyhdr}
\usepackage{color}
\usepackage{mathrsfs}
\usepackage{fancyhdr}
\theoremstyle{remark}
\newtheorem{remark}[theorem]{Remark}

\newcommand{\pa}{\partial}
\allowdisplaybreaks[3]

\begin{document}

\title[SMCF in $\mathbb{R}^d$, $d\ge3$]{Global well-posedness in the critical Besov space of the skew mean curvature flow in $\mathbb{R}^d$ : $d \ge 3$}

\author{Jie Shao}
\address{School of Mathematics and Statistics, Hunan Normal University, Changsha
	410081, China}
 \email{shaojiehn@foxmail.com(J. Shao)}

\author{Zexian Zhang}
\address{School of Mathematical Sciences, Fudan University, Shanghai 200433, China}
\email{23110840019@m.fudan.edu.cn}
 \author{Yi Zhou}
\address{School of Mathematical Sciences, Fudan University, Shanghai 200433, China}
\email{yizhou@fudan.edu.cn(Y. Zhou)}



\begin{abstract}
The skew mean curvature flow (SMCF) is a quasilinear Schr\"odinger-type geometric flow describing the motion of codimension-two submanifolds of Euclidean space; at the scaling-critical regularity the natural space for its second fundamental form is the Besov space $\dot B^{d/2-1}_{2,1}$. In this paper we prove global well-posedness of the SMCF for small initial data in this critical space, for all dimensions $d\ge3$. For $d\ge4$ the difference of two solutions obeys a weak Lipschitz bound in $\dot B^{d/2-2}_{2,1}$, while for $d=3$ the analogous bound holds in the weaker space $\dot B^{d/2-2}_{2,\infty}$. This extends the critical regularity theory, previously available only in dimensions $d\ge5$ \cite{lai_global_2025}, to the remaining cases $d=3$ and $d=4$. In a related direction, local well-posedness for quasilinear Schr\"odinger equations was obtained in \cite{ShZ26} through a generalized energy method based on a div-curl lemma, which differs from the Tataru--Kenig approach. The proof rests on two ingredients: a family of div-curl-type bilinear estimates, based on the integrated momentum law for the Schr\"odinger variable and on a div-curl lemma, which control the high-low interactions with a half-derivative gain, together with the interaction Morawetz estimate for the high-high interactions; and modified heat coordinates and a modified heat gauge, in which two correction terms are chosen so that the curvature part of the nonlinearity satisfies a cancellation, reducing it to product-type terms with a hidden derivative and thereby allowing the critical regularity to be reached.
\end{abstract}

\maketitle

\tableofcontents

\section{Introduction}

The skew mean curvature flow (SMCF) describes the motion of an oriented codimension-two
submanifold. For curves in \(\mathbb{R}^3\), it reduces to the binormal flow of vortex filaments, which models the dynamics of vortex lines in superfluids and superconductors \cite{Da06,Jer02,Lin98,LinT00,Vega15}; in this integrable setting the flow is equivalent, through the Hasimoto transform \cite{Has72}, to the focusing cubic nonlinear Schr\"odinger equation. In higher dimensions the SMCF provides a natural model for the evolution of vortex membranes \cite{HaVi04,Khe12}, and it has been studied both as a geometric flow and as a Hamiltonian system on a Grassmannian \cite{HaVi04,Khe12}; see also \cite{Ban20,Li21} concerning the binormal flow of polygonal lines and the transversal stability of Euclidean planes.

Let $F:[0,T]\times \mathbb{R}^d \to (\mathbb{R}^{d+2},g_{\mathbb{R}^{d+2}})$ be a one-parameter family of immersions. Denote by $g$ the induced metric, by $\mathbf{H}(F)$ the mean curvature, and by $J(F)$ the complex structure on the normal bundle $N\Sigma_t$. The SMCF reads
\begin{equation}\label{SMCF}
\begin{cases}
(\partial_t F)^{\perp} = J(F)\, \mathbf{H}(F), \\[2pt]
F|_{t=0} = F_0 .
\end{cases}
\end{equation}
Up to a tangential reparametrization, \eqref{SMCF} is equivalent to $\partial_t F = J(F)\mathbf{H}(F)$. The flow is invariant under the scaling $F(t,x)\mapsto \lambda^{-1} F(\lambda^2 t,\lambda x)$; accordingly, the second fundamental form $\lambda$ (introduced in Section~2) lives at the scaling-critical regularity $\dot B^{d/2-1}_{2,1}$.

From the point of view of partial differential equations, the SMCF in dimensions $d\ge2$ is a genuinely quasilinear Schr\"odinger-type flow, whose analysis requires a careful treatment of the derivative losses in the parabolic directions. Local well-posedness for the SMCF was first established by Song and Sun \cite{Song19} for compact surfaces in dimension $d=2$, and was extended by Song \cite{Song21} to compact oriented manifolds in all dimensions $d\ge3$, where both existence and uniqueness were proved. In the non-compact setting, Li \cite{Li21,Li22} studied transversal perturbations of Euclidean planes and obtained global regularity for small data as well as local well-posedness for large data. Huang and Tataru \cite{HuT22} proved small-data local well-posedness in dimensions $d\ge4$ at regularity $s>d/2$, one derivative above the scaling-critical exponent, in harmonic and Coulomb coordinates; they subsequently introduced a heat gauge which extends the small-data theory to all dimensions $d\ge2$ at the same regularity \cite{HuT24}, and removed the small-data restriction in \cite{HuT25}. On the global side, Huang, Li and Tataru \cite{HLT24} proved global regularity for small data in dimensions $d\ge4$ at regularity $s>s_d$, with $s_d=3$ for $d=4$ and $s_d>\frac{d+1}{2}+\frac{1}{2(d-1)}$ for $d\ge5$, and in \cite{lai_global_2025} this threshold was lowered to the scaling-critical Besov space $\dot B^{d/2-1}_{2,1}$ for all dimensions $d\ge5$. Consequently, at the scaling-critical regularity, global well-posedness was previously known only in dimensions $d\ge5$. The present paper addresses the remaining cases $d=3$ and $d=4$.

These developments are part of a broader program on quasilinear Schr\"odinger equations of the form
\[
\sqrt{-1}\,u_t+g^{\alpha\beta}(u)\partial^2_{\alpha\beta}u=N(u,\partial_x u),
\]
where the metric $g$ is at least quadratic at the origin and the nonlinearity $N$ is at least cubic at the origin and at most quadratic in $\partial_x u$. Local well-posedness for regular localized data was established in \cite{KPV,KPV2,KPV3,KPRV1,KPRV2}, and the theory in translation-invariant Sobolev spaces was later developed in \cite{MMT1,MMT2,MMT3}, with the ultrahyperbolic case treated in \cite{PT23}. For cubic nonlinearities these results give local well-posedness at regularity $s>(d+3)/2$. The regularity threshold has recently been lowered substantially by Ifrim and Tataru. For the cubic quasilinear Schr\"odinger equation in dimensions $d\ge2$, they proved local well-posedness at regularity $s>(d+1)/2$ \cite{ifrim_Global_2024}---a full derivative below the previous results and only half a derivative above the scaling exponent $d/2$, which is conjectured to be sharp. For $d\ge3$ they further obtained small-data global well-posedness and scattering at the same regularity, without additional structural assumptions; in two dimensions the global result holds at $s>3/2$ under a conservative condition and at $s\ge2$ without it. In one dimension they proved local well-posedness at the sharp threshold $s>1$, together with small-data global well-posedness under additional phase rotation, conservative and defocusing assumptions \cite{ifrim_global_2025a}; together with Pineau, they also extended these results to the ultrahyperbolic setting \cite{ifrim_Global_2025}. Using an approach different from those of Tataru et al. and Kenig et al., Shao and Zhou  \cite{ShZ26} proved local well-posedness for quasilinear Schr\"odinger equations  by developing a generalized energy method based on a div-curl lemma. While this threshold is expected to be sharp for the general class, the SMCF possesses additional geometric cancellations, described below, which allow us to reach the scaling-critical regularity $d/2-1$.

The main result of this paper is the global well-posedness of the SMCF in the critical Besov space for all dimensions $d\ge3$, under the modified heat gauge. The precise statement is given in Theorem~\ref{thm:boot}, which is formulated for the differentiated system \eqref{eq_lam}, \eqref{eq_g}--\eqref{eq_A} in the modified heat coordinates and modified heat gauge: for initial data satisfying the compatibility conditions and the harmonic coordinate and Coulomb gauge conditions, with
\[
\|(\pa_x h_0,\lambda_0,A_0)\|_{\dot B^{\frac d2-1}_{2,1}}\le \varepsilon\ll1,
\]
the system admits a global-in-time solution obeying
\[
\|(\pa_x h,\lambda,A)\|_{C_t\dot B^{\frac d2-1}_{2,1}}\lesssim \varepsilon .
\]
For $d\ge4$ the difference of two solutions satisfies the weak Lipschitz bound in $\dot B^{\frac d2-2}_{2,1}$, while for $d=3$ the analogous bound holds in the weaker space $\dot B^{\frac d2-2}_{2,\infty}$. In particular, the SMCF is well-posed in the sense of Hadamard for all dimensions $d\ge3$. The passage from the SMCF itself to the differentiated system, including the choice of the initial harmonic coordinates and Coulomb gauge, is described in Remark~\ref{rem:initial_data}.

Two ideas are responsible for reaching the critical regularity uniformly in $d\ge3$. The first is a family of div-curl-type bilinear estimates. In the Bony paraproduct decomposition, the nonlinear terms split into high-low and high-high interactions. For the high-low interactions, we use a div-curl lemma \cite{Zhou23,WaZ24} to convert the pairing of the integrated momentum law for the Schr\"odinger variable with the integrated balance laws for the low-frequency variables into bilinear space-time estimates of the form
\[
\begin{aligned}
\mathfrak{B}(j,k;(\pa_x h,A),\lambda)&\lesssim 2^{-(\frac d2-1)j}2^{-(\frac{d-1}{2}+\sigma)k}c_jc_k(\sigma),\qquad k\ge j,\\
\mathfrak{B}(j,k;\lambda,\lambda)&\lesssim 2^{-(\frac d2-1)j}2^{-(\frac{d-1}{2}+\sigma)k}c_jc_k(\sigma),\qquad k\ge j+50,
\end{aligned}
\]
where $\mathfrak{B}(j,k;\cdot,\cdot)$ denotes the bilinear norm introduced in Section~4 and $\{c_k(\sigma)\}$ the frequency envelope of \cite{Tao01}; the distinguished $x_\gamma$-direction provides the half-derivative gain. For the high-high interactions, we use the interaction Morawetz estimate \cite{LiS78,RyV07}, in the spirit of \cite{ifrim_Global_2024}, which yields
\[
\|\lambda_k\|_{L^4_{t,y_J}L^2_{z_J}}\lesssim 2^{-(\frac d2-\frac54+\sigma)k}c_k(\sigma).
\]
The two mechanisms are complementary and, together with the frequency-envelope formulation, provide the bilinear estimates needed in the bootstrap uniformly for all $d\ge3$.

The second idea is the modified heat coordinates and the modified heat gauge, which we introduce in this paper. Building on, but differing essentially from, the heat coordinates and heat gauge of Huang and Tataru \cite{HuT24,HuT25}, we add two correction terms,
\[
V^\gamma=g^{\alpha\beta}\Gamma^\gamma_{\alpha\beta}+W^\gamma,\qquad
B=\nabla^\alpha A_\alpha+\mathcal{F},
\]
with, schematically,
\[
W\sim\bigl(\bar\lambda_{\mathrm{hi}}\cdot|D|^{-1}\lambda_{\mathrm{hi}}\bigr)_{\mathrm{low}},\qquad
\mathcal{F}\sim\mathscr{R}^2(\lambda\cdot\lambda).
\]

Both corrections are introduced for the difference variables; the correction $W$ is the key to the weak Lipschitz bound. Under the modified heat coordinate condition, the metric satisfies the parabolic equation
\[
(\partial_t-\partial^\gamma\partial_\gamma)g
=\operatorname{Ric}+\operatorname{Im}(\lambda\wedge\bar\lambda)+\nabla W+\text{lower-order terms}.
\]
The difficult term occurs in the endpoint case $d=3$. To localize at frequency $k$, we apply the Littlewood--Paley operator $P_k$ (or $P_{k,\gamma}$), as in the derivation of \eqref{exp_dlamk}. The difference of the quasilinear term $g^{\alpha\beta}\partial^2_{\alpha\beta}\lambda$ then contains the term $L(2^{2k}\delta h,\lambda_k)$: the coefficient difference $\delta h$ is at frequency much lower than $k$, while the two derivatives act on the factor of frequency $2^k$, so an extra factor $2^{2k}$ appears. Estimating it requires substituting $\delta g$ through Duhamel's formula, which leads to a double bilinear estimate; the obstruction is the term $ \operatorname{Ric}+\operatorname{Im}(\lambda\wedge\bar\lambda) $, in which the two factors have comparable frequencies of order $2^k$; for comparable frequencies the div-curl lemma gives no gain. This is where the correction $W$ enters: together with an auxiliary tensor $\mathcal{H}$, it turns the curvature contribution into a single cancellation identity (Lemma~\ref{hidden_curvature}), which for the high-high interactions takes the simple form
\[
\operatorname{Ric}+\operatorname{Im}(\lambda\wedge\bar\lambda)+\nabla W
\sim(\partial_t-\partial^\gamma\partial_\gamma)\mathcal{H},
\qquad
\mathcal{H}\sim\bigl(\bar\lambda_{\mathrm{hi}}\cdot|D|^{-2}\lambda_{\mathrm{hi}}\bigr)_{\mathrm{low}}.
\]
Thus the difficult term becomes the parabolic term $(\partial_t-\partial^\gamma\partial_\gamma)\mathcal{H}$, which is moved to the left-hand side. It then suffices to subtract $\mathcal{H}$ from $h$: the two inverse derivatives in $\mathcal{H}$ cancel the factor $2^{2k}$, and for $\tilde h:=h-\mathcal{H}$ the equation has the parabolic left-hand side $(\partial_t-\Delta)\tilde h$, while the remaining nonlinearities are harmless. This is what allows the critical Besov regularity to be reached (and prevents the logarithmic loss for $d=3$); it is also the structure used in the difference estimates, where the high-high contribution is absorbed into $\delta\mathcal{H}$ and the remaining terms are controlled by the heat kernel through Duhamel's formula (see Section~5).

The role of $\mathcal{F}$ is analogous for the connection coefficients. It cancels the curl-free part of the nonlinear terms in the parabolic equation for $A$; consequently the equation for $A_{\mathrm{cf}}$ contains no term with an unbalanced derivative, and all its nonlinearities are of product type and are controlled by the $L^1_tL^2_x$ estimates. The difference $\delta A_{\mathrm{cf}}$ inherits the same structure, which is essential for the weak Lipschitz bound for the difference variables. The choices of $W$ and $\mathcal{F}$ are the key new ingredients of the construction; when $W=\mathcal{F}=0$, the modified gauge reduces exactly to the original heat gauge of \cite{HuT24,HuT25}.

The proof combines the frequency-envelope formulation of the bootstrap (Section~4), the curvature structure and the schematic form of the equations (Section~5), the div-curl and interaction Morawetz estimates (Section~6), and the nonlinear estimates (Section~7), and is concluded by closing the bootstrap in Section~8. 

The paper is organized as follows. Section~2 derives the differentiated system and the modified heat gauge; Section~3 collects the notation and the Littlewood--Paley and multilinear calculi; Section~4 states the frequency-envelope form of the bootstrap; Section~5 establishes the curvature cancellation and the schematic form of the equations; Section~6 develops the div-curl and interaction Morawetz estimates; Section~7 proves the nonlinear estimates; and Section~8 closes the bootstrap and completes the proof of Theorem~\ref{thm:boot}.
\section{The differentiated equations and modified heat gauge}

This section is devoted to the derivation of the differentiated equations for the skew mean curvature flow in $\mathbb{R}^{d+2}$, where the dimension $d\ge 3$ is fixed throughout. The main dynamic variable is the complexified second fundamental form $\lambda$, together with the auxiliary variables $(g,A)$, where $g$ is the induced metric on the immersed surface and $A$ denotes the connection coefficients of the normal bundle. The gauge formulation we use is the heat gauge of Huang and Tataru \cite{HuT24,HuT25}, modified by two correction terms: a term $W$ in the heat coordinate condition and a term $\mathcal{F}$ in the heat gauge condition (see \eqref{heat_coordinate} and \eqref{heat_gauge} below). It consists of modified heat coordinates together with a modified heat flow for the orthonormal frame of the normal bundle, and it reduces to the original heat gauge of \cite{HuT24,HuT25} when $W=\mathcal{F}=0$. We follow the derivation of \cite{HuT24,HuT25} closely and include the details below only to make the paper self-contained and to fix the notation used in the sequel.

\subsection{Compatibility conditions and evolution equations}
We begin by introducing the necessary notation, referring to \cite[Section 2]{HuT24} for more details. Let $F$ denote a one-parameter family of immersions
\[
F : \mathbb{R}_+ \times \mathbb{R}^d \rightarrow (\mathbb{R}^{d+2}, g_{\mathbb{R}^{d+2}}).
\]
For each $t \in \mathbb{R}^+$, we denote the submanifold by $\Sigma_t = F(t,\mathbb{R}^d)$, its tangent bundle by $T\Sigma_t$, and its normal bundle by $N\Sigma_t$. We use $\alpha,\beta,\gamma,\sigma\in\{1,\dots,d\}$ for coordinate indices. The induced metric on $\Sigma_t$ and the corresponding Christoffel symbols are given by
\begin{equation}\label{def_g}
	g_{\alpha\beta} = \partial_{\alpha}F \cdot \partial_{\beta}F,\qquad \Gamma^{\gamma}_{\alpha \beta} = g^{\gamma \sigma}\pa^2_{\alpha \beta}F \cdot \pa_{\sigma}F.
\end{equation}
We use $\nabla_\alpha$ to denote the associated Levi-Civita connection on $\Sigma_t$. The mean curvature $\mathbf{H}(F)$ can be naturally identified with a section of the normal bundle, and by the Gauss relation,
\begin{equation}
  \mathbf{H}(F) = \Delta_g F,
\end{equation}
where $\Delta_g$ is the Laplace--Beltrami operator associated with $g$.

Since the normal bundle is a rank-two vector bundle, it is naturally equipped with the complex structure $J(F)$, which rotates a vector in the normal space counterclockwise by $\pi/2$. We choose $\{\nu_1,\nu_2\}$ to be an orthonormal frame of $N\Sigma_t$ such that
\[
  \nu_2 = J(F) \nu_1,\qquad \nu_1 = -J(F) \nu_2,
\]
and complexify it as $m := \nu_1 + \sqrt{-1} \nu_2$. Define the tensors $\kappa,\tau$, the connection coefficients $A$ and the temporal component $B$ of the connection in the normal bundle by
\begin{equation}\label{def_A_B}
\kappa_{\alpha \beta} := \partial_{\alpha}F_{\beta} \cdot \nu_1, \quad \tau_{\alpha \beta} := \partial_{\alpha}F_{\beta} \cdot \nu_2, \quad A_{\alpha} := \partial_{\alpha}\nu_1 \cdot \nu_2, \quad B := \partial_{t}\nu_1 \cdot \nu_2,
\end{equation}
where $F_\alpha := \pa_\alpha F$, and complexify the second fundamental form as
\begin{equation}\label{def_lam}
\lambda_{\alpha \beta} := \kappa_{\alpha \beta} +\sqrt{-1} \tau_{\alpha \beta}, \qquad \psi := \operatorname{Tr} \lambda = g^{\alpha\beta}\lambda_{\alpha\beta}.
\end{equation}
The choice of the orthonormal frame is not unique: for a real-valued function $\theta$, the $\mathrm{U}(1)$ gauge group acts by
\begin{equation}\label{gauge_action}
  \psi \to  e^{\sqrt{-1}\theta}\psi,\qquad \lambda\to e^{\sqrt{-1}\theta}\lambda,\qquad m \to  e^{\sqrt{-1}\theta}m,\qquad A_\alpha \to  A_\alpha - \pa_\alpha \theta.
\end{equation}

Differentiating the frame yields the structure equations
\begin{equation}\label{struc}
\begin{cases}
	\nabla_{\alpha}F_{\beta} = \operatorname{Re}(\overline{m} \lambda_{\alpha \beta}), \\
	\partial_{\alpha}^A m = -\lambda_{\alpha}^\gamma F_\gamma, \qquad (\partial_{\alpha}^A := \partial_{\alpha} + \sqrt{-1} A_{\alpha}).
\end{cases}
\end{equation}
The compatibility conditions, in the complexified form of the classical Gauss--Codazzi--Ricci equations, then read
\begin{align}
	R_{\sigma \gamma \alpha \beta} &= \operatorname{Re}(\lambda_{\beta \gamma} \overline{\lambda}_{\alpha \sigma} - \lambda_{\alpha \gamma} \overline{\lambda}_{\beta \sigma}), \label{gauss}\\
	\nabla_{\alpha}^A \lambda_{\beta\gamma } &= \nabla_{\beta}^A \lambda_{\alpha\gamma }, \qquad (\nabla_{\alpha}^A := \nabla_{\alpha} + \sqrt{-1} A_{\alpha}), \label{codazzi}\\
	\nabla_{\alpha}A_{\beta} - \nabla_{\beta}A_{\alpha} &= \operatorname{Im}(\lambda_{\alpha}^\gamma \overline{\lambda}_{\beta \gamma}). \label{Ricci_eq}
\end{align}

We next record the evolution equations. Denoting $V_\gamma := \partial_t F\cdot F_\gamma$, we rewrite the SMCF equation \eqref{SMCF} as
\begin{equation}\label{SMCF_dF}
  \pa_t F = -\operatorname{Im}(\psi \bar{m}) + V^\gamma F_\gamma.
\end{equation}
Together with the structure equations \eqref{struc} and the orthogonality relation $m\perp F_\alpha$, this gives the equations of motion for the frame,
\begin{equation}\label{eq_frame}
	\left\{
\begin{aligned}
	&\partial_{t}F_{\alpha} = -\operatorname{Im}( \partial_{\alpha}^{A} \psi \cdot \overline{m} - \sqrt{-1} \lambda_{\alpha \gamma} V^ \gamma  \cdot \overline{m} ) + \operatorname{Im}( \psi \bar{\lambda}_{\alpha}^{\gamma} + \nabla_{\alpha} V^\gamma ) F_{\gamma}, \\
	&\partial_{t}^{B} m = -\sqrt{-1} ( \partial^{A,\alpha} \psi - \sqrt{-1} \lambda_{\gamma}^{\alpha} V^\gamma ) F_{\alpha}, \\
	&\pa_t^B \lambda^\sigma_\alpha + \lambda_\alpha^\gamma(\operatorname{Im}(\psi \bar{\lambda}_\gamma^\sigma)+\nabla_\gamma V^\sigma) = \sqrt{-1}\nabla^A_\alpha(\pa^{A,\sigma}\psi-\sqrt{-1} \lambda^\sigma_\gamma V^\gamma),
\end{aligned}
\right.
\end{equation}
where $\partial_{t}^{B} = \partial_{t} + \sqrt{-1}B$. By the definition \eqref{def_g} of $g$ and the first equation in \eqref{eq_frame}, the induced metric evolves according to
\begin{equation}\label{dt_g}
  \partial_{t} g_{\alpha \beta} = 2 \operatorname{Im} ( \bar{\lambda}_{\alpha \beta}\psi  )+ \nabla_{\alpha }  V_\beta + \nabla_{\beta}  V_\alpha,
\end{equation}
and, using the relation $[\pa_t^B, \pa_\alpha^A] m = \sqrt{-1}(\pa_t A_\alpha - \pa_\alpha B)m$, we obtain the compatibility condition for $A$,
\begin{equation}\label{dt_A}
  	\partial_{t}A_{\alpha} - \partial_{\alpha}B = \operatorname{Re}(\lambda_{\alpha}^\gamma \partial_\gamma^A \overline{\psi}) - \operatorname{Im}(\lambda_{\alpha}^\gamma \overline{\lambda}_{\gamma \sigma}) V^\sigma.
\end{equation}
Finally, the last equation in \eqref{eq_frame}, combined with the Codazzi equations \eqref{codazzi} and the Gauss equations \eqref{gauss}, yields the gauge-invariant Schr\"odinger equation for $\lambda$:
\begin{equation}\label{eq_lam}
  \begin{aligned}
(\sqrt{-1} \partial_t+\nabla^\sigma\nabla_\sigma ) \lambda_{\alpha \beta}= & \sqrt{-1}(V-2 A)^\sigma \nabla_\sigma \lambda_{\alpha \beta}-\sqrt{-1} \nabla^\sigma A_\sigma \lambda_{\alpha \beta}\\
		&+\sqrt{-1} \lambda_\alpha^\gamma \nabla_\beta V_\gamma+\sqrt{-1} \lambda_\beta^\gamma \nabla_\alpha V_\gamma \\
		& +(B+ A^\sigma A_\sigma-V^\sigma A_\sigma) \lambda_{\alpha \beta}+\psi \operatorname{Re}(\lambda_{\alpha \delta} \bar{\lambda}_\beta^\delta) \\
		& -\operatorname{Re}(\lambda_{\sigma \delta} \bar{\lambda}_{\alpha \beta}-\lambda_{\sigma \beta} \bar{\lambda}_{\alpha \delta}) \lambda^{\sigma \delta}-\lambda_{\alpha \mu} \bar{\lambda}_\sigma^\mu \lambda_\beta^\sigma.
\end{aligned}
\end{equation}

\subsection{The modified heat gauge}
The gauge we adopt originates from the heat gauge introduced by Huang and Tataru \cite{HuT24,HuT25}. The original heat gauge consists of heat coordinates on $\Sigma_t$ together with a heat flow for the orthonormal frame of the normal bundle, and it works uniformly in all dimensions $d\ge 2$. This gauge has since become a standard tool in the low-regularity theory of the SMCF, and our analysis relies on it in an essential way. Based on the formulation of Huang and Tataru, we introduce a modified gauge --- both for the coordinates and for the orthonormal frame --- which incorporates two correction terms beyond the original formulation: a term $W$ in the heat coordinate condition and a term $\mathcal{F}$ in the heat gauge condition, so that
\[
  V^\gamma = g^{\alpha\beta}\Gamma^\gamma_{\alpha\beta} + W^\gamma,\qquad B = \nabla^\alpha A_\alpha + \mathcal{F}.
\]
The correction $W$ is designed to cancel certain unfavorable terms in the parabolic equation for the metric $g$, which is essential to avoid the logarithmic loss in the nonlinear estimates in the endpoint case $d=3$.  The correction $\mathcal{F}$ cancels the curl-free part of certain terms in the parabolic equation for $A$, which is essential for the low-regularity bounds for the difference variables. When $W=\mathcal{F}=0$, the gauge reduces exactly to the original one of \cite{HuT24,HuT25}. The exact expressions for $W$ and $\mathcal{F}$ are given in Lemma~\ref{hidden_curvature} and \eqref{heat_gauge}, respectively. The derivation below is included only to make the paper self-contained and to fix notation. At the initial time we use the harmonic coordinates and the Coulomb gauge of \cite{HuT22}.

At $t=0$, we require the coordinate functions $\{x_1,\dots,x_d\}$ to be globally Lipschitz solutions of the elliptic equations
\[
  \Delta_g x_\alpha = 0,\qquad \alpha = 1,\dots,d,
\]
which can be expressed in terms of the Christoffel symbols as
\begin{equation}\label{Gam}
  g^{\alpha\beta} \Gamma_{\alpha\beta}^{\gamma} = 0.
\end{equation}
Using \eqref{Gam} and the definition of the Ricci curvature, equation \eqref{Gam} can be interpreted as an elliptic system for the metric components $g_{\alpha\beta}$:
\begin{equation}\label{eq_g_0}
  \begin{aligned}
    g^{\alpha\beta} \partial_{\alpha\beta}^2 g_{\gamma\sigma}
    = &-\partial_{\gamma} g^{\alpha\beta} \partial_{\beta} g_{\alpha\sigma}
      - \partial_{\sigma} g^{\alpha\beta} \partial_{\beta} g_{\alpha\gamma}
      + \partial_{\gamma} g_{\alpha\beta} \partial_{\sigma} g^{\alpha\beta} \\
    &+ 2 g^{\alpha\beta} \Gamma_{\sigma\alpha,\nu} \Gamma_{\beta\gamma}^{\nu}
      - 2\operatorname{Re}( \lambda_{\gamma\sigma} \overline{\psi} - \lambda_{\alpha\gamma} \overline{\lambda}_{\sigma}^{\alpha} ).
  \end{aligned}
\end{equation}

We also impose the Coulomb condition on the connection coefficients $A|_{t=0} = A_0$,
\begin{equation}\label{Cou}
  \nabla^\alpha (A_0)_\alpha = 0,
\end{equation}
which, together with \eqref{Ricci_eq}, forms an elliptic system for $A_0$.

For $t>0$, we introduce the modified heat coordinates, where we require the coordinate functions $\{x_1,\dots,x_d\}$ to be globally Lipschitz solutions of the heat equations
\begin{equation}
  (\partial_t - \Delta_g - (V-W)^\gamma \partial_\gamma)x^\alpha = 0,\qquad \alpha = 1,\dots,d,
\end{equation}
which can also be expressed in terms of the Christoffel symbols as
\begin{equation}\label{heat_coordinate}
  V^\gamma = g^{\alpha\beta} \Gamma_{\alpha\beta}^{\gamma} + W^\gamma.
\end{equation}
Here $W$ denotes a correction term of the schematic form
\[
  W \sim \bigl(\bar{\lambda}_{\mathrm{hi}}\cdot |D|^{-1}\lambda_{\mathrm{hi}}\bigr)_{\mathrm{low}},
\]
whose exact expression is specified in Lemma~\ref{hidden_curvature}.  When $W=0$, this reduces exactly to the original heat coordinates introduced in \cite{HuT24,HuT25}.
With the advection field $V$ fixed by \eqref{heat_coordinate}, DeTurck's trick for the Ricci tensor turns the evolution equation \eqref{dt_g} into the parabolic equation for the metric
\begin{equation}\label{eq_g}
  \begin{aligned}
    \partial_t g_{\mu\nu} - g^{\alpha\beta} \partial_{\alpha\beta}^2 g_{\mu\nu}
    &= \nabla_\mu W_\nu + \nabla_\nu W_\mu +  2\operatorname{Im}(\psi \bar{\lambda}_{\mu\nu}) \\
    &\quad + 2 \operatorname{Re}(\lambda_{\mu \nu} \overline{\psi}  - \lambda_{\mu \sigma} \overline{\lambda}^\sigma_{\nu})
       - 2g^{\alpha\beta} \Gamma_{\mu\beta, \sigma} \Gamma^\sigma_{\alpha\nu}\\
    &\quad + \partial_\mu g^{\alpha\beta} \Gamma_{\alpha\beta, \nu} + \partial_\nu g^{\alpha\beta} \Gamma_{\alpha\beta, \mu},
  \end{aligned}
\end{equation}
see \cite[Lemma 2.4]{HuT24} and \cite[Section 2]{HuT25}.  The role of $W$ is precisely to cancel certain unfavorable high-high interactions in the term $2\operatorname{Im}(\psi\bar{\lambda}_{\mu\nu})$ on the right-hand side; see Lemma~\ref{hidden_curvature}.

For the orthonormal frame in the normal bundle, we modify the heat gauge of \cite{HuT24,HuT25} by adding a correction term to the condition $B=\nabla^\alpha A_\alpha$:
\begin{equation}\label{heat_gauge}
  \begin{aligned}
    &B = \nabla^{\alpha} A_{\alpha} + \mathcal{F},\\
    &\mathcal{F} := \delta^{\alpha \alpha'}\mathscr{R}_{\alpha'}\bigl(\mathscr{R}_{\gamma}\operatorname{Re}(\lambda^\gamma_\alpha \bar{\psi}) + \mathscr{R}_{\sigma}\operatorname{Im}(\lambda_\alpha^\gamma \bar{\lambda}^\sigma_\gamma) \bigr) + \frac{1}{2}|\psi|^2,
  \end{aligned}
\end{equation}
where $\mathscr{R}_\alpha = (-\Delta)^{-\frac{1}{2}}\partial_\alpha$ are the Riesz operators. The purpose of the correction $\mathcal{F}$ is to cancel the curl-free part of certain terms in the parabolic equation for $A$, which is essential for the weak Lipschitz bounds for the difference variables (see below). In the special case $\mathcal{F} = 0$, the gauge condition reduces exactly to the original heat gauge condition introduced in \cite{HuT24,HuT25}.

We now derive the parabolic equation for $A$. First, using the Ricci equations \eqref{Ricci_eq}, we obtain
\begin{align}
  &~\partial_\alpha (\nabla^\sigma A_\sigma) \notag\\
  =& ~ \partial_\alpha (g^{\sigma \sigma'}\partial_{\sigma'} A_\sigma) - \partial_\alpha(\Gamma^\sigma A_\sigma) \notag \\
 = & ~ \partial_\alpha g^{\sigma \sigma'}\partial_{\sigma'} A_\sigma + g^{\sigma \sigma'}\partial_{\sigma'} (\partial_\alpha A_\sigma - \partial_\sigma A_\alpha) \notag\\
  & + g^{\sigma \sigma'}\partial^2_{\sigma \sigma'}A_\alpha - \partial_\alpha(\Gamma^\sigma A_\sigma) \notag \\
  = & ~ g^{\sigma \sigma'}\partial^2_{\sigma \sigma'}A_\alpha + \partial_\alpha g^{\sigma \sigma'}\partial_{\sigma'} A_\sigma - \partial_{\sigma'}g^{\sigma \sigma'} \operatorname{Im}(\lambda_{\alpha}^\gamma \overline{\lambda}_{\sigma \gamma}) \notag \\
  &  + \partial_{\sigma'}\operatorname{Im}(\lambda_\alpha^\gamma \bar{\lambda}^{\sigma'}_\gamma) - \partial_\alpha(\Gamma^\sigma A_\sigma) \label{paraA_derive_1}.
\end{align}

Next, we use the Codazzi equations to compute
\begin{align}
  \operatorname{Re}(\lambda_\alpha^\gamma \partial_\gamma^A \bar{\psi}) &= \operatorname{Re}(\lambda_\alpha^\gamma \nabla_\gamma^A \bar{\psi}) \notag \\
  &= \nabla_\gamma \operatorname{Re}(\lambda_\alpha^\gamma \bar{\psi}) - \operatorname{Re}(\nabla_\gamma^A \lambda_\alpha^\gamma \bar{\psi}) \notag \\
  &= \nabla_\gamma \operatorname{Re}(\lambda_\alpha^\gamma \bar{\psi}) - \operatorname{Re}(\nabla_\alpha^A \psi \bar{\psi}) \notag\\
  &= \nabla_\gamma \operatorname{Re}(\lambda_\alpha^\gamma \bar{\psi}) - \frac{1}{2}\partial_\alpha |\psi|^2 \label{paraA_derive_2}.
\end{align}

Substituting \eqref{paraA_derive_1} and \eqref{paraA_derive_2} into the evolution equation \eqref{dt_A} for $A$ yields the parabolic equation for $A$:
\begin{equation}\label{eq_A}
  \begin{aligned}
    &~\partial_t A_\alpha - g^{\sigma \sigma'}\partial^2_{\sigma \sigma'}A_\alpha \\
    =& ~\partial_\alpha g^{\sigma \sigma'}\partial_{\sigma'} A_\sigma - \partial_{\sigma'}g^{\sigma \sigma'} \operatorname{Im}(\lambda_{\alpha}^\gamma \overline{\lambda}_{\sigma \gamma}) - \partial_\alpha(\Gamma^\sigma A_\sigma) - \operatorname{Im}(\lambda_{\alpha}^\gamma \overline{\lambda}_{\gamma \sigma}) V^\sigma \\
    &  + \partial_\alpha \mathcal{F} + \partial_{\sigma'}\operatorname{Im}(\lambda_\alpha^\gamma \bar{\lambda}^{\sigma'}_\gamma) + \nabla_\gamma \operatorname{Re}(\lambda_\alpha^\gamma \bar{\psi}) - \frac{1}{2}\partial_\alpha |\psi|^2.
  \end{aligned}
\end{equation}
As the preceding computation shows, the term $\partial_\alpha \mathcal{F}$ is precisely chosen to cancel the curl-free part of the combination $\partial_{\sigma'}\operatorname{Im}(\lambda_\alpha^\gamma \bar{\lambda}^{\sigma'}_\gamma) + \nabla_\gamma \operatorname{Re}(\lambda_\alpha^\gamma \bar{\psi}) - \frac{1}{2}\partial_\alpha |\psi|^2$, thereby enabling the low-regularity bounds.

\subsection{The resulting system}
Under the modified gauge, the advection field $V$ and the connection coefficient $B$ are determined by the other variables through the modified heat coordinate condition \eqref{heat_coordinate} and the modified heat gauge condition \eqref{heat_gauge}, respectively. The skew mean curvature flow is thus expressed as the quasilinear Schr\"odinger equation \eqref{eq_lam} for $\lambda$, coupled with the parabolic equations \eqref{eq_g} and \eqref{eq_A} for $(g,A)$. The initial data are prescribed as
\begin{equation}
  g(0,x) = g_0(x),\qquad A(0,x) = A_0(x),
\end{equation}
where $g_0$ satisfies the harmonic coordinate condition \eqref{Gam} and $A_0$ the Coulomb gauge condition \eqref{Cou}. To ensure the solvability of the parabolic system, we impose the asymptotic conditions
\begin{equation}
  g(\infty) = I_d,\qquad A(\infty)=0.
\end{equation}
It is convenient to write
\begin{equation}
  g = h + I_d,
\end{equation}
so that $h$ vanishes at infinity.

We conclude this section by stating the global well-posedness of the differentiated system \eqref{eq_lam}, \eqref{eq_g}--\eqref{eq_A} in the modified heat coordinates and modified heat gauge; the proof is deferred to Section~\ref{sec:proof}.

\begin{theorem}[Global well-posedness for small data under the modified gauge]\label{thm:boot}
Let $d\ge 3$. There exists $\varepsilon>0$ sufficiently small such that, for all initial data $(g_0,\lambda_0,A_0)$ satisfying the compatibility conditions \eqref{gauss}--\eqref{Ricci_eq}, the harmonic coordinate condition \eqref{Gam} and the Coulomb condition \eqref{Cou}, with
\begin{equation}
   \|(\pa_x h_0, \lambda_0,A_0)\|_{\dot B^{\frac{d}{2}-1}_{2,1}} \le \varepsilon \ll 1,
\end{equation}
the quasilinear Schr\"odinger equation \eqref{eq_lam}, coupled with the parabolic system \eqref{eq_g}--\eqref{eq_A} under the modified heat coordinate and heat gauge conditions \eqref{heat_coordinate} and \eqref{heat_gauge}, is globally well-posed in the critical Besov space regarding the variables $(h,\lambda,A)$. More precisely, it admits a unique global-in-time solution $(h,\lambda,A)$ satisfying
\begin{equation}
    \|(\pa_x h, \lambda,A)\|_{C_t \dot B^{\frac{d}{2}-1}_{2,1}} \lesssim \varepsilon.
\end{equation}
Moreover, for $d\ge4$, the difference of two solutions obeys the weak Lipschitz bound
\begin{equation}
  \|(\pa_x h^{(1)}-\pa_x h^{(2)}, \lambda^{(1)}-\lambda^{(2)},A^{(1)}-A^{(2)})\|_{C_t \dot B^{\frac{d}{2}-2}_{2,1}} \lesssim \|(\pa_x h_0^{(1)}-\pa_x h_0^{(2)}, \lambda_0^{(1)}-\lambda_0^{(2)},A_0^{(1)}-A_0^{(2)})\|_{ \dot B^{\frac{d}{2}-2}_{2,1}} .
\end{equation}
When $d=3$, the analogous bound holds with $\dot B_{2,1}$ replaced by the weaker space $\dot B_{2,\infty}$:
\begin{equation}
 \|(\pa_x h^{(1)}-\pa_x h^{(2)}, \lambda^{(1)}-\lambda^{(2)},A^{(1)}-A^{(2)})\|_{C_t \dot B^{\frac{d}{2}-2}_{2,\infty}} \lesssim \|(\pa_x h_0^{(1)}-\pa_x h_0^{(2)}, \lambda_0^{(1)}-\lambda_0^{(2)},A_0^{(1)}-A_0^{(2)})\|_{ \dot B^{\frac{d}{2}-2}_{2,\infty}}.
\end{equation}
\end{theorem}

\begin{remark}[Initial data sets]\label{rem:initial_data}
The gauge conditions imposed on the initial data in Theorem~\ref{thm:boot} are not restrictive. For any immersion whose induced metric $g=I_d+h$ and mean curvature $\mathbf H$ with coordinate $x$ are small in the original coordinates, the construction of Huang and Tataru \cite[Section 4]{HuT24,HuT25} (adapted to the present critical regularity by a standard regularization argument) yields a change of coordinates and a rotation of the normal frame such that the new initial data $(h_0,\lambda_0,A_0)$ satisfy the harmonic coordinate condition \eqref{Gam} and the Coulomb condition \eqref{Cou} in the coordinate $y$, and obey
\[
 \|(\pa_x h_0, \lambda_0,A_0)\|_{ \dot B^{\frac{d}{2}-1}_{2,1,y} } \lesssim  \|(\pa_x h,\mathbf H)\|_{  \dot B^{\frac{d}{2}-1}_{2,1,x}}.
\]
Under these conditions $(g_0,A_0)$ is determined by $\lambda_0$ through the elliptic equations obtained from \eqref{Gam} and \eqref{Cou}, so $\lambda_0$ may be regarded as the fundamental initial datum. The harmonic/Coulomb condition at $t=0$ enters the analysis only through these bounds, while the evolution for $t>0$ is governed by the modified heat coordinates and heat gauge.
\end{remark}

\section{Notations and Preliminaries}

\subsection{Notations}
\begin{enumerate}
  \item Greek indices such as $\alpha, \beta, \gamma, \ldots$ range over $\{1, \ldots,d\}$.
  \item We write $a \lesssim b$ to denote that $a$ is bounded by an absolute constant times $b$.
  \item We introduce constant $\delta'$ satisfying $  \delta' \ll 1$.
  \item We adopt the convention that expressions involving $h$, $A$, or $\lambda$ are understood to be contracted with the metric $g$ whenever the context demands it. For example, $g^{\alpha\beta}\lambda_{\alpha \beta}$ denotes such a contracted form.
  \item For a function $u(t,x)$, we denote its spatial translation by $y \in \mathbb{R}^d$ as $u^y(t,x) := u(t, x - y)$.
  \item We define the Japanese bracket $\langle a \rangle := \sqrt{1 + |a|^2}$.
  \item Let $\mathscr{R}_\alpha := (-\Delta)^{-1/2}\partial_\alpha$, $\alpha = 1,\dots,d$, denote the Riesz operators, and write $\mathscr{R} := (\mathscr{R}_1,\dots,\mathscr{R}_d)$. Then for any vector field $v$ on $\mathbb{R}^d$, the Helmholtz decomposition is given by
    \[
    v = \mathbb{P}_{\mathrm{df}}v + \mathbb{P}_{\mathrm{cf}}v,\qquad
    \mathbb{P}_{\mathrm{cf}}v := -\mathscr{R}(\mathscr{R}\cdot v),\qquad
    \mathbb{P}_{\mathrm{df}}v := v + \mathscr{R}(\mathscr{R}\cdot v),
    \]
    where $\mathbb{P}_{\mathrm{cf}}$ and $\mathbb{P}_{\mathrm{df}}$ denote the curl-free and divergence-free projections, respectively.
  \item For a fixed $\alpha$, we write $\hat{x}_\alpha := (x_\beta)_{\beta \neq \alpha}$, i.e., the set of all spatial coordinates except $x_\alpha$.
  \item We denote the spatial gradient by $\partial_x = (\partial_{1}, \ldots,\partial_{d})$. The operators $|D|$ and $\langle D \rangle$ are defined as Fourier multipliers via
    \[
    \widehat{|D|f}(\xi) := |\xi| \hat{f}(\xi), \qquad
    \widehat{\langle D \rangle f}(\xi) := \langle \xi \rangle \hat{f}(\xi) = \sqrt{1 + |\xi|^2}\hat{f}(\xi).
    \]
  \item Throughout the paper, we use the compact notation
    \[
    \mathscr{A} := (\partial_x h,\, A,\, W)
    \]
    to denote the collection of auxiliary variables that appear in the estimates.
\end{enumerate}
\subsection{Frequency decomposition and $\mathcal{T}^{(i)}_{k}$ operators}

\subsubsection{Littlewood-Paley decomposition}
We work with the standard homogeneous Littlewood-Paley frequency decomposition.

Let $\varphi$ be a smooth cutoff function supported in the ball $\{|\xi| \leq 1\}$ and equal to $1$ on $\{|\xi| \leq 1/2\}$. For $k \in \mathbb{Z}$, we define the low-frequency projections $P_{\leq k}$ and the dyadic blocks $P_k$ by
\[
  \widehat{P_{\leq k} f}(\xi) = \varphi(2^{-k}\xi) \hat{f}(\xi), \qquad
  P_k := P_{\leq k} - P_{\leq k-1}.
\]
Equivalently, $\widehat{P_k f}(\xi) = \chi_k(\xi) \hat{f}(\xi)$ with
\[
  \chi_k(\xi) := \varphi(2^{-k}\xi) - \varphi(2^{-k+1}\xi),
\]
where each $\chi_k$ is supported in the annulus $\{2^{k-2} \leq |\xi| \leq 2^k\}$ and therefore
 \[
  \sum_{k \in \mathbb{Z}} P_k= \text{Id}.
\]

For notational convenience, we introduce the following shorthand for any tensor $f$:
\begin{equation*}
  f_k := P_k f, \qquad
  f_{|\cdot - k| \leq j} := \sum_{k' = k-j}^{k+j} P_{k'} f,
\end{equation*}
and we write $f_{\leq k} := P_{\leq k} f$; the expressions $f_{j \leq \cdot \leq j'}$ and their obvious variants are defined analogously.

In view of the bilinear estimates below, we further decompose each dyadic block as
\begin{equation}
  P_k = \sum_{\gamma = 1}^d P_{k,\gamma}, \qquad
  P_{k,\gamma} := \Xi_\gamma(D) P_k,
\end{equation}
where the $\Xi_\gamma(D)$ are Fourier multipliers whose symbols $\Xi_\gamma(\xi)$ satisfy
\[
  \sum_{\gamma=1}^{d} \Xi_\gamma(\xi) = 1 \qquad \text{for all } \xi,
\]
and each $\Xi_\gamma(\xi)$ is supported in the region $\{ |\xi| \lesssim |\xi_\gamma| \}$. For convenience, we also use the shorthand
\[
  f_{k,\gamma} := P_{k,\gamma} f.
\]

\subsubsection{Bony decomposition of products}

For Schwartz functions $f$ and $g$, the Bony paraproduct decomposition reads
\begin{equation}
  \begin{aligned}
    P_k(fg) = P_k(P_{\le k} f \cdot P_{|\cdot - k| \le 2} g)
            + P_k(P_{\le k} g \cdot P_{|\cdot - k| \le 2} f)
            + \sum_{\substack{|k_1 - k_2| \le 2  \\ \max\{k_1, k_2 \}\ge k-2}} P_k( f_{k_1} g_{k_2} ),
  \end{aligned}\label{Bony_decomp}
\end{equation}
where $f_{k} := P_{k} f$ and $g_{k} := P_{k} g$.

We schematically denote the low-high and high-low interaction terms by
\begin{equation}
  L(f_{\le k}, g_{k}), \  L(g_{\le k}, f_{k})
\end{equation}
and the high-high interaction by
\begin{equation}
  \Pi (f, g) := \sum_{|k_1 - k_2| \le 2} L( f_{k_1} g_{k_2} ),\quad \Pi_k(f, g) := \sum_{\substack{|k_1 - k_2| \le 2, \\ \max\{k_1, k_2 \}\ge k-2}} P_k L( f_{k_1} g_{k_2} ),
\end{equation}
which may be rewritten schematically as
\begin{equation}
  \Pi_k(f, g) = \sum_{k' \ge k} P_k L( f_{k'} g_{k'} ),
\end{equation}
where the latter expression employs a slight abuse of notation: here and below, $L$ denotes a generic multilinear expression of product type (Section~\ref{Lnotation}), which may vary from line to line and absorbs the additional frequency factors that arise in the nonlinear decompositions below.  In the same spirit, for $j\in\mathbb{Z}$ we write
\begin{equation}
  \Pi_{\ge j}(f, g) := \sum_{k \ge j} \Pi_k(f, g),
\end{equation}
so that $\Pi_{\ge j}(f,g)$ denotes the high-high interaction part whose frequencies are at least $2^j$.
\subsubsection{The $\mathcal{T}^{(i)}_{k}$ operators}

For $i = 1, \ldots ,d$ and $k \ge 0$, we define the operators $\mathcal{T}^{(i)}$ and $\mathcal{T}^{(i)}_{k}$ by
\begin{equation}
  \mathcal{T}^{(i)} = \Delta^{-1} \partial_i, \qquad
  \mathcal{T}^{(i)}_{k} = \mathcal{T}^{(i)} \circ P_k.
\end{equation}
It is clear that
\begin{equation}
  \mathcal{T}^{(i)} = \sum_{k \in \mathbb{Z}} \mathcal{T}^{(i)}_{k}, \qquad
  \mathcal{T}_k^{(i)} \partial_i = P_k,
\end{equation}
where the summation over $i$ in the second identity is taken directly (not the Einstein summation of an index raised with the metric $g$).
\subsection{Multilinear expressions} \label{Lnotation}

Following \cite{ifrim_Global_2024, Tao01}, we introduce a convenient notation for describing multilinear expressions of product type.

Throughout the paper, we denote by $L$ any multilinear operator of the form
\begin{equation}
  L(u_1, \dots, u_m)(x) = \int_{\mathbb{R}^{dm}} K(y_1, \dots, y_m) \, u_1^{y_1}(x) \cdots u_m^{y_m}(x) \, \mathrm{d}y_1 \cdots \mathrm{d}y_m,
\end{equation}
where $K$ is either an integrable kernel or, more generally, a bounded measure (thus encompassing, for example, product-type expressions). The specific kernel may vary from line to line, but we always assume that its total mass is uniformly bounded.

The introduction of this $L$-notation serves two main purposes.

First, it provides a unified framework for expressing all nonlinear terms appearing in the sequel, thereby simplifying the presentation of the equations and the corresponding estimates.

Second, and more importantly, the notation naturally accommodates spatial translations of individual factors via the shifts $u^{y_j}(x)=u(x-y_j)$. This is particularly useful because, unlike linear estimates, which are translation-invariant, bilinear and multilinear estimates are not invariant under independent translations of the factors. By keeping the translation parameters explicit in the kernel representation, the $L$-notation allows us to reduce the analysis to a setting where each factor is suitably translated. In this way, the $L$-notation streamlines the treatment of estimates that would otherwise be complicated by the lack of independent translation invariance.

This $L$-notation also proves useful for expressing compositions with Littlewood-Paley multipliers when the fine structure of the multipliers is not needed. For example, we can write
\begin{equation}
  \begin{aligned}
    L(P_k u_1, \dots, u_m) &= L(u_1, \dots, u_m),\\
    L(\partial_x P_{\le k} u_1, \dots, u_m) &= 2^k L(P_{\le k} u_1, \dots, u_m).
  \end{aligned}
\end{equation}
The same identities hold for $P_{>k}$, $P_{<k}$, etc. In particular, for derivatives of the homogeneous Littlewood-Paley operators $P_k$ and $P_{k,\gamma}$ (with $k \ge 1$), we have
\begin{equation}
  \begin{aligned}
    L(P_k u_1, \dots, u_m) &= 2^{-k} L(\partial_x P_k u_1, \dots, u_m),\\
    L(P_{k,\gamma} u_1, \dots, u_m) &= 2^{-k} L(\partial_\gamma P_{k,\gamma} u_1, \dots, u_m).
  \end{aligned}
\end{equation}

As noted above, multilinear estimates are not invariant under independent translations of the factors. To obtain bounds for $L$ in terms of such estimates, we consider translation-adjusted hypotheses. Specifically, suppose we have an estimate of the form
\[
\sup_{y_2,\dots,y_m} \| u_1 \, u_2^{y_2} \cdots u_m^{y_m} \|_X \le \mathcal{C},
\]
where $\|\cdot\|_X$ is a translation-invariant norm. Then for any multilinear form $L$ with an integrable kernel, we obtain
\begin{align*}
  \|L(u_1,\dots,u_m)\|_X
  &\le \int |K(y_1,\dots,y_m)| \, \| u_1^{y_1} u_2^{y_2} \cdots u_m^{y_m} \|_X \, \mathrm{d}y_1\cdots\mathrm{d}y_m \\
  &= \int |K(y_1,\dots,y_m)| \, \| u_1 u_2^{y_2} \cdots u_m^{y_m} \|_X \, \mathrm{d}y_1\cdots\mathrm{d}y_m \\
  &\le \mathcal{C} \cdot \|K\|_{L^1}.
\end{align*}
Thus, as long as the kernel $K$ has uniformly bounded total mass, the estimate for $L$ follows from the translation-invariant bound with suprema over translations.

A typical place where the $L$ notation is useful is in commutator expressions.

 \begin{lemma} Let $f,g$ be Schwartz functions.  Then
	\begin{align}\label{commu}
		 [Q_k, f]g = 2^{-k}L(\pa_x f, g),\ \text{for}\ Q_{k} = P_k,P_{k,\gamma}.
	\end{align}
\end{lemma}
 \begin{proof}
	The conclusion follows from the direct computation
		\begin{align}
		\left[P_{k}, f\right] g=-2^{-k} \int_0^1 \int_{\mathbb{R}^d}  z \chi(z) \cdot\partial_x f(x-2^{-k} s z)  g(x-2^{-k} z) \mathrm{d} z \mathrm{d} s,
	\end{align}
where $\chi$ is the kernel of the operator $Q_k$. The claim follows from the rapid decay of the kernel $\chi$.
 \end{proof}

\section{Frequency envelope formulation of bootstrap arguments}
The proof of Theorem \ref{thm:boot} is based on establishing \textit{a priori} bounds for the complex fields. To unify the treatment of linear and bilinear estimates in the bootstrap argument, we employ the language of frequency envelopes. This approach packages all the required bounds into a single set of bootstrap assumptions, enabling us to handle the nonlinear terms in a coherent way without switching between different types of estimates. In this section, we define the frequency envelope and restate our results in frequency envelope form.

\subsection{Frequency envelope}
Frequency envelopes were introduced by Tao~\cite{Tao01} and serve as a convenient tool for tracking the evolution of the energy of solutions across dyadic frequency shells; we refer to~\cite{IfT23} for a more detailed explanation.

We say that $\{c_k\}_{k \in \mathbb{Z}}$ is an $\ell^1$-frequency envelope for $u\in B_{2,1}^0$ if it satisfies the following properties:
\begin{enumerate}
    \item Slow variation: for a sufficiently small fixed $\delta' > 0$,
    \[
    c_k \lesssim 2^{\delta'|k-j|} c_j, \quad \forall k,j \in \mathbb{Z}.
    \]
    \item $\ell^p$ bounded:
    \[
   \| c_k\|_{\ell_k^p} \le  \|u\|_{L^2}.
    \]
\end{enumerate}

\subsection{The frequency envelope form of the results}

Now given  the  initial data set $(\lambda_0,h_0,A_0)$ as in Theorem \ref{thm:boot}, we define the $\ell^1$-frequency envelope $c_k(\sigma)$  by
\begin{equation}
  c_k(\sigma) = \sup_{j}2^{-\delta'|k-j|}2^{(\frac{d}{2}-1+\sigma)j}\|P_j (\pa_x h_0, \lambda_0,  A_0)\|_{L^2_x},\ \text{for}\ \sigma\ge 0.
\end{equation}
We also let $c_k := c_k(0)$ to denote the critical regularity.

For the bilinear estimates below it is convenient to record the time translation structure in a compact notation.  For $j,k\in\mathbb{Z}$ and functions $u,v$ we set
\begin{equation}\label{def_B}
  \begin{aligned}
  \mathfrak{B}(j,k;u,v):= \sup_{1\le \gamma\le d}\sup_{s\in\mathbb{R}_+}\sup_y \big(&
    \big\|\boldsymbol{1}_{\{t\ge s\}}\|u_j(t-s)\|_{L^2_{\hat x_\gamma}}\|v_{k,\gamma}^y(t)\|_{L^2_{\hat x_\gamma}}\big\|_{L^2_{t,x_\gamma}}\\
    &+ \big\|\boldsymbol{1}_{\{t\ge s\}}\|u_j(t)\|_{L^2_{\hat x_\gamma}}\|v_{k,\gamma}^y(t-s)\|_{L^2_{\hat x_\gamma}}\big\|_{L^2_{t,x_\gamma}} \big).
  \end{aligned}
\end{equation}

To prove Theorem \ref{thm:boot}, we make the following bootstrap assumptions, with $C = \varepsilon^{-\frac{1}{16}}$ (so that $1\ll C \le \varepsilon^{-1}$):
\begin{enumerate}
	\item Energy bounds:
	\begin{equation}
	  \| P_k(\partial_x h, A,\lambda) \|_{L^\infty_tL^2_x}
  \leq C 2^{-(\frac{d}{2}-1+\sigma) k}c_k(\sigma), \label{BA1}\tag{BA1}
	\end{equation}
	\item $L^2_{t,x}$ bounds for elliptic and parabolic variables:
	\begin{equation}
	   	2^k\| P_k(\partial_x h, A) \|_{L^2_t L^2_x}
	 \leq C 2^{-(\frac{d}{2}-1+\sigma) k}c_k(\sigma). \label{BA2}\tag{BA2}
	\end{equation}
	\item The bilinear estimates for the pair $(\mathscr{A},\lambda)$:
	\begin{equation}
	  \mathfrak{B}(j,k;(\pa_x h, A),\lambda)\leq C^2 2^{-(\frac{d}{2}-1)j} 2^{-(\frac{d-1}{2}+\sigma)k} c_j c_k(\sigma), ~~~\text{ for } k\geq j. \label{BA3}\tag{BA3}
	\end{equation}
	\item The bilinear estimates for Schr\"odinger variables:
	\begin{equation}
	  \mathfrak{B}(j,k;\lambda,\lambda)\leq C^2 2^{-(\frac{d}{2}-1)j} 2^{-(\frac{d-1}{2}+\sigma)k} c_j c_k(\sigma), ~~~\text{ for } k\geq j+50. \label{BA4}\tag{BA4}
	\end{equation}
	\item $L^4_{t,y_J} L^2_{z_J}$ bounds for Schr\"odinger variables:
	\begin{equation}
	\sup_J   \|    \lambda_k\|_{L^4_{t,y_J} L^2_{z_J}}   \leq  C2^{-(\frac{d}{2}-\frac{5}{4} + \sigma)k} c_k(\sigma). \label{BA5}\tag{BA5}
	\end{equation}

\end{enumerate}
Here $J$ ranges over increasing $3$-tuples of indices in $\{1,\ldots,d\}$, and we write $x=(y_J,z_J)\in\mathbb{R}^3\times\mathbb{R}^{d-3}$; thus the norm in \eqref{BA5} is defined by
\[
\| \cdot \|_{L^4_{t,y_J} L^2_{z_J}} := \big\| \| \cdot \|_{L_{z_J}^2} \big\|_{L_{t,y_{J}}^4}.\]
When $d=3$, the variable $z_J$ is absent and the norm in \eqref{BA5} is the usual $L^4_{t,x}$ norm.

To prove uniqueness and continuous dependence on the initial data, we consider two solutions $(h^{(a)},\lambda^{(a)},A^{(a)})$, $a=1,2$, of the system \eqref{eq_lam}, \eqref{eq_g}--\eqref{eq_A} in the modified heat gauge whose initial data satisfy the smallness condition in Theorem~\ref{thm:boot}, and denote their difference by
\begin{equation}\label{def_delta}
  \delta h := h^{(1)}-h^{(2)},\qquad \delta A := A^{(1)}-A^{(2)},\qquad \delta\lambda := \lambda^{(1)}-\lambda^{(2)}.
\end{equation}
The difference $\delta V$ is determined by $\delta h$ through the modified heat coordinate condition \eqref{heat_coordinate}. Subtracting the equations for the two solutions yields the linearized system: $\delta\lambda$ satisfies a linear Schr\"odinger equation whose coefficients are those of the first solution and whose source terms involve $(\delta h,\delta A)$, while $(\delta h,\delta A)$ satisfy linear parabolic equations with sources involving $\delta\lambda$. Following \cite[Section 7]{HuT24} and \cite[Section 9]{lai_global_2025}, we measure the difference one derivative below the critical regularity.

We measure the initial separation of the two solutions by the weak norm of the initial difference. For $d\ge4$, let $\{d_k\}$ be its frequency envelope; for $d=3$, the sequence $d_k$ is constant and equal to the weak norm of the initial difference itself:
\begin{equation}\label{def_d}
  d_k =\begin{cases}\sup_j 2^{-\delta'|k-j|}2^{(\frac{d}{2}-2)j}\|(\pa_x\delta h_0,\delta\lambda_0,\delta A_0)\|_{L^2_x} + \varepsilon^{-1}\|(\pa_x\delta h_0,\delta\lambda_0,\delta A_0)\|_{\dot{B}^{\frac{d}{2}-1}_{2,1}}c_k, & d \ge 4,\\ \| (\pa_x\delta h_0,\delta\lambda_0,\delta A_0)\|_{\dot{B}^{-\frac{1}{2}}_{2,\infty}}, & d=3.   \end{cases}
\end{equation}
  As in \cite[Section 9]{lai_global_2025}, the corresponding elliptic differences at $t=0$ are controlled by the same quantity, since the elliptic constraint map is locally Lipschitz in the weak topology. In addition to \eqref{BA1}--\eqref{BA5}, we assume the following bounds for the difference variables, which are the analogues of \eqref{BA1}--\eqref{BA5} with one derivative lost for the difference factor and with no additional regularity parameter:
\begin{align}
  &\|P_k(\pa_x\delta h,\delta A,\delta\lambda)\|_{L^\infty_tL^2_x} \leq C2^{-(\frac{d}{2}-2)k}d_k, \label{dBA1}\tag{\ensuremath{\delta}BA1}\\
  &2^k\|P_k(\pa_x\delta h,\delta A)\|_{L^2_{t,x}} \leq C^2 2^{-(\frac{d}{2}-2)k}d_k, \label{dBA2}\tag{\ensuremath{\delta}BA2}\\
  &\mathfrak{B}(j,k; (\pa_x \delta h, \delta A),\lambda)\leq C^2 2^{-(\frac{d}{2}-2)j}2^{-\frac{d-1}{2}k}d_jc_k, ~~~\text{for } k\geq j, \label{dBA3a}\tag{\ensuremath{\delta}BA3a}\\
  &\mathfrak{B}(j,k;(\pa_x h, A),\delta\lambda)\leq C^2 2^{-(\frac{d}{2}-1)j}2^{-\frac{d-3}{2}k}c_jd_k, ~~~\text{for } k\geq j, \label{dBA3b}\tag{\ensuremath{\delta}BA3b}\\
  &\mathfrak{B}(j,k;\delta\lambda,\lambda)\leq C^2 2^{-(\frac{d}{2}-2)j}2^{-\frac{d-1}{2}k}d_jc_k, ~~~\text{for } k\geq j+50, \label{dBA4a}\tag{\ensuremath{\delta}BA4a}\\
  &\mathfrak{B}(j,k;\lambda,\delta\lambda)\leq C^2 2^{-(\frac{d}{2}-1)j}2^{-\frac{d-3}{2}k}c_jd_k, ~~~\text{for } k\geq j+50, \label{dBA4b}\tag{\ensuremath{\delta}BA4b}\\
  &\sup_J\| \delta\lambda_k\|_{L^4_{t,y_J}L^2_{z_J}} \leq C2^{-(\frac{d}{2}-\frac{9}{4})k}d_k, \label{dBA5}\tag{\ensuremath{\delta}BA5}
\end{align}

Our goal is to improve the constants in the bounds \eqref{BA1}--\eqref{BA5} and \eqref{dBA1}--\eqref{dBA5}: it suffices to show that each of them holds with $C$ replaced by $C/2$. Since the initial data are bounded independently of $C$, it remains to control the nonlinear terms, each of which gains an extra factor $\varepsilon$; a standard continuity argument then closes the bootstrap. Once it is closed, \eqref{BA1}--\eqref{BA5} give the bounds in Theorem~\ref{thm:boot}, while \eqref{dBA1}--\eqref{dBA5} give the weak Lipschitz dependence
\begin{equation}\label{dependence}
  \|(\pa_x \delta h,\delta\lambda,\delta A)\|_{L^\infty_t\dot B^{\frac{d}{2}-2}_{2,p}} \lesssim\|(\pa_x\delta h_0,\delta\lambda_0,\delta A_0)\|_{\dot B^{\frac{d}{2}-2}_{2,p}}, \quad p = \begin{cases}
    1,\ d\ge 4, \\ \infty,\ d = 3,
  \end{cases}
\end{equation}
which yields uniqueness and continuous dependence in the weak topology for $d\ge3$. We summarize this reduction as follows.
\begin{corollary}\label{cor:boot}
	It suffices to prove Theorem~\ref{thm:boot} under the bootstrap assumptions \eqref{BA1}--\eqref{BA5} and \eqref{dBA1}--\eqref{dBA5}.
\end{corollary}

\section{Curvature structure and schematic form of the equations}
In this section we uncover the derivative structure hidden in the curvature bilinears and use it to construct the correction term $W$ in the modified heat coordinate condition.  At the high-frequency level, the bilinear form $\lambda\wedge\bar{\lambda}$ can be rewritten as a divergence together with lower-order terms (Lemma~\ref{hidden_1}); moreover, its real part, which accounts for the curvature contribution to the evolution of $g$, admits a double-divergence structure, while the imaginary part $\operatorname{Im}(\bar{\lambda}_{\alpha\beta}\psi)$ can be converted into a time derivative.  Combining these two structures, we obtain a single cancellation identity in which the correction $W$ absorbs the symmetric gradient terms and an auxiliary tensor $\mathcal{H}$ the remaining $(\pa_t-\Delta)$-term (Lemma~\ref{hidden_curvature}).  The section closes by collecting the schematic form of the equations.
\subsection{Derivative structure in curvature tensors}
We introduce the shorthand notation
\begin{equation}\label{wedge}
  \lambda \wedge \bar{\lambda} = \lambda_{\alpha \beta}\bar{\lambda}_{\gamma \delta} - \lambda_{\alpha \delta} \bar{\lambda}_{\gamma \beta},
\end{equation}
and observe that the key geometric quantities can be expressed compactly in terms of this operation:
\begin{align}
  & R_{\alpha \beta \sigma \gamma}  = \operatorname{Re}(\lambda_{\beta \gamma} \overline{\lambda}_{\alpha \sigma} - \lambda_{\beta\sigma } \overline{\lambda}_{\alpha \gamma}) = \operatorname{Re}(\lambda \wedge \bar{\lambda}),\label{wedge_R}\\
  &\operatorname{Im}(\lambda^\gamma_\alpha \bar{\lambda}_{\gamma \beta})  =\frac 12 \operatorname{Im}(\lambda^\gamma_\alpha \bar{\lambda}_{\gamma \beta}-\lambda^\gamma_\beta \bar{\lambda}_{\gamma \alpha}) = \operatorname{Im}(\lambda \wedge \bar{\lambda}),\label{wedge_Im}\\
  & \operatorname{Im}(\bar{\lambda}_{\alpha\beta}\psi ) = \operatorname{Im}(\bar{\lambda}_{\alpha\beta}\psi  - \bar{\lambda}_{\alpha}^\gamma \lambda_{\gamma \beta}) + \operatorname{Im}(\bar{\lambda}_{\alpha}^\gamma \lambda_{\gamma \beta}) = \operatorname{Im}(\lambda \wedge \bar{\lambda}).\label{wedge_psi}
\end{align}
These identities show that the curvature tensor $R_{\alpha \beta \gamma \delta}$, the quantities $\operatorname{Im}(\bar{\lambda}_{\alpha\beta}\psi)$ appearing in the metric evolution, and the term $\operatorname{Im}(\lambda^\gamma_\alpha\bar{\lambda}_{\gamma\beta})$ appearing in the compatibility condition can all be expressed in terms of the single bilinear form $\lambda \wedge \bar{\lambda}$. The lemmas below extract an additional derivative from this form at the high-frequency level.
\begin{lemma}\label{hidden_1}
  Let $\mathcal{S}$ be any symmetric subset of $\{(k_1,k_2)\mid k_1,k_2\in\mathbb{Z}\ \text{and}\ |k_1-k_2|\le 4\}$. Then
  \begin{equation}
    \begin{aligned}
     &~\sum_{(k_1,k_2)\in\mathcal{S}} \lambda_{k_1} \wedge \bar{\lambda}_{k_2} \\
     =  &~\sum_{(k_1,k_2)\in\mathcal{S}} \Bigl[
            \pa_x  L(|D|^{-1}\lambda_{k_1}\cdot  \lambda_{k_2})
           + L\big(|D|^{-1}\lambda_{k_1}, P_{k_2} ( \lambda \cdot\mathscr{A} ) \big) \Bigr].
    \end{aligned}
  \end{equation}
\end{lemma}
\begin{proof}
It suffices to analyze the dyadic interaction
\begin{equation}
    P_{k_1}\lambda_{\alpha \beta}\ P_{k_2}\bar{\lambda}_{\gamma \delta} - P_{k_1}\lambda_{\alpha \delta}\ P_{k_2}\bar{\lambda}_{\gamma \beta}
\end{equation}
for frequencies $|k_1-k_2|\le 4$. The main tool is the Codazzi equation, which we write schematically as
\begin{equation}
    \pa_\alpha \lambda_{\beta \gamma } - \pa_{\beta}\lambda_{\alpha \gamma } = \lambda \cdot \mathscr{A}.
\end{equation}
Using the identity $P_k = \mathcal{T}^{(i)}_k\pa_i$ (with summation over $i$), where $\mathcal{T}^{(i)}_k := \Delta^{-1}\pa_i P_k$, we compute
\begin{align*}
    &~P_{k_1}\lambda_{\alpha \beta}\ P_{k_2}\bar{\lambda}_{\gamma \delta} - P_{k_1}\lambda_{\alpha \delta}\ P_{k_2}\bar{\lambda}_{\gamma \beta}\\
    = & ~\mathcal{T}_{k_1}^{(\mu)}\pa_\mu \lambda_{\alpha \beta  }\ P_{k_2} \bar{\lambda}_{\gamma \delta} -\mathcal{T}_{k_1}^{(\mu)}\pa_\mu \lambda_{\alpha \delta}\ P_{k_2} \bar{\lambda}_{\gamma \beta}\\
    = & ~\mathcal{T}_{k_1}^{(\mu)}\pa_\beta \lambda_{ \mu\alpha}\ P_{k_2} \bar{\lambda}_{\gamma \delta} -\mathcal{T}_{k_1}^{(\mu)} \pa_\delta\lambda_{\mu\alpha}\ P_{k_2} \bar{\lambda}_{\gamma \beta} +\mathcal{T}_{k_1}(\mathscr{A}\cdot \lambda)\cdot \lambda_{k_2} \\
    = & ~ \pa_\beta (\mathcal{T}^{(\mu)}_{k_1} \lambda_{\mu \alpha} P_{k_2}\bar{\lambda}_{\gamma \delta}) - \pa_\delta(\mathcal{T}_{k_1}^{(\mu)}  \lambda_{\mu\alpha}\ P_{k_2} \bar{\lambda}_{\gamma \beta})+\mathcal{T}_{k_1}(\mathscr{A}\cdot \lambda)\cdot  \lambda_{k_2}\\
    &  - \mathcal{T}^{(\mu)}_{k_1} \lambda_{\mu \alpha} P_{k_2}(\pa_\beta \bar{\lambda}_{\gamma \delta} - \pa_\delta \bar{\lambda}_{\gamma \beta})\\
   = & ~\pa_\beta (\mathcal{T}^{(\mu)}_{k_1} \lambda_{\mu \alpha} P_{k_2}\bar{\lambda}_{\gamma \delta}) - \pa_\delta(\mathcal{T}_{k_1}^{(\mu)}  \lambda_{\mu\alpha}\ P_{k_2} \bar{\lambda}_{\gamma \beta})+\mathcal{T}_{k_1}(\mathscr{A}\cdot \lambda)\cdot \lambda_{k_2}\\
    &   + \mathcal{T}_{k_1} \lambda \ P_{k_2}(\mathscr{A}\cdot \lambda).
\end{align*}
The first two terms are in divergence form, while the remaining terms involve the auxiliary variables $\mathscr{A}$. Summing over $(k_1,k_2)\in\mathcal{S}$ and using the symmetry of $\mathcal{S}$ to combine the contributions from $(k_1,k_2)$ and $(k_2,k_1)$ yields the claimed identity.
\end{proof}

\subsection{Cancellation via a modified gauge condition}
In this subsection we determine the correction term $W$ in the modified heat coordinate condition \eqref{heat_coordinate}.  The goal is to eliminate the unfavorable high-high interactions in $\operatorname{Im}(\bar{\lambda}_{\alpha \beta} \psi)$ that cannot be controlled by the critical norms and, at the same time, to exhibit the divergence structure of the real part, which accounts for the curvature contributions to the evolution of the metric.  As will be seen below, the derivative structure uncovered in the previous subsection plays a central role in this construction.

A further iteration of the above argument reveals a double-divergence structure for the real part, while the imaginary part can be converted into a time derivative; the two structures combine into a single cancellation identity for the sum of the two contributions, which is precisely the combination of curvature terms appearing in the evolution of $g$ (see \eqref{eq_g}).  We now state the main lemma of this subsection, which provides the precise choices of $W$ and $\mathcal{H}$ as well as the resulting identity.
\begin{lemma}\label{hidden_curvature}
   Given $\mathcal{S}$ as in Lemma~\ref{hidden_1}, there exist functions $W_\alpha[k_1,k_2]$ and $\mathcal{H}_{\alpha \beta}[k_1,k_2]$, with $W_\alpha[k_1,k_2]$ of the form
   \begin{equation}
     W [k_1,k_2] = L(|D|^{-1} \lambda_{k_1}, \bar{\lambda}_{k_2}),\quad \mathcal{H}[k_1,k_2] = L(|D|^{-2}\lambda_{k_1},\bar{\lambda}_{k_2}),
   \end{equation}
such that
\begin{equation}\label{hidden_curvature_formula}
\begin{aligned}
    &~\sum_{(k_1,k_2)\in \mathcal{S}}\Big[2\operatorname{Im}(P_{k_1}\bar{\lambda}_{\alpha \beta}P_{k_2}\psi) +   2 g^{\gamma \delta}\operatorname{Re}(P_{k_1}\lambda_{\alpha \beta}\ P_{k_2}\bar{\lambda}_{\gamma \delta} - P_{k_1}\lambda_{\alpha \delta}\ P_{k_2}\bar{\lambda}_{\gamma \beta})\\
    & \qquad\qquad\ +  \pa_\alpha W_\beta[k_1,k_2] + \pa_\beta W_\alpha[k_1,k_2]\Big] \\
     =&~ \sum_{(k_1,k_2)\in \mathcal{S}}\Big[ L\big(|D|^{-1}\lambda_{k_1}, P_{k_2}(\mathscr{A}\cdot \lambda)\big)  +  L \big( |D|^{-2}\lambda_{k_1}, P_{k_2}[(\mathscr{A}^2,\pa_x\mathscr{A},\lambda^2)\cdot  \lambda] \big)   \\
      &\qquad\qquad\ +  (\pa_t -\pa^\gamma\pa_\gamma )\mathcal{H}_{\alpha \beta}[k_1,k_2]\Big],
\end{aligned}
\end{equation}
where $\mathscr{A} := (\pa_x h, A, W)$ is the collection of auxiliary variables introduced in Section~3.  The choice  of $W$ in the modified heat gauge is then given by
\begin{equation}
  W = \sum_{|k_1-k_2|\le 4} W[k_1,k_2] .
\end{equation}
\end{lemma}

\begin{proof}
Let RHS denote the collection of lower-order terms appearing on the right-hand side of \eqref{hidden_curvature_formula}.  We begin with the expression
$$2\operatorname{Im}(P_{k_1}\bar{\lambda}_{\alpha\beta}P_{k_2}\psi) + 2 g^{\gamma \delta}\operatorname{Re}(P_{k_1}\lambda_{\alpha \beta}\ P_{k_2}\bar{\lambda}_{\gamma \delta} - P_{k_1}\lambda_{\alpha \delta}\ P_{k_2}\bar{\lambda}_{\gamma \beta})$$
summed over $(k_1,k_2)\in\mathcal{S}$; for brevity we suppress the summation in what follows.
Using Lemma~\ref{hidden_1} and the symmetry in $(\alpha,\beta)$, we obtain
\begin{align}
  &~2 g^{\gamma \delta}\operatorname{Re}(P_{k_1}\lambda_{\alpha \beta}\ P_{k_2}\bar{\lambda}_{\gamma \delta} - P_{k_1}\lambda_{\alpha \delta}\ P_{k_2}\bar{\lambda}_{\gamma \beta}) \notag\\
  = & ~g^{\gamma \delta}\operatorname{Re}[\pa_\alpha (\mathcal{T}^{(\mu)}_{k_1} \lambda_{\mu \beta} P_{k_2}\bar{\lambda}_{\gamma \delta})+\pa_\beta (\mathcal{T}^{(\mu)}_{k_1} \lambda_{\mu \alpha} P_{k_2}\bar{\lambda}_{\gamma \delta}) ]\label{item_1}\\
  &   - \pa^\gamma\operatorname{Re}(\mathcal{T}_{k_1}^{(\mu)}  \lambda_{\mu\alpha}\ P_{k_2} \bar{\lambda}_{\gamma \beta} + \mathcal{T}_{k_1}^{(\mu)}  \lambda_{\mu\beta}\ P_{k_2} \bar{\lambda}_{\gamma \alpha} )\label{item_2}\\
    & + \operatorname{RHS},\notag
\end{align}
The terms in \eqref{item_1} will be canceled by the symmetric gradient of $W$.  For \eqref{item_2} we integrate by parts repeatedly, using the Codazzi equation each time to rotate the derivative, until it acts on $\lambda$ with the summed index $\mu$, so that $\mathcal{T}^{(\mu)}_{k_1}\pa_\mu = P_{k_1}$ produces a single product:
\begin{align}
  &~\mathcal{T}^{(\mu)}_{k_1} \lambda_{\mu \alpha } P_{k_2} \bar{\lambda}_{\gamma \beta}\notag\\
   =  &~\mathcal{T}^{(\mu)}_{k_1} \lambda_{\mu \alpha } \mathcal{T}^{(\nu)}_{k_2} \pa_\nu\bar{\lambda}_{\gamma \beta}\notag\\
  =&~ \mathcal{T}^{(\mu)}_{k_1} \lambda_{\mu \alpha } \mathcal{T}^{(\nu)}_{k_2} \pa_\gamma\bar{\lambda}_{\nu \beta} + \mathcal{T}_{k_1} \lambda\cdot \mathcal{T}_{k_2}(\mathscr{A}\cdot \lambda)\notag\\
  =&~ \pa_\gamma(\mathcal{T}^{(\mu)}_{k_1} \lambda_{\mu \alpha } \mathcal{T}^{(\nu)}_{k_2}\bar{\lambda}_{\nu \beta}) - \mathcal{T}^{(\mu)}_{k_1} \pa_\gamma\lambda_{\mu \alpha } \mathcal{T}^{(\nu)}_{k_2}\bar{\lambda}_{\nu \beta} + \mathcal{T}_{k_1} \lambda\cdot \mathcal{T}_{k_2}(\mathscr{A}\cdot \lambda) \notag\\
  =&~ \pa_\gamma(\mathcal{T}^{(\mu)}_{k_1} \lambda_{\mu \alpha } \mathcal{T}^{(\nu)}_{k_2}\bar{\lambda}_{\nu \beta}) - \mathcal{T}^{(\mu)}_{k_1} \pa_\alpha \lambda_{\mu \gamma } \mathcal{T}^{(\nu)}_{k_2}\bar{\lambda}_{\nu \beta}\notag\\
  & + \mathcal{T}_{k_1} \lambda\cdot \mathcal{T}_{k_2}(\mathscr{A}\cdot \lambda) +  \mathcal{T}_{k_1}(\mathscr{A}\cdot \lambda)\cdot \mathcal{T}_{k_2} \lambda \notag\\
  =&~ \pa_\gamma(\mathcal{T}^{(\mu)}_{k_1} \lambda_{\mu \alpha } \mathcal{T}^{(\nu)}_{k_2}\bar{\lambda}_{\nu\beta}) - \pa_\alpha(\mathcal{T}^{(\mu)}_{k_1}  \lambda_{\mu \gamma } \mathcal{T}^{(\nu)}_{k_2}\bar{\lambda}_{\nu \beta}) +\mathcal{T}^{(\mu)}_{k_1}  \lambda_{\mu \gamma } \mathcal{T}^{(\nu)}_{k_2}\pa_\alpha \bar{\lambda}_{\nu\beta} \notag\\
  & + \mathcal{T}_{k_1} \lambda\cdot \mathcal{T}_{k_2}(\mathscr{A}\cdot \lambda) +  \mathcal{T}_{k_1}(\mathscr{A}\cdot \lambda)\cdot \mathcal{T}_{k_2} \lambda\notag\\
  =&~ \pa_\gamma(\mathcal{T}^{(\mu)}_{k_1} \lambda_{\mu \alpha } \mathcal{T}^{(\nu)}_{k_2}\bar{\lambda}_{\nu\beta}) - \pa_\alpha(\mathcal{T}^{(\mu)}_{k_1}  \lambda_{\mu \gamma } \mathcal{T}^{(\nu)}_{k_2}\bar{\lambda}_{\nu \beta}) +\pa_\beta(\mathcal{T}^{(\mu)}_{k_1}  \lambda_{\mu \gamma } \mathcal{T}^{(\nu)}_{k_2} \bar{\lambda}_{\nu\alpha} )\notag\\
  & -  \mathcal{T}^{(\mu)}_{k_1}  \pa_\beta\lambda_{\mu \gamma } \mathcal{T}^{(\nu)}_{k_2} \bar{\lambda}_{\nu\alpha}  + \mathcal{T}_{k_1} \lambda\cdot \mathcal{T}_{k_2}(\mathscr{A}\cdot \lambda) +  \mathcal{T}_{k_1}(\mathscr{A}\cdot \lambda)\cdot \mathcal{T}_{k_2} \lambda\notag\\
  =&~ -  P_{k_1}  \lambda_{  \gamma \beta} \mathcal{T}^{(\nu)}_{k_2} \bar{\lambda}_{\nu\alpha}  +\pa_\gamma(\mathcal{T}^{(\mu)}_{k_1} \lambda_{\mu \alpha } \mathcal{T}^{(\nu)}_{k_2}\bar{\lambda}_{\nu\beta})  - \pa_\alpha(\mathcal{T}^{(\mu)}_{k_1}  \lambda_{\mu \gamma } \mathcal{T}^{(\nu)}_{k_2}\bar{\lambda}_{\nu \beta}) \label{item_3}\\
  &+\pa_\beta(\mathcal{T}^{(\mu)}_{k_1}  \lambda_{\mu \gamma } \mathcal{T}^{(\nu)}_{k_2} \bar{\lambda}_{\nu\alpha} )   + \mathcal{T}_{k_1} \lambda\cdot \mathcal{T}_{k_2}(\mathscr{A}\cdot \lambda) +  \mathcal{T}_{k_1}(\mathscr{A}\cdot \lambda)\cdot \mathcal{T}_{k_2} \lambda,\notag
\end{align}
Taking real parts, the first term of \eqref{item_3} is the negative of a term already appearing on the left-hand side, so it can be moved to the left-hand side; substituting the resulting expression for \eqref{item_2}, we obtain
\begin{align}
 \eqref{item_2} =   -\pa^\gamma \pa_\gamma\operatorname{Re}(\mathcal{T}^{(\mu)}_{k_1} \lambda_{\mu \alpha } \mathcal{T}^{(\nu)}_{k_2}\bar{\lambda}_{\nu\beta})  + \operatorname{RHS}.\label{item_4}
\end{align}
Adding up the above identities, we obtain
\begin{equation}\label{item_5}
  \begin{aligned}
     &~ 2\operatorname{Im}(P_{k_1}\bar{\lambda}_{\alpha\beta}P_{k_2}\psi) + 2 g^{\gamma \delta}\operatorname{Re}(P_{k_1}\lambda_{\alpha \beta}\ P_{k_2}\bar{\lambda}_{\gamma \delta} - P_{k_1}\lambda_{\alpha \delta}\ P_{k_2}\bar{\lambda}_{\gamma \beta})  \notag\\
  = &~ g^{\gamma \delta}\operatorname{Re}[\pa_\alpha (\mathcal{T}^{(\mu)}_{k_1} \lambda_{\mu \beta} P_{k_2}\bar{\lambda}_{\gamma \delta})+\pa_\beta (\mathcal{T}^{(\mu)}_{k_1} \lambda_{\mu \alpha} P_{k_2}\bar{\lambda}_{\gamma \delta}) ] + (\pa_t -\pa^\gamma\pa_\gamma)\operatorname{Re}(\mathcal{T}^{(\mu)}_{k_1} \lambda_{\mu \alpha } \mathcal{T}^{(\nu)}_{k_2}\bar{\lambda}_{\nu\beta}) \\
  &+ 2\operatorname{Im}(P_{k_1}\bar{\lambda}_{\alpha\beta}P_{k_2}\psi)  -\pa_t \operatorname{Re}(\mathcal{T}^{(\mu)}_{k_1} \lambda_{\mu \alpha } \mathcal{T}^{(\nu)}_{k_2}\bar{\lambda}_{\nu\beta})
  \end{aligned}
\end{equation}

It remains to compute the term $-\pa_t \operatorname{Re}(\mathcal{T}^{(\mu)}_{k_1} \lambda_{\mu \alpha } \mathcal{T}^{(\nu)}_{k_2}\bar{\lambda}_{\nu\beta})$ appearing in \eqref{item_5}.  We recall the quasilinear Schr\"odinger equation \eqref{eq_lam} for $\lambda_{\alpha\beta}$ and, to isolate the leading-order elliptic part, expand the covariant Laplacian as
\begin{equation}
\begin{aligned}
    \nabla^A_\sigma \nabla^{A,\sigma} \lambda_{\alpha \beta}
    &= \nabla^A_\sigma \nabla_\alpha^A \lambda_\beta^\sigma
    = [\nabla_\sigma^A, \nabla_\alpha^A]\lambda^{\sigma}_\beta + \nabla_\alpha^A \nabla_\beta^A \psi \\
    &= (\lambda \wedge \bar{\lambda})\cdot \lambda + (\pa_\alpha + (\pa_x h, A))(\pa_\beta + (\pa_x h, A))\psi.
\end{aligned}
\end{equation}
Inserting this expansion into \eqref{eq_lam} and keeping only the schematic form of the lower-order corrections, we obtain the compact expression
\begin{equation}\label{lambdadt}
\begin{aligned}
    &\sqrt{-1}\pa_t\lambda_{\alpha\beta}   =- \pa^2_{\alpha \beta}\psi +  \pa_x(\mathscr{A}\cdot \lambda)
       + ( \pa_x \mathscr{A},\, \mathscr{A}^2,\,\lambda^2  )\cdot \lambda.
\end{aligned}
\end{equation}

Consequently,
\begin{align}
  &~ -\pa_t\operatorname{Re} (\mathcal{T}^{(\mu)}_{k_1} \lambda_{\mu \alpha } \mathcal{T}^{(\nu)}_{k_2}\bar{\lambda}_{\nu\beta})\notag \\
   = & -\operatorname{Im} (\sqrt{-1}\mathcal{T}^{(\mu)}_{k_1} \pa_t\lambda_{\mu \alpha } \mathcal{T}^{(\nu)}_{k_2}\bar{\lambda}_{\nu\beta} + \sqrt{-1}\mathcal{T}^{(\mu)}_{k_1} \lambda_{\mu \alpha } \mathcal{T}^{(\nu)}_{k_2}\pa_t\bar{\lambda}_{\nu\beta}) \notag\\
    =&~ \operatorname{Im} (\pa_\alpha P_{k_1}\psi\mathcal{T}^{(\nu)}_{k_2}\bar{\lambda}_{\nu\beta} - \mathcal{T}^{(\mu)}_{k_1} \lambda_{\mu \alpha } \pa_\beta P_{k_2}\bar{\psi}) + \operatorname{RHS}\notag\\
   =&~ \pa_\alpha\operatorname{Im} ( P_{k_1}\psi\mathcal{T}^{(\nu)}_{k_2}\bar{\lambda}_{\nu\beta}) - \pa_\beta\operatorname{Im} (\mathcal{T}^{(\mu)}_{k_1} \lambda_{\mu \alpha }  P_{k_2}\bar{\psi}) -  \operatorname{Im} ( P_{k_1}\psi P_{k_2} \bar{\lambda}_{\alpha\beta}  -  P_{k_1} \lambda_{ \alpha \beta }  P_{k_2}\bar{\psi}) + \operatorname{RHS} \notag\\
    = & ~ \pa_\alpha\operatorname{Im} ( P_{k_1}\psi\mathcal{T}^{(\nu)}_{k_2}\bar{\lambda}_{\nu\beta})+ \pa_\beta\operatorname{Im} (\mathcal{T}^{(\mu)}_{k_1} \bar\lambda_{\mu \alpha }  P_{k_2} \psi )  - 2 \operatorname{Im} ( P_{k_1}\psi P_{k_2} \bar{\lambda}_{\alpha\beta}  ) + \operatorname{RHS},
\end{align}
Substituting this into \eqref{item_5}, we obtain
\begin{equation}
    \begin{aligned}
     &~ 2\operatorname{Im}(P_{k_1}\bar{\lambda}_{\alpha\beta}P_{k_2}\psi) + 2 g^{\gamma \delta}\operatorname{Re}(P_{k_1}\lambda_{\alpha \beta}\ P_{k_2}\bar{\lambda}_{\gamma \delta} - P_{k_1}\lambda_{\alpha \delta}\ P_{k_2}\bar{\lambda}_{\gamma \beta})  \notag\\
  = &~ g^{\gamma \delta}\operatorname{Re}[\pa_\alpha (\mathcal{T}^{(\mu)}_{k_1} \lambda_{\mu \beta} P_{k_2}\bar{\lambda}_{\gamma \delta})+\pa_\beta (\mathcal{T}^{(\mu)}_{k_1} \lambda_{\mu \alpha} P_{k_2}\bar{\lambda}_{\gamma \delta}) ] + (\pa_t -\pa^\gamma\pa_\gamma)\operatorname{Re}(\mathcal{T}^{(\mu)}_{k_1} \lambda_{\mu \alpha } \mathcal{T}^{(\nu)}_{k_2}\bar{\lambda}_{\nu\beta}) \\
  &+ \pa_\alpha\operatorname{Im} ( P_{k_1}\psi\mathcal{T}^{(\nu)}_{k_2}\bar{\lambda}_{\nu\beta})+ \pa_\beta\operatorname{Im} (\mathcal{T}^{(\mu)}_{k_1} \bar\lambda_{\mu \alpha }  P_{k_2} \psi ),
  \end{aligned}
\end{equation}
We therefore choose
\begin{equation}
  \begin{aligned}
    &W_{\alpha}[k_1,k_2] =-g^{\gamma \delta} \operatorname{Re}(\mathcal{T}^{(\mu)}_{k_1} \lambda_{\mu \alpha} P_{k_2}\bar{\lambda}_{\gamma \delta})- \operatorname{Im} ( P_{k_1}\psi\mathcal{T}^{(\nu)}_{k_2}\bar{\lambda}_{\nu\alpha}),\\
    & \mathcal{H}_{\alpha \beta}[k_1,k_2] = \operatorname{Re}(\mathcal{T}^{(\mu)}_{k_1} \lambda_{\mu \alpha } \mathcal{T}^{(\nu)}_{k_2}\bar{\lambda}_{\nu\beta}),
  \end{aligned}
\end{equation}
This completes the proof.
\end{proof}

\subsection{Schematic form of the equations}

\subsubsection{Equations for $(h_k,A_k)$ and $(\delta h_k, \delta A_k)$}
We rewrite \eqref{eq_g} and \eqref{eq_A} schematically in the form
\begin{align}
(\pa_t - \Delta) h = &\,h \cdot  \partial_x^2 h  +  \partial_x h \cdot  \partial_x h + \pa_x W+ \lambda\wedge \bar{\lambda}, \label{gen_eq_h}\\
(\pa_t - \Delta) A = &\,h \cdot  \partial_x^2 A +  \partial_x h\cdot \pa_x A   + \lambda^2 \cdot(\partial_x h, A) +\partial_x (\lambda^2)   + \partial_x \mathcal{F}. \label{gen_eq_A}
\end{align}
 Applying $P_k$ to both sides of \eqref{gen_eq_h} and using the modified heat gauge \eqref{heat_gauge} together with the Bony decomposition, we obtain
 \begin{equation}\label{exp_hk}
   \begin{aligned}
(\pa_t - \Delta) h_k  = 2^{ k}L(h_{\le k}, \pa_x h_k) + \Pi_k(\pa_x h, \pa_x h) + P_{k}(\pa_x W + \lambda \wedge \bar{\lambda})=: \mathcal{K}(h_k).
\end{aligned}
 \end{equation}

Similarly,
\begin{equation}\label{exp_Ak}
\begin{aligned}
(\pa_t - \Delta) A_k = &\,2^{2k}L(h_{\leq k}, A_k)  + L(\pa_x h_k,\pa_x A_{\leq k}) +  \Pi_k ( \pa_x h,\pa_x A)    \\
&+ P_k( 2^k\lambda^2+ \lambda^2 \cdot \mathscr{A})  \\
=&: \mathcal{K}(A_k).
\end{aligned}
\end{equation}

We can write both equations in the unified form
\begin{equation}
  \label{exp_scrA}
  \begin{aligned}
   (\pa_t -\Delta) \mathscr{A}_k &= 2^{2k}L(h_{\le k},\mathscr{A}_k) + L( \pa_x\mathscr{A}_{\le k}, \mathscr{A}_k) + \Pi_{k}(\mathscr{A},\pa_x \mathscr{A}) + P_k(2^k \lambda^2 + \lambda^2 \cdot \mathscr{A}) \\
   & =: \mathcal{K}( \mathscr{A}_k)
  \end{aligned}
\end{equation}

We also need an alternative form of the equation for $h_k$, which will be used to derive the equations for $\delta h_k$.
It is obtained by applying Lemma~\ref{hidden_1} to the high-high interaction of $\lambda \wedge \bar{\lambda}$:
 \begin{equation}\label{exp_hk_2}
   \begin{aligned}
(\pa_t - \Delta) h_k  =  &\,  2^{ k}L(h_{\le k}, \pa_x h_k) + \Pi_k(\pa_x h, \pa_x h) + L(\lambda_{\le k}, \lambda_k)\\
& + \Pi_k(|D|^{-1}\lambda,  2^{k}\lambda +\lambda\cdot \mathscr{A} )  .
\end{aligned}
 \end{equation}

We also decompose $A$ into its divergence-free and curl-free components:
\begin{equation}
  A = A_{\mathrm{df}} + A_{\mathrm{cf}}.
\end{equation}
A direct computation shows that the nonlinear terms in \eqref{eq_A} admit the following decomposition:
\begin{align}
  &\partial_\alpha \mathcal{F} + \partial_{\sigma}\operatorname{Im}(\lambda_\alpha^\gamma \bar{\lambda}^{\sigma}_\gamma) + \pa_\gamma \operatorname{Re}(\lambda_\alpha^\gamma \bar{\psi}) - \frac{1}{2}\partial_\alpha |\psi|^2\notag\\
=&~ \partial_\alpha( \mathcal{F}-\frac{1}{2}|\psi|^2)   + \mathbb{P}_{\mathrm{cf}}\big(\partial_{\sigma}\operatorname{Im}(\lambda_\alpha^\gamma \bar{\lambda}^{\sigma}_\gamma) + \pa_\gamma \operatorname{Re}(\lambda_\alpha^\gamma \bar{\psi})\big)\notag\\
&+ \mathbb{P}_{\mathrm{df}}\big(\partial_{\sigma}\operatorname{Im}(\lambda_\alpha^\gamma \bar{\lambda}^{\sigma}_\gamma) + \pa_\gamma \operatorname{Re}(\lambda_\alpha^\gamma \bar{\psi})\big)\notag\\
=&~ \partial_\alpha\big( \mathcal{F}-\frac{1}{2}|\psi|^2 - \delta^{\alpha \alpha'}\mathscr{R}_{\alpha'}\bigl(\mathscr{R}_{\gamma}\operatorname{Re}(\lambda^\gamma_\alpha \bar{\psi}) + \mathscr{R}_{\sigma}\operatorname{Im}(\lambda_\alpha^\gamma \bar{\lambda}^\sigma_\gamma) \bigr)\big)\notag\\
&+ \mathbb{P}_{\mathrm{df}}\big(\partial_{\sigma}\operatorname{Im}(\lambda_\alpha^\gamma \bar{\lambda}^{\sigma}_\gamma) + \pa_\gamma \operatorname{Re}(\lambda_\alpha^\gamma \bar{\psi})\big)\notag\\
=& ~\mathbb{P}_{\mathrm{df}}\big(\partial_{\sigma}\operatorname{Im}(\lambda_\alpha^\gamma \bar{\lambda}^{\sigma}_\gamma) + \pa_\gamma \operatorname{Re}(\lambda_\alpha^\gamma \bar{\psi})\big),\label{A_cancel}
\end{align}
where the last equality follows from the definition \eqref{heat_gauge} of $\mathcal{F}$: the curl-free part of $\partial_{\sigma}\operatorname{Im}(\lambda_\alpha^\gamma \bar{\lambda}^{\sigma}_\gamma) + \pa_\gamma \operatorname{Re}(\lambda_\alpha^\gamma \bar{\psi})$ is exactly canceled by $\partial_\alpha\mathcal{F}$, so that the pure gradient term vanishes. Consequently, applying the curl-free projection $\mathbb{P}_{\mathrm{cf}}$ to the parabolic equation \eqref{eq_A} yields the following schematic form for $A_{\mathrm{cf}}$:
\begin{align}
  (\pa_t   - \Delta) A_{k,\mathrm{cf}} &= 2^{2k}L(h_{\le k}, A_k)+L(\pa_x \mathscr{A}_{\le  k},  \mathscr{A}_k)  + \Pi_k(\mathscr{A}, \pa_x\mathscr{A})+ P_k (   \mathscr{A}\cdot\lambda^2 ), \label{para_A_cf}
\end{align}
Furthermore, the Ricci equation \eqref{Ricci_eq}, together with Lemma~\ref{hidden_1} applied to the high-high interactions, yields the relation
\begin{equation}\label{Adf_lambda}
  \begin{aligned}
  \pa_x A_{k,\mathrm{df}}& = P_k(\lambda \wedge \bar{\lambda})   = L(\lambda_{\le k},\lambda_k)+  \pa_x\Pi_k(|D|^{-1}\lambda,  \lambda  ) +  \Pi_k(|D|^{-1}\lambda, \lambda\cdot \mathscr{A} )
  \end{aligned}
\end{equation}

Subtracting the equations for two solutions, we obtain the parabolic equation for the difference $\delta h_k$:
 \begin{equation}\label{exp_dhk}
   \begin{aligned}
(\pa_t - \Delta) \delta h_k  =  &\,  2^{ k}\big[L(\delta h_{\le k}, \pa_x h_k) +L( h_{\le k}, \pa_x \delta h_k) \big]+ \Pi_k(\pa_x \delta h, \pa_x h) \\
& +   L(\delta \lambda_{\le k}, \lambda_k)  + L(\lambda_{\le k}, \delta \lambda_k)+\Pi_k(|D|^{-1}\delta\lambda,\lambda \cdot \mathscr{A}  +2^k\lambda )\\
& +   \Pi_k(|D|^{-1}\lambda,  \delta\lambda\cdot \mathscr{A}+ \lambda\cdot \delta\mathscr{A})\\
=:&~  \mathcal{K}(\delta h_k).
\end{aligned}
 \end{equation}

For $\delta A$, we similarly obtain the equation for $\delta A_{k,\mathrm{cf}}$:
\begin{equation} \label{para_dA_cf}
  \begin{aligned}
 (\pa_t  - \Delta )\delta A_{k,\mathrm{cf}} =&~ 2^{k}(L(\delta h_{\le k},\pa_x A_k) + L( h_{\le k},\pa_x \delta A_k))+L(\pa_x \delta \mathscr{A}_{\le  k}, \mathscr{A}_k) \\
 &+ L(\pa_x \mathscr{A}_{\le  k}, \delta\mathscr{A}_k)+ \Pi_k(\delta \mathscr{A}, \pa_x\mathscr{A})+ P_k (   \delta\mathscr{A}\cdot\lambda^2  + \lambda \cdot \delta \lambda \cdot \mathscr{A})\\
=:&~  \mathcal{K}(\delta A_{k,\mathrm{cf}}).
  \end{aligned}
\end{equation}
and $\delta A_{k,\mathrm{df}}$ satisfies the relation
\begin{equation}\label{dAdf_lambda}
  \begin{aligned}
  \pa_x \delta A_{k,\mathrm{df}} = &~L(\delta \lambda_{\le k},\lambda_k) +  L( \lambda_{\le k},\delta\lambda_k) + \Pi_k(|D|^{-1}\delta\lambda,  2^{k}\lambda +\lambda\cdot \mathscr{A} )\\
  & +  \Pi_k(|D|^{-1}\lambda,   \delta\lambda\cdot \mathscr{A}+\lambda\cdot \delta \mathscr{A} )
  \end{aligned}
\end{equation}

\subsubsection{Equations for $(\lambda_k,\lambda_{k,\gamma})$ and $(\delta\lambda_k,\delta\lambda_{k,\gamma})$}
We consider the variables $\varphi_k := Q_k\lambda$, where $Q_k = P_k$ or $P_{k,\gamma}$ are Littlewood-Paley operators. Set $U^\alpha = \boldsymbol{1}_{\{d=3\}}(2A-V)^\alpha$, we rewrite the quasilinear Schr\"odinger equation \eqref{eq_lam} schematically in the form
 \begin{equation}
    \begin{aligned}\label{eq_lamk}
 		&\sqrt{-1}( \partial_t + U^\alpha_{\le k} \pa_\alpha )\varphi_k +\partial_{\alpha}( g^{\alpha \beta} \partial_{\beta}  \varphi_k) \\
	= 	&   -[Q_k,g^{\alpha \beta} \partial_{\alpha \beta}^2 \lambda] + \pa_x h \cdot \pa_x \lambda_k \\
   & - \boldsymbol{1}_{\{d =3\}}( [Q_k,U^\alpha_{\le k}\pa_\alpha]\lambda + Q_k(\mathscr{A}_{\ge k}\cdot \pa_x \lambda)) + \boldsymbol{1}_{\{d=4\}}Q_k(\mathscr{A}\cdot \pa_x\lambda) \\
   &   + Q_k (  \pa_x \mathscr{A}\cdot \lambda + (\lambda,\mathscr{A})^2 \cdot \lambda) \\
		=& ~   L(2^k  \pa_x h, \lambda_k) + \boldsymbol{1}_{\{d=3\}} L(\pa_x \mathscr{A}_{\le k+2},\lambda_k) + \boldsymbol{1}_{\{d=4\}}L(2^k \mathscr{A}_{\le k+2}, \lambda_k ) \\
    &+ \Pi_k(\pa_x \mathscr{A},\lambda) + Q_k [(\partial_x \mathscr{A}+(\lambda,\mathscr{A})^2) \cdot \lambda]  \\
 		= &: \mathcal{N}(\phi_k).
 \end{aligned}
 \end{equation}
Similarly, the difference $\delta\varphi_k$ satisfies
\begin{equation}\label{exp_dlamk}
\begin{aligned}
&\sqrt{-1}  \partial_t   \delta\varphi_k +\partial_{\alpha}( g^{\alpha \beta} \partial_{\beta}  \delta\varphi_k) \\
		=& \, L(2^k \delta \mathscr{A}_{\le k} + 2^k \pa_x \delta h + 2^{2k}\delta h,\lambda_{k})    + \Pi_k(\delta \mathscr{A},\pa_x \lambda)    + \Pi_k(\pa_x  \mathscr{A},\delta\lambda)\\
    &+ L(2^k \pa_x h,\delta \lambda_k) +  \boldsymbol{1}_{\{d=3\}}   L(\pa_x \mathscr{A}_{\le k+2},\delta\lambda_k)  + \boldsymbol{1}_{\{d=4\}}  L(2^k \mathscr{A}_{\le k+2}, \delta\lambda_k )  \\
		& +Q_k [(\partial_x \mathscr{A},(\lambda,\mathscr{A})^2) \cdot \delta \lambda]   + Q_k [(\partial_x \delta \mathscr{A}+(\delta\lambda,\delta \mathscr{A})\cdot (\lambda,\mathscr{A}) ) \cdot \lambda]\\
 		= & :\, \mathcal{N} (\delta \phi_k),
\end{aligned}
\end{equation}
where the coefficient-difference terms involving $\delta h$ and $\delta\mathscr{A}$ carry one extra derivative, and are therefore responsible for the one-derivative loss in the linearized analysis and in the weak Lipschitz bound.

This completes the schematic form of the equations that will be used in the estimates below.

\section{Balance laws, div-curl lemma and interaction Morawetz estimates}
\subsection{Balance laws for Schr\"odinger variables}
In this subsection we establish the balance laws for $\phi_k \in \{\varphi_k, 2^{-k}\delta \varphi_k\}$.
\subsubsection{Energy balance law}

Multiplying~\eqref{eq_lamk} or \eqref{exp_dlamk} by $ \overline{\phi}_k $ and taking the imaginary part, we compute
\begin{align}
 &~\operatorname{Im}\bigl(	\sqrt{-1}( \partial_t+ U_{\leq k}^\alpha \partial_{\alpha})  \phi_k \cdot\overline{\phi}_k \bigr)\notag\\
= & ~ \operatorname{Re}\bigl(	( \partial_t+ U_{\leq k}^\alpha \partial_{\alpha})  \phi_k \cdot\overline{\phi}_k \bigr)\notag\\
 =&   ~\frac 12 \partial_t |\phi_k|^2+\frac 12\partial_{\alpha} (U_{\leq k}^\alpha   |\phi_k|^2)-\frac 12\partial_{\alpha}U_{\leq k}^\alpha   |\phi_k|^2,
 \end{align}
and
 \begin{equation}
   \operatorname{Im}\bigl( \partial_{\alpha}( g^{\alpha \beta} \partial_{\beta}  \phi_k)\cdot\overline{\phi}_k \bigr) = \operatorname{Im}  \partial_{\alpha}( \overline{\phi}_k    \partial^{\alpha}  \phi_k).
 \end{equation}
Combining these identities yields the mass balance law
\begin{equation}\label{energy}
     \begin{aligned}
 	&~\frac 12 \partial_t |\phi_k|^2+\partial_{\alpha} \bigl(\frac 12 U_{\leq k}^\alpha   |\phi_k|^2+\operatorname{Im}    (\overline{\phi}_k   \partial^{\alpha}  \phi_k)\bigr)  \\
  =&~ \operatorname{Im} ( \mathcal{N}(\phi_k)\overline{\phi}_k ) +\frac 12\partial_{\alpha}U_{\leq k}^\alpha   |\phi_k|^2=:\mathcal{E} (\phi_k ),
 \end{aligned}
\end{equation}
where $\mathcal{E}(\phi_k)$ denotes the nonlinear source terms of the mass balance law.

 \subsubsection{Momentum balance law}

Multiplying the equation by $\partial_\gamma \overline{\phi}_k$ and taking the real part, we obtain
 \begin{align}
 &~	\operatorname{Re}(\partial_\gamma \overline{\phi}_k\cdot 	\sqrt{-1}( \partial_t+ U_{\leq k}^\alpha \partial_{\alpha})  \phi_k )\notag\\
  =&\, \frac{1}{2} 	\operatorname{Re}(\partial_\gamma \overline{\phi}_k\cdot 	\sqrt{-1}( \partial_t+ U_{\leq k}^\alpha \partial_{\alpha})  \phi_k ) -\frac 12 	\operatorname{Re}(\partial_\gamma \phi_k\cdot 	\sqrt{-1}( \partial_t+ U_{\leq k}^\alpha \partial_{\alpha})  \overline{\phi}_k ) \notag\\
 	=&\,   \frac{1}{2} \partial_\gamma	\operatorname{Re}( \overline{\phi}_k\cdot 	\sqrt{-1}( \partial_t+ U_{\leq k}^\alpha \partial_{\alpha})  \phi_k ) -\frac 12\partial_t \operatorname{Re}( \sqrt{-1}\, \overline{\phi}_k\cdot 	 \partial_\gamma   \phi_k )\notag \\
 	&\,-\frac 12\partial_\alpha	\operatorname{Re}( \sqrt{-1}  U_{\leq k}^\alpha  \overline{\phi}_k  \partial_{\gamma} \phi_k )-\frac 12 \operatorname{Re}( \overline{\phi}_k  \sqrt{-1}   \partial_{\gamma}  U_{\leq k}^\alpha  \partial_{\alpha} \phi_k )\label{Mom_com_1}\\
 	&\,+\frac 12\operatorname{Re}( \sqrt{-1} \partial_\alpha	  U_{\leq k}^\alpha  \overline{\phi}_k  \partial_{\gamma} \phi_k ) \notag\\
 	=&\, \frac12 \partial_\gamma \operatorname{Re} (\overline{\phi}_k\mathcal{N}(\phi_k)-\overline{\phi}_k \partial_{\alpha}( g^{\alpha \beta} \partial_{\beta}  \phi_k))+\frac 12\partial_t \operatorname{Im} (\overline{\phi}_k \partial_\gamma \phi_k )+\frac 12 \partial_ \alpha \operatorname{Im} (    U_{\leq k}^\alpha  \overline{\phi}_k   \partial_{\gamma} \phi_k )\notag\\
 	&\,+\frac 12 \operatorname{Im} (   \overline{\phi}_k    \partial_\gamma U_{\leq k}^\alpha  \partial_{\alpha} \phi_k -   \overline{\phi}_k    \partial_ \alpha U_{\leq k}^\alpha  \partial_{\gamma} \phi_k ). \notag
 \end{align}
 and
\begin{align}
    \operatorname{Re}(\partial_\gamma \overline{\phi}_k \partial_\alpha (g^{\alpha \beta} \partial_\beta \phi_k))
  =& ~\partial_\alpha \operatorname{Re}(\partial_\gamma \overline{\phi}_k g^{\alpha \beta}  \partial_\beta \phi_k) - \frac{1}{2} \partial_\gamma (g^{\alpha \beta} \partial_\alpha \overline{\phi}_k \partial_\beta \phi_k   )  \notag\\
  & + \frac{1}{2} \partial_\gamma g^{\alpha\beta}\partial_\alpha \phi_k \partial_\beta \phi_k,
   \end{align}

For the term appearing in \eqref{Mom_com_1}, we note that
 \begin{equation}
   \begin{aligned}
 -\frac{1}{2} \partial_\gamma \operatorname{Re}(\overline{\phi}_k \partial_\alpha (g^{\alpha \beta} \partial_\beta \phi_k)) =&\, -\frac{1}{2} \partial_\gamma \partial_\alpha \operatorname{Re}(\overline{\phi}_k g^{\alpha \beta} \partial_\beta \phi_k) + \frac{1}{2} \partial_\gamma (g^{\alpha \beta} \partial_\alpha \overline{\phi}_k \partial_\beta \phi_k)\notag\\
 = &\,-\frac{1}{4} \partial_\gamma \partial^2_{\alpha \beta} (g^{\alpha \beta}  |\phi_k|^2) + \frac{1}{4} \partial_\gamma \partial_\alpha (|\phi_k|^2 \partial_\beta g^{\alpha \beta})\notag\\
 &+ \frac{1}{2} \partial_\gamma (g^{\alpha \beta} \partial_\alpha \overline{\phi}_k \partial_\beta \phi_k).
 \end{aligned}
 \end{equation}

Collecting the above identities, we obtain the momentum balance law for $\phi_k$:
\begin{equation}\label{momentum}
  \begin{aligned}
&\frac{1}{2} \partial_t \operatorname{Im}(\overline{\phi}_k \partial_\gamma \phi_k) + \partial_\alpha [\operatorname{Re}(\partial_\gamma \overline{\phi}_k g^{\alpha \beta} \partial_\beta \phi_k) + \frac{1}{2} \operatorname{Im}(U^\alpha_{\leq k} \overline{\phi}_k \partial_\gamma \phi_k)]\\
&+ \frac{1}{4} \partial_\gamma \partial_\alpha (|\phi_k|^2 \partial_\beta g^{\alpha\beta}) - \frac{1}{4} \partial_\gamma \partial^2_{\alpha \beta} (|\phi_k|^2 g^{\alpha\beta})  \\
=&\,\operatorname{Re}(\partial_\gamma \overline{\phi}_k \mathcal{N}(\phi_k)) + \frac{1}{2} \operatorname{Im}(\overline{\phi}_k \partial_\alpha U^\alpha_{\leq k} \partial_\gamma \phi_k) - \frac{1}{2} \operatorname{Im}(\overline{\phi}_k \partial_\gamma U^\alpha_{\leq k} \partial_\alpha \phi_k) \\
 & - \frac{1}{2} \partial_\gamma g^{\alpha\beta}\partial_\alpha \phi_k \partial_\beta \phi_k   -  \frac{1}{2} \partial_\gamma (\operatorname{Re}(\overline{\phi}_k \mathcal{N}(\phi_k)))\\
 =:&~ \mathcal{M}_\gamma(\phi_k).
\end{aligned}
\end{equation}
Or equivalently,
\begin{equation}\label{momentum2}
  \begin{aligned}
&  \partial_t \operatorname{Im}(\overline{\phi}_k \partial_\gamma \phi_k) + \partial_\alpha [ 2\operatorname{Re}(\partial_\gamma \overline{\phi}_k   \partial^\alpha \phi_k) -\frac{1}{2} \pa_\gamma \partial^\alpha |\phi_k|^2 +\operatorname{Im}(U^\alpha_{\leq k} \overline{\phi}_k \partial_\gamma \phi_k)] \\
=&~2 \mathcal{M}_\gamma(\phi_k),
\end{aligned}
\end{equation}
where we use the shorthand $\pa^\alpha = g^{\alpha \beta}\pa_\beta$.  Raising the index in \eqref{momentum2} with the help of \eqref{eq_g}, we obtain
\begin{equation}\label{momentum3}
  \begin{aligned}
&  \partial_t \operatorname{Im}(\overline{\phi}_k \partial^\gamma \phi_k) + \partial_\alpha [  \operatorname{Re}(\partial^\gamma \overline{\phi}_k   \partial^\alpha \phi_k) -\operatorname{Re}( \overline{\phi}_k   \partial^\gamma\partial^\alpha \phi_k) +\operatorname{Im}(U^\alpha_{\leq k} \overline{\phi}_k \partial^\gamma \phi_k) ]\\
=&~2 \mathcal{M}^\gamma(\phi_k) + L (2^{2k}\pa_x h, \phi_k, \bar{\phi}_k ) +  2^k  L\big((\pa_x \mathscr{A},\mathscr{A}^2,\lambda^2), \phi_k, \bar{\phi}_k \big) \\
=& ~2\widetilde{\mathcal{M}}^\gamma(\phi_k).
\end{aligned}
\end{equation}

 \subsection{Div-curl lemma and bilinear estimates}\label{sec:divcurl}
The following div-curl lemma, originally due to the third author \cite{Zhou23}, is the central analytic tool that converts the integrated balance laws~\eqref{energy}, \eqref{momentum2} into bilinear space-time estimates.

\begin{lemma}[Div-curl lemma]\label{lemdiv-curl}
	Suppose that the functions $f^{ij}$ ($i,j=1,2$) satisfy
	$$
	\begin{gathered}
		\left\{\begin{array}{l}
			\partial_t f^{11}+\partial_x f^{12}=G^1, \\
			\partial_t f^{21}-\partial_x f^{22}=G^2,\\
		\end{array}\right. \\
		f^{11}, f^{12}, f^{21}, f^{22} \rightarrow 0, \text{ as  }x \rightarrow \infty,
	\end{gathered}
	$$
	then, for every $s\in\mathbb{R}_+$, we have
	\begin{equation}\label{divcurl}
	\begin{aligned}
		&\int_s^{+\infty} \int_{-\infty}^{+\infty} f^{11}(t-s) f^{22}(t)+f^{12}(t-s) f^{21}(t)  \mathrm{ d} x \mathrm{ d} t \\
        \leq&~ 2\sup_{t\in\mathbb{R}_+} \|f^{11}(t) \|_{L^1_x}\cdot\sup_{t\in\mathbb{R}_+} \|f^{21}(t) \|_{L^1_x}\\
		&+\left|\int_s^{+\infty} \int_{-\infty}^{+\infty}   G^2(t,x)\int_{-\infty}^x f^{11} (t-s,y) \mathrm{ d} y  \mathrm{ d} x \mathrm{ d} t\right| \\
		&+\left|\int_s^{+\infty} \int_{-\infty}^{+\infty}  G^1(t-s,x) \int^{+\infty}_x f^{21} (t,y) \mathrm{ d} y  \mathrm{ d} x \mathrm{ d} t\right|,
	\end{aligned}
	\end{equation}
	provided the right-hand side is finite.
\end{lemma}
\begin{proof}
	 The same computation as in \cite{WaZ24} yields
	\begin{align*}
	&	\frac{\partial}{\partial t} \int_{x<y} f^{11} (t-s,x) f^{21} (t,y) \mathrm{ d} x\mathrm{ d} y+\int_{-\infty}^{+\infty}\left(f^{11}(t-s) f^{22}(t)+f^{12}(t-s) f^{21}(t)\right) \mathrm{ d} x \\
	=&~ \int_{-\infty}^{+\infty} \left( \int_{-\infty}^x f^{11} (t-s,y) \mathrm{ d} y\right)G^2(t,x) \mathrm{ d} x  +\int_{-\infty}^{+\infty} \left( \int^{+\infty}_x f^{21} (t,y) \mathrm{ d} y\right)G^1(t-s,x) \mathrm{ d} x,
	\end{align*}
    then \eqref{divcurl} follows by integrating the above inequality with respect to $t$ over $[s, +\infty)$.
	\end{proof}

To obtain the bilinear bounds~\eqref{BA4}--\eqref{BA5}, we pair the \emph{integrated} momentum law for $\phi_{k,\gamma}$ --- namely~\eqref{momentum4}, obtained by integrating the pointwise momentum balance law~\eqref{momentum2} over the transverse directions $\widehat{x}_\gamma$ --- with various integrated balance laws for the low-frequency variables.  The three lemmas below treat the three types of low-frequency counterparts:
\begin{enumerate}
  \item[\textbf{(i)}] a free heat flow $e^{t\Delta}v_j$ (Lemma~\ref{lem:heat_divcurl});
  \item[\textbf{(ii)}] the parabolic variables $\mathscr{A}_j$ (Lemma~\ref{lem:Ah_divcurl});
  \item[\textbf{(iii)}] the Schr\"odinger variable $\phi_j$ itself (Lemma~\ref{lem:lambda_divcurl}).
\end{enumerate}
The first lemma is presented with complete details; the remaining two follow the same template and we only indicate the necessary modifications.

For the reader's convenience, we record here the integrated momentum law for $\phi_{k,\gamma} \in \{\varphi_k, 2^{-k}\delta \varphi_k\}$, obtained by integrating~\eqref{momentum2} over $\widehat{x}_\gamma$ and translating by $y$:
\begin{equation}\label{momentum4}
  \begin{aligned}
&~ \partial_t \int \operatorname{Im}(\overline{\phi}_{k,\gamma} \partial_\gamma \phi_{k,\gamma})^{y }\, \mathrm{d} \hat{x}_\gamma\\
& + \partial_\gamma \int \delta_{\alpha \gamma} [ 2\operatorname{Re}(\partial_\gamma \overline{\phi}_{k,\gamma}   \partial^\alpha \phi_{k,\gamma}) -\frac{1}{2} \pa_\gamma \partial^\alpha |\phi_{k,\gamma}|^2 +\operatorname{Im}(U^\alpha_{\leq k} \overline{\phi}_k \partial_\gamma \phi_k)  ]^{y }\, \mathrm{d} \hat{x}_\gamma \\
=&~\int  2 \mathcal{M}_\gamma(\phi_{k,\gamma})^{y } \, \mathrm{d}\hat{x}_\gamma,
\end{aligned}
\end{equation}
which serves as the common high-frequency counterpart in all three applications of the div-curl lemma below.

\begin{lemma}\label{lem:heat_divcurl}
  Given $v \in L^1_{x_\gamma}L^2_{\hat{x}_\gamma}$ and set $v_j := P_j v$.  For $j+50 \le k$ we have
  \begin{equation}
    \begin{aligned}
    &2^k \sup_{y}\sup_{s\in\mathbb{R}_+}\Big\|\boldsymbol{1}_{\{t\ge s\}}\|e^{(t-s)\Delta}v_j\|_{L^2_{\hat{x}_\gamma}}\|\phi^y_{k,\gamma}\|^2_{L^2_{\hat{x}_\gamma}}\Big\|_{L^1_{t,x_\gamma}}\\
  \lesssim  &~ \big(\|\phi_{k}\|_{L^\infty_t L^2_x}^2 + 2^{-k}\|\mathcal{M}_\gamma(\phi_{k,\gamma})\|_{L^1_{t,x}}\big)\|v_j\|_{L^1_{x_\gamma}L^2_{\hat{x}_\gamma}}.
    \end{aligned}
  \end{equation}
\end{lemma}
\begin{proof}
  We set $s = 0$; the general case follows by time translation.  Apply Lemma~\ref{lemdiv-curl} to the following pair of integrated laws:
  \begin{itemize}
    \item the integrated momentum law~\eqref{momentum4} for $\phi_{k,\gamma}$ (high-frequency);
    \item the (trivial) integrated balance law for the low-frequency heat flow:
    \begin{equation}
      \Big|\partial_t \|e^{t\Delta}v_j\|_{L^2_{\hat{x}_\gamma}}\Big|\le
    \| \Delta e^{t\Delta}v_j \|_{L^2_{\hat{x}_\gamma}}
    \end{equation}
    which follows by differentiating $\|e^{t\Delta}v_j\|_{L^2_{\hat{x}_\gamma}}$ in time and using $\partial_t e^{t\Delta}v_j=\Delta e^{t\Delta}v_j$.
  \end{itemize}
  To relate the target quantity to the div-curl output, we use the elementary   estimate
  \begin{equation}\label{elem_est}
	\begin{aligned}
	& \sup_{y} \left\| \|e^{t\Delta}v_j \|_{L^2_{\widehat{x}_\gamma}}  \|   \phi_{k,\gamma}^y\|_{L^2_{\widehat{x}_\gamma}}^2  \right\|_{L^2_{t,x_\gamma}}  \\
	\lesssim &~\sup_y  2^{-2k} \left\| \|e^{t\Delta}v_j \|_{L^2_{\widehat{x}_\gamma}}  \|  \pa_\gamma   \phi^y_{k,\gamma}\|^2_{L^2_{\widehat{x}_\gamma}}  \right\|_{L^2_{t,x_\gamma}}.
	\end{aligned}
\end{equation}
Indeed, on the cone $\{|\xi|\sim 2^k,\ |\xi|\lesssim|\xi_\gamma|\}$ we have $\phi_{k,\gamma}=2^{-k}K\ast\pa_\gamma\phi_{k,\gamma}$ with $\|K\|_{L^1_x}\lesssim1$ (cf. Section~3). This mechanism is used throughout all three lemmas below.

  Applying Lemma~\ref{lemdiv-curl} to the above pair of integrated laws and transferring the $\partial_\gamma$ derivatives from $\pa_\gamma\pa_\alpha|\phi_{k,\gamma}|^2$ onto the low-frequency heat flow via integration by parts where necessary, we obtain
  \begin{align}
    & \sup_{y }\iint  \|e^{t\Delta}v_j\|_{L^2_{\hat{x}_\gamma}}\cdot\|\pa_\gamma \phi^{y }_{k,\gamma}\|^2_{L^2_{\hat{x}_\gamma}}\ \mathrm{d}t \mathrm{d}x_\gamma \notag\\
    \lesssim&~  \sup_{y } \Big|\iint \|e^{t\Delta}v_j\|_{L^2_{\hat{x}_\gamma}} \Big(\int \operatorname{Re}(h \cdot \pa_x \phi_{k,\gamma}\cdot \pa_x\bar{\phi}_{k,\gamma})^{y }   \notag\\
    &\qquad \qquad \qquad\qquad \qquad\quad+ \operatorname{Im}(U_{\le k}\cdot \bar{\lambda}_{k,\gamma}\cdot \pa_x \lambda_{k,\gamma})^{y }   \mathrm{d} \hat{x}_{\gamma}\Big)\ \mathrm{d}t \mathrm{d}x_\gamma\Big|\notag\\
    &\quad + \Big|\iint \pa_\gamma\|e^{t\Delta}v_j\|_{L^2_{\hat{x}_\gamma}} \cdot \int  \operatorname{Re}(g \cdot\pa_x\phi_{k,\gamma} \cdot \bar{\phi}_{k,\gamma})^{y }   \mathrm{d} \hat{x}_{\gamma}\ \mathrm{d}t \mathrm{d}x_\gamma\Big|\notag\\
    &\quad + \bigl(2^k\|\phi_{k,\gamma}\|_{L^\infty_t L^2_x}^2 + \|\mathcal{M}_\gamma(\phi_{k,\gamma})\|_{L^1_{t,x}}\bigr)\notag\\
    &\qquad \cdot \bigl(\|e^{t\Delta}v_j\|_{L^\infty_t L^1_{x_\gamma}L^2_{\hat{x}_\gamma}} + \|e^{t\Delta }\Delta v_j\|_{L^1_t L^1_{x_\gamma}L^2_{\hat{x}_\gamma}}\bigr)\notag\\
    =: &\operatorname{Err}_1 +\operatorname{Err}_2 +\operatorname{Err}_3 .
  \end{align}
  This decomposition reflects the common architecture of all three proofs: $\operatorname{Err}_1$ and $\operatorname{Err}_2$ are lower-order geometric corrections involving $h$, which are absorbed using the frequency separation $|j-k|\ge 50$ together with the smallness of the metric and gauge coefficients; $\operatorname{Err}_3$ collects the contributions from the nonlinear source terms  via the $L^1_{t,x}$ norm.

  For $\operatorname{Err}_{1,2}$, we exploit the frequency gap $j +50\le k$ to obtain
  \begin{align}
    \operatorname{Err}_{1,2}
    &\lesssim \bigl(2^{j-k} + \|h\|_{L^\infty_{t,x}} + 2^{-k}\|U_{\le k}\|_{L^\infty_{t,x}}  \bigr)\cdot \text{LHS} \notag\\
    &\lesssim (2^{-50} + C\varepsilon \langle  C \varepsilon  \rangle)\cdot \text{LHS},
  \end{align}
  which is absorbed into the left-hand side by the smallness of $\varepsilon$.

  For $\operatorname{Err}_3$, the standard regularity estimate for the heat semigroup gives
  \begin{align}
    \operatorname{Err}_3 \lesssim \bigl(2^k\|\phi_{k}\|_{L^\infty_t L^2_x}^2 + \|\mathcal{M}_\gamma(\phi_{k,\gamma})\|_{L^1_{t,x}}\bigr)\cdot \|v_j\|_{L^1_{x_\gamma}L^2_{\hat{x}_\gamma}}.
  \end{align}
  Combining the above estimates yields the desired bound.
\end{proof}

\begin{lemma}\label{lem:Ah_divcurl}
	Under the bootstrap assumptions~\eqref{BA1}--\eqref{BA5}, for $k \ge j$ and $s\in\mathbb{R}_+$ we have
		\begin{align}
	  \sup_{y} &~ 2^{\frac{k}{2}}\left\|\boldsymbol{1}_{\{t\ge s\}} \| \mathscr{A}_j (t-s) \|_{L^2_{\widehat{x}_\gamma}}  \|   \phi_{k,\gamma}^y\|_{L^2_{\widehat{x}_\gamma}}\right\|_{L^2_{t,x_\gamma}} \notag\\
    &   +2^{\frac{k}{2}}\left\|\boldsymbol{1}_{\{t\ge s\}}\|\mathscr{A}_j(t) \|_{L^2_{\widehat{x}_\gamma}}  \|   \phi_{k,\gamma}^y(t-s)\|_{L^2_{\widehat{x}_\gamma}}   \right\|_{L^2_{t,x_\gamma}} \notag\\
       \lesssim&~    \big(\|\phi_{k }\|_{L^\infty_{t}L^2_{x}}^2+2^{-k}\| \mathcal{M}_\gamma (\phi_{k,\gamma})  \|_{L^1_{t,x}}\big)^{\frac{1}{2}}\cdot \langle C \varepsilon  \rangle 2^{-(\frac{d}{2}-1) j }c_j .
	\end{align}
\end{lemma}

\begin{proof}
  We set $s = 0$ without loss of generality.   We first construct the integrated balance laws for the low-frequency parabolic variables. Using \eqref{exp_scrA}, we obtain the integrated balance law for $\mathscr{A}_j$:
\begin{equation}\label{Ascr_div}
  \begin{aligned}
    \pa_t \int |\mathscr{A}_j|^2\, \mathrm{d} \hat{x}_\gamma  =  \int  2^{2j} L(\mathscr{A}_j,\mathscr{A}_j) + \mathscr{A}_j\cdot \mathcal{K} (\mathscr{A}_j)\, \mathrm{d} \hat{x}_\gamma .
  \end{aligned}
\end{equation}

  We now apply the same div-curl procedure as in Lemma~\ref{lem:heat_divcurl}, pairing the integrated momentum law~\eqref{momentum4} for $\phi_{k,\gamma}$ with \eqref{Ascr_div}.
  \begin{align}
  &\iint \|\mathscr{A}_j\|_{L^2_{\hat{x}_\gamma}}^2 \|\pa_\gamma \phi^y_{k,\gamma}\|^2_{L^2_{\hat{x}_\gamma}} \mathrm{d} t \mathrm{d} x_\gamma \notag\\
  \lesssim&~  \sup_{y } \Big|\iint\|\mathscr{A}_j\|_{L^2_{\hat{x}_\gamma}}^2\Big(\int \operatorname{Re}(h \cdot \pa_x \phi_{k,\gamma}\cdot \pa_x\bar{\phi}_{k,\gamma})^{y }  \notag\\
    &\qquad \qquad \qquad\qquad \qquad+ \operatorname{Im}(U_{\le k}\cdot \bar{\lambda}_{k,\gamma}\cdot \pa_x \lambda_{k,\gamma})^{y }   \mathrm{d} \hat{x}_{\gamma}\Big)\ \mathrm{d}t \mathrm{d}x_\gamma\Big|\notag\\
  &\quad + \Big|\iint \pa_\gamma\|\mathscr{A}_j\|_{L^2_{\hat{x}_\gamma}}^2 \cdot \int  \operatorname{Re}(g \cdot\pa_x\phi_{k,\gamma} \cdot \bar{\phi}_{k,\gamma})^{y }   \mathrm{d} \hat{x}_{\gamma}\ \mathrm{d}t \mathrm{d}x_\gamma\Big|\notag\\
  &\quad + (2^k\|\phi_{k,\gamma}\|_{L^\infty_t L^2_x}^2 + \|\mathcal{M}(\phi_{k,\gamma})\|_{L^1_{t,x}})  \cdot (\|P_j\mathscr{A}_0 \|_{L^2_x}^2   + \|\pa_x \mathscr{A}_j\|^2_{L^2_{t,x}} +   \|2^{-j}\mathcal{K} (\mathscr{A}_j )\|^2_{L^2_{t,x}}     )\notag\\
     = &\operatorname{Err}_1 +\operatorname{Err}_2 +\operatorname{Err}_3,
  \end{align}

  The error terms $\operatorname{Err}_1$ and $\operatorname{Err}_2$ are estimated exactly as in Lemma~\ref{lem:heat_divcurl} and are therefore omitted.  The only new feature is $\operatorname{Err}_3$, which now involves the nonlinear sources $\mathcal{K}(\mathscr{A}_j)$; these are controlled by the estimates of Proposition~\ref{prop:K}, which completes the proof.  We note that the proof of Proposition~\ref{prop:K} uses only the bootstrap assumptions and Lemmas~\ref{prim_bound}--\ref{delta_prim_bound}; it does not rely on the div-curl estimates of the present section, so the apparent forward reference creates no circularity.
\end{proof}

\begin{lemma}\label{lem:lambda_divcurl}
	Under the bootstrap assumptions~\eqref{BA1}--\eqref{BA5}, for $k \ge j+50$ we have
		\begin{align}
	  \sup_{y}  & ~2^{\frac{k}{2}}\left\|\boldsymbol{1}_{\{t \ge s\}} \|  \phi_j(t-s) \|_{L^2_{\widehat{x}_\gamma}}  \|   \phi_{k,\gamma}^y(t)\|_{L^2_{\widehat{x}_\gamma}}\right\|_{L^2_{t,x_\gamma}} \notag\\
    &   +2^{\frac{k}{2}}\left\|\boldsymbol{1}_{\{t \ge s\}}\|  \phi_j (t) \|_{L^2_{\widehat{x}_\gamma}}  \|   \phi_{k,\gamma}^y(t-s)\|_{L^2_{\widehat{x}_\gamma}}   \right\|_{L^2_{t,x_\gamma}} \notag\\
       \lesssim&~    \big(\| \phi_j \|_{L_t^\infty L_x^2}^2  +  \| \mathcal{E} (\phi_j)  \|_{L^1_{t,x}}\big)^{\frac{1}{2}} \big(\|\phi_{k,\gamma }\|_{L^\infty_t L^2_{ x}}^2+2^{-k}\| \mathcal{M}_\gamma (\phi_{k,\gamma})  \|_{L^1_{t,x}}\big)^{\frac{1}{2}} .
	\end{align}
\end{lemma}
\begin{proof}
  We set $s = 0$.  Apply Lemma~\ref{lemdiv-curl} to the following pair of integrated laws:
  \begin{itemize}
    \item the integrated mass law for $\phi_j$,
     \begin{equation}\label{int_energy}
      \partial_t \int \frac12 |\lambda_j|^2 \, \mathrm{d}\widehat{x}_\gamma + \partial_\gamma \int \bigl(\operatorname{Im}(\overline{\lambda}_j \partial^{\gamma} \lambda_j) + \tfrac12 U^\gamma_{\le j} |\lambda_j|^2\bigr) \, \mathrm{d}\widehat{x}_\gamma = \int \mathcal{E}(\lambda_j) \, \mathrm{d}\widehat{x}_\gamma,
    \end{equation}
    obtained by integrating the pointwise mass balance law~\eqref{energy} over $\widehat{x}_\gamma$;
    \item the integrated momentum law~\eqref{momentum4} for $\phi_{k,\gamma}$.
  \end{itemize}

  Using the  estimate~\eqref{elem_est} to convert $\|\phi^y_{k,\gamma}\|_{L^2_{\hat{x}_\gamma}}$ into $2^{-k}\|\partial_\gamma\phi^y_{k,\gamma}\|_{L^2_{\hat{x}_\gamma}}$, and applying Lemma~\ref{lemdiv-curl} to the above pair of integrated laws, we obtain the following decomposition (cf.~the proof of Lemma~\ref{lem:heat_divcurl}):
\begin{align}
	&\sup_y \iint \| \phi_j \|_{L_{\widehat{x}_\gamma}^2}^2 \| \partial_\gamma \phi_{k,\gamma}^y \|_{L_{\widehat{x}_\gamma}^2}^2 \ \mathrm{d} x_\gamma \mathrm{d}  t  \notag\\
  \lesssim&~  \sup_{y } \Big|\iint \| \phi_j \|_{L_{\widehat{x}_\gamma}^2}^2 \Big(\int \operatorname{Re}(h \cdot \pa_x \phi_{k,\gamma}\cdot \pa_x\bar{\phi}_{k,\gamma})^{y } \notag\\
    &\qquad \qquad \qquad\qquad \qquad+ \operatorname{Im}(U_{\le k}\cdot \bar{\lambda}_{k,\gamma}\cdot \pa_x \lambda_{k,\gamma})^{y }   \mathrm{d} \hat{x}_{\gamma}\Big)\ \mathrm{d}t \mathrm{d}x_\gamma\Big|\notag\\
  &\quad + \Big|\iint \pa_\gamma\| \phi_j \|_{L_{\widehat{x}_\gamma}^2}^2 \cdot \int  \operatorname{Re}(g \cdot\pa_x\phi_{k,\gamma} \cdot \bar{\phi}_{k,\gamma})^{y }   \mathrm{d} \hat{x}_{\gamma}\ \mathrm{d}t \mathrm{d}x_\gamma\Big|\notag\\
	&\quad + \left| \iint \int  \partial_x  \phi_j \cdot  \phi_j \, \mathrm{d} \widehat{x}_\gamma \int( \phi_{k,\gamma}  \cdot \partial_x \phi_{k,\gamma})^y \mathrm{d} \widehat{x}_\gamma  \ \mathrm{d} x_\gamma \mathrm{d}  t  \right| \notag \\
	&\quad + ( 2^k \|   \phi_{k,\gamma}  \|^2_{L_t^\infty L_x^2} +  \| \mathcal{M}_\gamma(\phi_{k,\gamma} )  \|_{L^1_{t,x}} )  ( \| \phi_j \|_{L_t^\infty L_x^2}^2 +  \| \mathcal{E} (\phi_j)  \|_{L^1_{t,x}} ) \notag \\
	&= \mathrm{ Err}_1 + \mathrm{ Err}_2 + \mathrm{ Err}_3 + \mathrm{ Err}_4 .
\end{align}
  Among these, the leading term is the left-hand side, which arises from the pairing of $\|\partial_\gamma\phi_{k,\gamma}\|^2_{L^2_{\hat{x}_\gamma}}$ (from the integrated momentum law) with $\frac12\|\phi_j\|^2_{L^2_{\hat{x}_\gamma}}$ (from the integrated mass law).  The error analysis follows the template of Lemma~\ref{lem:heat_divcurl}: $\operatorname{Err}_1$ and $\operatorname{Err}_2$ are the geometric corrections, absorbed as before using the frequency gap $k \ge j+50$; $\operatorname{Err}_3$ comes from the cross-coupling of the momentum fluxes $\operatorname{Im}(\bar{\phi}_j\partial^\alpha\phi_j)$ and $\operatorname{Im}(\bar{\phi}_{k,\gamma}\partial_\gamma\phi_{k,\gamma})$ and is controlled analogously; $\operatorname{Err}_4$ collects the nonlinear source terms and is of the required form.  This completes the proof.
\end{proof}

\subsection{Interaction Morawetz estimate}\label{sec:morawetz}

In this subsection we follow the strategy of Ifrim and Tataru \cite{ifrim_Global_2024} to prove the interaction Morawetz estimates \eqref{BA5} and \eqref{dBA5} for $\phi_k$.

\begin{proposition}\label{pro:Morawetz}
  Under the bootstrap assumptions \eqref{BA1}--\eqref{BA5}, we have
  \begin{equation}
    \begin{aligned}
     \|\phi_k\|_{L^4_{t,y_J}L^2_{z_J}}^4 \lesssim &~ 2^k  \|\phi_k\|_{L^\infty_t L^2_x}^2  ( \|\phi_k\|_{L^\infty_t L^2_x}^2 +\|2^{-k}\widetilde{\mathcal{M}}^\gamma(\phi_k),  \mathcal{E}(\phi_k) \|_{L^1_{t,x}} )\\
     & + 2^k \|\phi_k\|_{L^\infty_t L^2_x}^2 \|L(2^k \pa_x h_{\le  k+2 } + \pa_x U_{\le k}, \phi_k, \phi_k)\|_{L^1_{t,x}}.
    \end{aligned}
  \end{equation}
\end{proposition}
\begin{proof}
We fix an increasing $3$-tuple $J$ of indices in $\{1,\dots,d\}$ and write $x=(y_J,z_J)\in\mathbb{R}^3\times\mathbb{R}^{d-3}$.  Without loss of generality we take $J=(1,2,3)$ and write $x=(y,z)$ with $y\in\mathbb{R}^3$ and $z\in\mathbb{R}^{d-3}$; the argument for the other tuples is identical.  For brevity we omit the integration domain and the volume elements in what follows.

We define the \textit{Morawetz interaction potential}
\begin{align} \label{interapoten}
 I(\phi_k) = -\int a_\gamma(y-y') \, |\phi_k(y,z)|^2 \, \operatorname{Im} (\bar{\phi}_k \pa^\gamma \phi_k) (y',z') \, \mathrm{d} y \mathrm{d}y' \mathrm{d}z \mathrm{d}z'
\end{align}
where the weight functions are given by
\begin{equation}
  a(y) = |y|,\quad a_\gamma(y) = \pa_\gamma a(y)=\frac{y_\gamma}{|y|},\quad a_{\alpha \gamma}(y) = \pa_\alpha\pa_\gamma a(y)=\frac{\delta_{\alpha \gamma}}{|y|}-\frac{y_\alpha y_\gamma}{|y|^3},
\end{equation}
so that
\begin{equation}
  (\Delta a)(y)=\frac{2}{|y|},\qquad (\Delta^2 a)(y)=-8\pi \delta(y).
\end{equation}
We note that for $m\ge1$, convolution with the kernel $\pa^m a$ is a homogeneous Fourier multiplier, namely the composition of $|D_y|^{m-4}$ with Riesz-type operators $\mathscr{R}^m$.

Taking the time derivative of the interaction Morawetz potential and using the balance laws \eqref{energy} and \eqref{momentum3}, we obtain
\begin{align}
&     \frac{\mathrm{d}}{\mathrm{d} t}I(\phi_k)\notag\\
    =& \int a_{\alpha \gamma}(y-y')\, |\phi_k(y,z)|^2 \, [2\mathrm{Re} (\pa^\gamma \overline{\phi}_k \pa^\alpha \phi_k) - \frac{1}{2}\pa^\gamma \pa^\alpha|\phi_k|^2 ](y',z')\notag\\
    &  - \int a_{\alpha \gamma}(y-y')\, 2 \mathrm{Im}(\overline{\phi}_k \pa^\gamma \phi_k)(y,z) \, \mathrm{Im}(\overline{\phi}_k \pa^\alpha \phi_k)(y',z')   + \operatorname{Err}_1 + \operatorname{Err}_2\notag\\
  =: & ~ J(\phi_k) + \operatorname{Err}_1 + \operatorname{Err}_2,
\end{align}
where
\begin{align}
   \operatorname{Err}_1  = & ~ \int a_\gamma(y-y')\, \bigl[\mathcal{E}(\phi_k)(y,z)\,\mathrm{Im}(\overline{\phi}_k \pa^\gamma \phi_k)(y',z')- 2|\phi_k(y,z)|^2 \widetilde{\mathcal{M}}^\gamma(\phi_k)(y',z')\bigr] . \notag\\
   \operatorname{Err}_2 =&~  -\int a_{\alpha \gamma}(y-y')\,(U^\alpha_{\le k}(y',z')-U^\alpha_{\le k}(y,z))|\phi_k|^2 (y,z)\,\mathrm{Im}(\overline{\phi}_k \pa^\gamma \phi_k)(y',z')\notag\\
   & = -\boldsymbol{1}_{\{d=3\}} \int_0^1\int (y-y') a_{\alpha \gamma}(y-y')\,|\phi_k|^2(y)\,\mathrm{Im}(\overline{\phi}_k \pa^\gamma \phi_k)(y')\notag\\
   &\qquad\qquad \qquad\qquad\cdot \pa_x U^\alpha_{\le k}\bigl(y'+ s(y-y')\bigr)\,\mathrm{d}s. \notag
\end{align}
In the second equality we used that $U_{\le k} = \boldsymbol{1}_{\{d=3\}}(2A-V)_{\le k}$, so this contribution is present only in the endpoint case $d=3$.  There the transverse variable $z$ is absent, and the difference $U^\alpha_{\le k}(y')-U^\alpha_{\le k}(y)$ can be expanded by the mean value theorem.

We now simplify the main term $J(\phi_k)$. Integrating by parts in the first integral to move one derivative off $|\phi_k|^2$, we find
\begin{align}
    &~\int a_{\alpha \gamma}(y-y')\, |\phi_k(y,z)|^2 \, [2\mathrm{Re} (\pa^\gamma \overline{\phi}_k \pa^\alpha \phi_k) - \frac{1}{2}\pa^\gamma \pa^\alpha|\phi_k|^2 ](y',z')\\
    =& ~ \int a_{\alpha \gamma}(y-y')\, \bigl[|\phi_k(y,z)|^2 \, \mathrm{Re} (2\pa^\gamma \overline{\phi}_k \pa^\alpha \phi_k)(y',z') + \frac{1}{2}\pa^\gamma|\phi_k(y,z)|^2\, \pa^\alpha|\phi_k(y',z')|^2\bigr]\notag\\
    &  + \operatorname{Err}_3,\notag
\end{align}
where the boundary term $\operatorname{Err}_3$, arising from the commutator between the covariant and Euclidean derivatives, is
\begin{equation}
  \begin{aligned}
    \operatorname{Err}_3
    =& \int (\pa^2 a)(y-y')\bigl[ (\pa_x h \cdot \pa_x|\phi_k|^2)(y,z)\, |\phi_k|^2 (y',z') + (\pa_x h \cdot |\phi_k|^2)(y,z)\, \pa_x|\phi_k(y',z')|^2 \bigr].
  \end{aligned}
\end{equation}
Substituting this back, $J(\phi_k)$ becomes
 \begin{align}
    J(\phi_k)
  =&~ 2\int a_{\alpha \gamma}(y-y')\, |\phi_k(y,z)|^2 \, \mathrm{Re} (\pa^\gamma \overline{\phi}_k \pa^\alpha \phi_k)(y',z')\notag\\
  &- 2\int a_{\alpha \gamma}(y-y')\, \mathrm{Im}(\overline{\phi}_k \pa^\gamma \phi_k)(y,z) \, \mathrm{Im}(\overline{\phi}_k \pa^\alpha \phi_k)(y',z')\notag\\
   &+ \frac{1}{2}\int a_{\alpha \gamma}(y-y')\,\pa^\gamma|\phi_k(y,z)|^2\,\pa^\alpha|\phi_k(y',z')|^2 + \operatorname{Err}_3. \label{eq:J_expanded}
 \end{align}

To reveal the coercive structure of $J(\phi_k)$, we employ the following algebraic identity, valid for any complex-valued function $u$:
 \begin{align}
    &\mathrm{Re}\bigl[( u(y,z) \pa^\gamma u(y',z')-\pa^\gamma u(y,z) u(y',z'))\cdot ( \bar{u}(y,z) \pa^\alpha \bar{u}(y',z')-\pa^\alpha \bar{u}(y,z) \bar{u}(y',z'))\bigr]\notag\\
    =& ~|u(y,z)|^2\mathrm{Re}(\pa^\gamma \bar{u} \pa^\alpha u)(y',z') +|u(y',z')|^2\mathrm{Re}(\pa^\gamma \bar{u} \pa^\alpha u)(y,z)\notag\\
    &  -2 \mathrm{Re}\bigl(\pa^\gamma u(y,z) u(y',z') \cdot \bar{u}(y,z) \pa^\alpha \bar{u}(y',z') \bigr)\notag\\
    =&~ |u(y,z)|^2\mathrm{Re}(\pa^\gamma \bar{u} \pa^\alpha u)(y',z') +|u(y',z')|^2\mathrm{Re}(\pa^\gamma \bar{u} \pa^\alpha u)(y,z) \notag\\
    &  -2\,\mathrm{Im}(\bar{u} \pa^\gamma u)(y,z)\, \mathrm{Im}(\bar{u} \pa^\alpha u)(y',z') - \frac{1}{2}\pa^\gamma|u(y,z)|^2 \pa^\alpha|u(y',z')|^2. \notag
 \end{align}
Thanks to the symmetry $a_{\alpha \gamma}=a_{\gamma \alpha}$ and the convexity of $a(y)=|y|$, the Hessian $a_{\alpha \gamma}$ is positive semi-definite.  Applying the above identity with $u=\phi_k$ and exploiting this positivity, we deduce from \eqref{eq:J_expanded} that
\begin{align}
    J(\phi_k)
   =&~ \int a_{\alpha \gamma}(y-y') \, F^\alpha(y,z;y',z') \overline{F^\gamma(y,z;y',z')} \notag\\
   &+ \int a_{\alpha \gamma}(y-y') \, \pa^\gamma|\phi_k(y,z)|^2 \, \pa^\alpha|\phi_k(y',z')|^2 + \operatorname{Err}_3\notag\\
   \ge&~ \int a_{\alpha \gamma}(y-y') \, \pa^\gamma|\phi_k(y,z)|^2 \, \pa^\alpha|\phi_k(y',z')|^2 + \operatorname{Err}_3,
\end{align}
where $F^\alpha(y,z;y',z') := \phi_k(y,z) \pa^\alpha \phi_k(y',z') - \pa^\alpha \phi_k(y,z) \phi_k(y',z')$, and we have discarded the nonnegative term $\int a_{\alpha \gamma}F^\alpha\bar{F}^\gamma$.

The dominant term is treated by a further integration by parts.  Recalling that $a_{\alpha \gamma} = \pa_\alpha\pa_\gamma a$ is the Hessian of $a(y) = |y|$ and that $(\Delta^2 a)(y)=-8\pi\delta(y)$, we replace the covariant derivatives $\pa^\gamma = g^{\gamma\beta}\pa_\beta$, $\pa^\alpha = g^{\alpha\beta}\pa_\beta$ by the Euclidean derivatives (equivalently, $g^{\alpha\beta}$ by $\delta^{\alpha\beta}$) and integrate by parts; the principal part yields
\begin{align}
     &\int a_{\alpha \gamma}(y-y') \, \pa^\gamma|\phi_k(y,z)|^2 \, \pa^\alpha|\phi_k(y',z')|^2 \notag\\
     =& ~8\pi   \|  \phi_k  \|_{L^4_yL^2_z}^4
        + \int (\pa^2 a)(y-y')\, (\pa_x h \cdot |\phi_k|^2)(y,z)\, \pa_x|\phi_k(y',z')|^2  \notag\\
     &  + \int (\pa^3 a)(y-y')\, (\pa_x h \cdot |\phi_k|^2)(y,z)\, |\phi_k(y',z')|^2  \notag\\
     &   + \int (\pa^4 a)(y-y')\, (h \cdot |\phi_k|^2)(y,z)\, |\phi_k(y',z')|^2 \notag\\
  = &~ 8\pi  \|\phi_k \|_{L^4_yL^2_z}^4 + \operatorname{Err}_4.
\end{align}
The first term on the right is precisely the desired coercive quantity; the remaining integrals arise from replacing the covariant derivatives $\pa^\alpha = g^{\alpha\beta}\pa_\beta$ by the Euclidean derivatives $\pa_\alpha$ (equivalently, $g^{\alpha\beta}$ by $\delta^{\alpha\beta}$), and are collected in $\operatorname{Err}_4$.

Integrating the above identities in time and rearranging, we arrive at the key estimate
\begin{equation}
    \begin{aligned}
     \|\phi_k \|_{L^4_{t,y_J}L^2_{z_J}}^4 \lesssim 2^{k}\|\phi_k\|_{L^\infty_t L^2_x}^4 +  \sum_{i=1}^4  \int |\operatorname{Err}_i| \, \mathrm{d}t.
    \end{aligned}
\end{equation}

It remains to bound each error term.

We first consider $\operatorname{Err}_1$ and $\operatorname{Err}_2$, which involve the nonlinear source terms $\mathcal{E}(\phi_k)$ and $\widetilde{\mathcal{M}}^\gamma(\phi_k)$ of the balance laws.  Directly from the definitions,
\begin{equation}
  \begin{aligned}
 \int|\operatorname{Err}_{1,2}|\, \mathrm{d}t
\lesssim &~ 2^k  \|\phi_k\|_{L^\infty_t L^2_x}^2   \|2^{-k}\widetilde{\mathcal{M}}^\gamma(\phi_k),  \mathcal{E}(\phi_k) \|_{L^1_{t,x}} \\
     & + 2^k \|\phi_k\|_{L^\infty_t L^2_x}^2 \|L(\pa_x U_{\le k}, \phi_k, \phi_k)\|_{L^1_{t,x}}.
\end{aligned}
\end{equation}

Next, we estimate $\operatorname{Err}_3$ and $\operatorname{Err}_4$, which involve higher derivatives of the weight $a$ coupled with the metric perturbation $h$.  Interpreting the convolution with the kernels $\pa^m a$ as homogeneous Fourier multipliers and transferring derivatives by Plancherel's theorem, we obtain
\begin{align}
   &~\int|\operatorname{Err}_{3,4}|\,\mathrm{d}t\notag\\
  \sim &~ \int_0^\infty   \Big| \int |D_y|^{-\frac{1}{2}}\mathscr{R}_y(|\phi_k|^2)(y,z) \cdot |D_{y}|^{-\frac{3}{2}}\mathscr{R}_y(\pa_x h_{\le k+ 2}\cdot \pa_x|\phi_k|^2 )(y,z') \Big|\mathrm{d}t \notag\\
  & + \int_0^\infty  \Big| \int |D_y|^{-\frac{1}{2}}\mathscr{R}_y(|\phi_k|^2)(y,z) \cdot |D_{y}|^{-\frac{1}{2}}\mathscr{R}_y(\pa_x h_{\le k+ 2} \cdot  |\phi_k|^2 )(y,z') \Big|\mathrm{d}t \notag\\
  & + \int_0^\infty\Big| \int \mathscr{R}_y^4  (h\cdot|\phi_k|^2)(y,z) \cdot  (|\phi_k|^2 )(y,z') \Big| \mathrm{d}t   \notag\\
  \lesssim &~ 2^k \|\phi_k\|_{L^\infty_t L^2_x}^2 \|L(2^k \pa_x h_{\le k+2}, \phi_k,\phi_k) \|_{L^1_{t,x}} + \|h\|_{L^\infty_{t,x}} \|\phi_k\|^4_{L^4_{t,y}L^2_z},
\end{align}
which is acceptable in view of the smallness assumption \eqref{BA1} on $h$.
This completes the proof.
\end{proof}

\section{Nonlinear estimates}\label{sec:nonlin}
\paragraph{Logical order of the estimates.}
We emphasize that the nonlinear estimates in this section are proved in acyclic order. Lemmas~\ref{W_bound}--\ref{delta_prim_bound} and Proposition~\ref{prop:K} use only the bootstrap assumptions \eqref{BA1}--\eqref{BA5} and \eqref{dBA1}--\eqref{dBA5}, together with the elementary Littlewood--Paley bounds of Section~3; in particular, Proposition~\ref{prop:K} does not use the div-curl estimates of Section~6, so its use in the proof of Lemma~\ref{lem:Ah_divcurl} is not circular. Lemma~\ref{flux} is proved next from the bootstrap assumptions and the curvature identity of Lemma~\ref{hidden_1}; Lemma~\ref{flux_delta} then uses, in addition, the div-curl estimate of Lemma~\ref{lem:heat_divcurl}. The estimate of the term \eqref{exp_Acf_3_1} in the proof of Lemma~\ref{flux_delta} uses only Lemma~\ref{lem:heat_divcurl}, the estimates already established in Lemmas~\ref{prim_bound}--\ref{flux}, and the bootstrap assumptions; it does not invoke Proposition~\ref{prop:source}. Finally, Proposition~\ref{prop:source} is deduced from Lemmas~\ref{flux} and~\ref{flux_delta}, and is used only in Section~\ref{sec:proof} to close the bootstrap.

We first prove the following estimates involving $W$ and $\delta W$, which show that $W$ obeys the same bounds as $(\pa_x h, A)$: compare \eqref{L2W} and \eqref{laW} with \eqref{BA2} and \eqref{BA3}, and likewise the difference bounds \eqref{L2dW}, \eqref{dlaW} and \eqref{ladW} with \eqref{dBA2}, \eqref{dBA3b} and \eqref{dBA3a}.

In the estimates below we shall absorb powers of $C$ without further mention: since $1\ll C\le\varepsilon^{-1}$ and $\varepsilon\ll1$, such factors are harmless.

We begin with the estimates for the correction $W$ and their difference analogues.
\begin{lemma}\label{W_bound}
  	Under the bootstrap assumptions \eqref{BA1}-\eqref{BA5} and \eqref{dBA1}--\eqref{dBA5}, for $k,j \in \mathbb{Z}$, we have the following bounds
  \begin{align}
        & \| \pa_x W_{k}\|_{L^2_{t,x}}\lesssim C 2^{-(\frac{d}{2}-1+\sigma)k}c_k(\sigma),\label{L2W}\\
        &\|\delta W_{k}\|_{L^2_{t,x}}\lesssim C 2^{-(\frac{d}{2}-1)k}d_k, \label{L2dW}\\
        &\mathfrak{B}(j,k; W,\lambda)\lesssim C^2 2^{-(\frac{d}{2}-1)j} 2^{-(\frac{d-1}{2}+\sigma)k} c_j c_k(\sigma), ~~~\text{ for } k\geq j,\label{laW}\\
        &\mathfrak{B}(j,k; W,\delta\lambda)\lesssim C^2 2^{-(\frac{d}{2}-1)j} 2^{-\frac{d-3}{2}k} c_j d_k, ~~~\text{ for } k\geq j,\label{dlaW}\\
        &\mathfrak{B}(j,k; \delta W,\lambda)\lesssim C^2 2^{-(\frac{d}{2}-2)j} 2^{-\frac{d-1}{2}k} d_j c_k, ~~~\text{ for } k\geq j.\label{ladW}
  \end{align}
\end{lemma}
\begin{proof}
The bound \eqref{L2W} follows directly from \eqref{BA5}, while \eqref{L2dW} follows in the same way from \eqref{dBA5}.

For the bound \eqref{laW} we may set $s=0$ without loss of generality, and it suffices to consider $j\le k-50$, since the complementary case follows from \eqref{BA1} and \eqref{L2W}.  We estimate the first term in the definition \eqref{def_B} of $\mathfrak{B}$; the second term is analogous.  As in Subsection~\ref{sec:morawetz}, we fix an increasing $3$-tuple $J$ of indices in $\{1,\ldots,d\}$ with $\gamma\in J$ and write $x=(y_J,z_J)\in\mathbb{R}^3\times\mathbb{R}^{d-3}$; we apply \eqref{BA4} and \eqref{BA5} to obtain
\begin{align}
  &~\sup_{y}\big\|\|W_j\|_{L^2_{\hat{x}}}\|\lambda^y_{k,\gamma}\|_{L^2_{\hat{x}}}\big\|_{L^2_{t,x_\gamma}}\notag\\
  \lesssim &~ \sum_{j\le j'\le k-50}2^{\frac{d-1}{2}j-j'}\sup_{y}\big\|\|\lambda_{j'}^2\|_{L^1_{\hat{x}}}\|\lambda^y_{k,\gamma}\|_{L^2_{\hat{x}}}\big\|_{L^2_{t,x_\gamma}}\notag\\
  & + \sum_{j'\ge k-50 } 2^{-j'}\|P_j L(\lambda_{j'},\lambda_{j'})\|_{L^2_{t,\hat{x}_\gamma}L^\infty_{x_\gamma} }\|\lambda_k\|_{L^\infty_t L^2_{x}}\notag\\
  \lesssim &~ \sum_{j\le j'\le k-50}2^{\frac{d-1}{2}j-j'}\|\lambda_{j'}\|_{L^\infty_{t,x_\gamma}L^2_{\hat{x}_\gamma}}\sup_{y}\big\|\|\lambda_{j'}\|_{L^2_{\hat{x}}}\|\lambda^y_{k,\gamma}\|_{L^2_{\hat{x}}}\big\|_{L^2_{t,x_\gamma}}\notag\\
  & + \sum_{j'\ge k-50 } 2^{\frac{d-2}{2}j-j'}\| \lambda_{j'} \|^2_{L^4_{t,y_J}L^2_{z_J }}\|\lambda_k\|_{L^\infty_t L^2_{x}}\notag \\
  \lesssim &~ \sum_{j\le j'\le k-50} 2^{-\frac{d-1}{2}(j'-j)}Cc_{j'} \cdot C^2 2^{-(\frac{d}{2}-1)j'} 2^{-(\frac{d-1}{2}+\sigma)k} c_{j'} c_k(\sigma)\notag\\
  & + \sum_{j'\ge k-50 } 2^{\frac{d-2}{2}j-j'}\cdot C^2 2^{-(d-\frac{5}{2})j'}c_{j'}^2 \cdot C 2^{-(\frac{d}{2}-1+\sigma)k}c_k(\sigma)\notag \\
  \lesssim &~ C^2 2^{-(\frac{d}{2}-1)j} 2^{-(\frac{d-1}{2}+\sigma)k} c_j c_k(\sigma),
\end{align}
which proves \eqref{laW}.

The bound \eqref{dlaW} follows as in the proof of \eqref{laW}, with \eqref{dBA4b} in place of \eqref{BA4}.
The bound \eqref{ladW} is analogous and requires only a slight modification of the proof of \eqref{laW}.  Again it is enough to consider $j\le k-50$, since the complementary case follows from \eqref{dBA1} and \eqref{L2dW}; again we estimate the first term in \eqref{def_B} of $\mathfrak{B}$, the other being analogous, and we obtain
\begin{align}
  &~\sup_{y}\big\|\|\delta W_j\|_{L^2_{\hat{x}}}\|\lambda^y_{k,\gamma}\|_{L^2_{\hat{x}}}\big\|_{L^2_{t,x_\gamma}}\notag\\
  \lesssim &~ \sum_{j\le j'\le k-50}2^{\frac{d-1}{2}j-j'}\sup_{y}\big\|\|\lambda_{j'}\cdot \delta\lambda_{j'}\|_{L^1_{\hat{x}}}\|\lambda^y_{k,\gamma}\|_{L^2_{\hat{x}}}\big\|_{L^2_{t,x_\gamma}}\notag\\
  & + \sum_{j'\ge k-50 } 2^{-j'}\|P_j L(\delta\lambda_{j'},\lambda_{j'})\|_{L^2_{t,\hat{x}_\gamma}L^\infty_{x_\gamma} }\|\lambda_k\|_{L^\infty_t L^2_{x}}\notag\\
  \lesssim &~ \sum_{j\le j'\le k-50}2^{\frac{d-1}{2}j-j'}\|\delta\lambda_{j'}\|_{L^\infty_{t,x_\gamma}L^2_{\hat{x}_\gamma}}\sup_{y}\big\|\|\lambda_{j'}\|_{L^2_{\hat{x}}}\|\lambda^y_{k,\gamma}\|_{L^2_{\hat{x}}}\big\|_{L^2_{t,x_\gamma}}\notag\\
  & + \sum_{j'\ge k-50 } 2^{\frac{d-2}{2}j-j'}\| \delta\lambda_{j'} \|_{L^4_{t,y_J}L^2_{z_J }}\|  \lambda_{j'} \|_{L^4_{t,y_J}L^2_{z_J }}\|\lambda_k\|_{L^\infty_t L^2_{x}}\notag \\
  \lesssim &~ \sum_{j\le j'\le k-50} 2^{j-\frac{d-3}{2}(j'-j)}Cd_{j'} \cdot C^2 2^{-(\frac{d}{2}-1)j'} 2^{- \frac{d-1}{2} k} c_{j'} c_k \notag\\
  & + \sum_{j'\ge k-50 } 2^{\frac{d-2}{2}j }\cdot C^2 2^{-(d-\frac{5}{2})j'}c_{j'}d_{j'} \cdot C 2^{-(\frac{d}{2}-1 )k}c_k \notag \\
  \lesssim &~ C^2 2^{-(\frac{d}{2}-2)j} 2^{- \frac{d-1}{2} k} d_j c_k .
\end{align}

This completes the proof.
\end{proof}
Combining the estimates above, we obtain the following bounds.
\begin{lemma}\label{prim_bound}
	Under the bootstrap assumptions \eqref{BA1}-\eqref{BA5}, for $k,j \in \mathbb{Z}$, we have the following bounds
\begin{align}
    &\sup_{y}\|(\lambda_{k},\mathscr{A}_k)\cdot (\lambda_{j},\mathscr{A}_j)^y\|_{L^2_{t,x}}\lesssim C^2   2^{-\frac{1}{2}|k-j|}  2^{-(\frac{d}{2}-1+\sigma)k}c_{j}c_k(\sigma),\ \text{for}\ k\ge j, \label{lala}\\
    & \|P_k [ (\lambda,\mathscr{A})^2 ]\|_{L^2_{t,x}}\lesssim C^2  \varepsilon 2^{-(\frac{d}{2}-1+\sigma)k} c_k(\sigma), \label{lala2}\\
	&\|L(h_{\le k}, \mathscr{A}_k), L(\mathscr{A}_{\le k}, h_k)\|_{L^2_{t,x}}\lesssim  C \varepsilon 2^{-(\frac{d}{2} +\sigma)k}c_k(\sigma)\label{hkA} \\
	&  \|P_k(\mathscr{A},\pa_x \mathscr{A})\|_{L^2_{t,x}\cap\ 2^{-k}L^1_t L^2_x}\lesssim C^2   \varepsilon 2^{-(\frac{d}{2}-2+\sigma)k} c_k(\sigma), \label{AdA}\\
	&\|P_k(\mathscr{A}\cdot (\lambda,\mathscr{A})^2)\|_{L^2_{t,x}\cap\ 2^{-k}L^1_t L^2_x}\lesssim C^3   \varepsilon^2 2^{-(\frac{d}{2}-2+\sigma)k} c_k(\sigma), \label{cubic}\\
	&\|\Pi_k\big(\lambda,(\pa_x \mathscr{A},(\mathscr{A},\lambda)^2)\cdot \lambda\big)\|_{L^1_t L^2_x}\lesssim C^3 \varepsilon^2 2^{-(\frac{d}{2}-2 + \sigma)k}c_k(\sigma).\label{quintic}
\end{align}
 where the norm of the scaled space $cX$ is defined via $\|\cdot\|_{cX} := c^{-1}\|\cdot\|_X$ for $c>0$.
\end{lemma}
\begin{proof}
The bound \eqref{lala} for $\lambda_k\cdot\lambda_j^y$ follows from applying \eqref{BA4} when $j\le k-50$ and \eqref{BA5} when $k-50\le j\le k$, while the bounds for the remaining pairs follow from \eqref{BA2}, \eqref{BA3} and \eqref{laW}.  Summing the resulting bounds over the frequencies yields \eqref{lala2}.

	The bound \eqref{hkA} follows by taking the $L^\infty_{t,x}$ norm of the low-frequency factor:
	\begin{equation}
		\begin{aligned}
		  &~ \|L(h_{\le k}, \mathscr{A}_k),\ L(\mathscr{A}_{\le k}, h_k)\|_{L^2_{t,x}}\\
		  \lesssim&~ \|h_{\le k}\|_{L^\infty_{t,x}}\|\mathscr{A}_k\|_{L^2_{t,x}} + \|\mathscr{A}_{\le k}\|_{L^\infty_{t,x}}\|h_k\|_{L^2_{t,x}}\\
		   \lesssim&~ \sum_{k'\le k}Cc_{k'} (\|\mathscr{A}_k\|_{L^2_{t,x}} + 2^{k'}\|h_k\|_{L^2_{t,x}})\lesssim C \varepsilon 2^{-(\frac{d}{2} +\sigma)k}c_k(\sigma).\\
		\end{aligned}
	\end{equation}

For \eqref{AdA} it suffices to consider $\Pi_k(\mathscr{A},\pa_x\mathscr{A})$; the low-high interactions $L(\mathscr{A}_{\le k},\pa_x\mathscr{A}_k)$ and $L(\mathscr{A}_k,\pa_x\mathscr{A}_{\le k})$ are analogous.  We estimate
\begin{align}
	&\|\Pi_k(\mathscr{A}, \pa_x \mathscr{A} )\|_{L^2_{t,x}\cap\ 2^{-k}L^1_t L^2_x}
	\lesssim  2^{\frac{d}{2}k}  ( \|\Pi_k(\mathscr{A}, \pa_x \mathscr{A} )\|_{L^2_{t}L^1_x} + 2^k \|\Pi_k(\mathscr{A}, \pa_x \mathscr{A} )\|_{L^1_{t,x}})  \notag\\
	\lesssim& ~2^{\frac{d}{2}k}   \sum_{j\ge k}2^j(\|  \mathscr{A}_j\|_{L^\infty_t L^2_x} + 2^k\|\mathscr{A}_j\|_{L^2_{t,x}} )\|\mathscr{A}_j\|_{L^2_{t,x}}  \notag
	\lesssim  2^{\frac{d}{2}k}\sum_{j\ge k}C^2 2^{-(\sigma +d-2 )j}c_jc_j(\sigma) \\
	\lesssim&~  C^2 \varepsilon 2^{-(\frac{d}{2}-2 + \sigma)k}c_k(\sigma).
\end{align}

For \eqref{cubic} we estimate
\begin{align}
	&~\|P_k[(\lambda,\mathscr{A})^2\cdot \mathscr{A}]\|_{L^2_{t,x} }\notag\\
	\lesssim&~\|P_{\le k}((\lambda,\mathscr{A})^2)\|_{L^\infty_{t,x}}\| \mathscr{A}_k\|_{L^2_{t,x} } +   \|P_{k}((\lambda,\mathscr{A})^2)\|_{L^2_{t,x}}\| \mathscr{A}_{\le k}\|_{L^\infty_{t,x}}\notag\\
	&\quad + 2^{\frac{d}{2}k}\sum_{j\ge k}\|P_j((\lambda,\mathscr{A})^2)\|_{L^2_{t,x}}\| \mathscr{A}_j\|_{L^\infty_t L^2_x}\notag\\
	\lesssim&~ \sum_{k'\le k} C^2 \varepsilon 2^{2k'}c_{k'} \cdot C 2^{-(\frac{d}{2}+\sigma)k}c_k(\sigma) + \sum_{k'\le k}C^2 \varepsilon 2^{-(\frac{d}{2}-1+\sigma)k}c_k(\sigma) \cdot C 2^{k'}c_{k'} \notag\\
	&\quad +  2^{\frac{d}{2}k}\sum_{j\ge k} C^2 \varepsilon 2^{\frac{j}{2}}c_j \cdot C 2^{-(\frac{3}{2}+\sigma)j}c_j(\sigma)\notag\\
	\lesssim&~ C^3 \varepsilon^2 2^{-(\frac{d}{2}-2+\sigma)k}c_k(\sigma),
\end{align}
and
\begin{align}
	&~2^{ k}\|P_k[(\lambda,\mathscr{A})^2\cdot \mathscr{A}]\|_{L^1_t L^2_x}\notag\\
	\lesssim& ~2^k\|P_{\le k}((\lambda,\mathscr{A})^2)\|_{L^2_t L^\infty_{x}}\|\mathscr{A}_k\|_{L^2_{t,x}} +  2^k  \|P_{ k}((\lambda,\mathscr{A})^2)\|_{L^2_{t,x}}\|\mathscr{A}_{\le k}\|_{L^2_t L^\infty_{x}}\notag\\
	&  + 2^{(\frac{d}{2}+1)k}\sum_{j\ge k}\|P_j((\lambda,\mathscr{A})^2)\|_{L^2_{t,x}}\|\mathscr{A}_j\|_{L^2_{t,x}}\notag\\
	\lesssim&~ \sum_{k'\le k} C^2 \varepsilon 2^{k'}c_{k'} \cdot C 2^{-(\frac{d}{2}-1+\sigma)k}c_k(\sigma) + \sum_{k'\le k}C^2 \varepsilon 2^{-(\frac{d }{2}-2+\sigma)k}c_k(\sigma) \cdot C c_{k'}\notag\\
	& +  2^{(\frac{d}{2}+1)k}\sum_{j\ge k} C^2 \varepsilon 2^{-(\frac{d}{2}-1+\sigma)j}c_j \cdot C 2^{-\frac d2 j}c_j\notag\\
	\lesssim& ~ C^3 \varepsilon^2 2^{-(\frac{d}{2}-2+\sigma)k}c_k(\sigma).
\end{align}

For \eqref{quintic}, it suffices to consider the term $\Pi_k(\lambda,\pa_x\mathscr{A}\cdot\lambda)$, since $(\mathscr{A},\lambda)^2$ obeys the same $L^2_{t,x}$ estimate as $\pa_x\mathscr{A}$.  We estimate
\begin{align}
  &~\|\Pi_k( \lambda, \pa_x \mathscr{A} \cdot \lambda)\|_{L^1_t L^2_x}\notag\\
  \lesssim &~ \sum_{\substack{k'\ge k \\ k_1\le k'}} 2^{\frac{d}{2}k}\|L(\lambda_{k'},\lambda_{k_1})\|_{L^2_{t,x}}\|\pa_x \mathscr{A}_{\ge k'}\|_{L^2_{t,x}}\notag\\
  & + \sum_{\substack{k'\ge k}}2^{k } \|L(\lambda_{k'},\lambda_{|\cdot-k'|\le 2})\|_{L^2_{t,x}}\|\pa_x\mathscr{A}_{\le k'}\|_{L^2_t L_x^d}\notag\\
  \lesssim &~ \sum_{\substack{k'\ge k \\ k_1\le k'}} 2^{\frac{d}{2}k} \cdot C^2 2^{-(\frac{d}{2}-1+\sigma)k'}    c_{k_1}c_{k'}(\sigma)\cdot C2^{-(\frac{d}{2}-1)k'}c_{k'}\notag\\
  & + \sum_{\substack{k'\ge k}}2^{k } \cdot C^2 2^{-(\frac{d}{2}-1+\sigma)k'}    c_{k'}c_{k'}(\sigma)\cdot C\varepsilon\notag\\
  \lesssim &~ C^3 \varepsilon^2 2^{-(\frac{d}{2}-2+\sigma)k}c_k(\sigma).
\end{align}
This completes the proof.
\end{proof}

We also need the difference analogue of Lemma~\ref{prim_bound} for the linearized variables.
\begin{lemma}\label{delta_prim_bound}
Under the bootstrap assumptions \eqref{BA1}--\eqref{BA5} and \eqref{dBA1}--\eqref{dBA5}, for $k,j\in\mathbb{Z}$ we have
\begin{align}
 &\sup_y\|(\delta\lambda_k,\delta\mathscr{A}_k)\cdot\lambda_j^y\|_{L^2_{t,x}}
   \lesssim C^2 2^{-\frac12|k-j|}2^{k-(\frac d2-1)j}d_kc_j,\ \text{for } j\ge k, \label{dlala}\\
 &\sup_y\|(\lambda_k, \mathscr{A}_k)\cdot\delta\lambda_j^y\|_{L^2_{t,x}}
   \lesssim C^2 2^{-\frac12|k-j|}2^{-(\frac d2-2)j}c_kd_j,\ \text{for } j\ge k, \label{dlala2}\\
 &\|L(\delta h_{\le k}, \mathscr{A}_k),\ L( h_{\le k}, \delta\mathscr{A}_k)\|_{L^2_{t,x}}\lesssim  C^2  \varepsilon 2^{-(\frac{d}{2}-1)k}d_k, \label{dhkA}\\
   &\|L(\delta\lambda_{\le k},\lambda_k),\ L(\lambda_{\le k},\delta \lambda_k)\|_{L^\infty_t L^2_x \cap \ 2^{-k}L^2_{t,x}}\lesssim  C^2 \varepsilon 2^{-(\frac{d}{2}-3)k} d_k, \label{dlala3}\\
      &\| 2^k\Pi_k(|D|^{-1}\delta \lambda, \lambda)\|_{L^\infty_t L^2_x \cap \ 2^{-k}L^2_{t,x}}\lesssim  C^2 \varepsilon 2^{-(\frac{d}{2}-3)k} d_k,\label{dlala4}\\
 &\| \Pi_k(\delta\mathscr{A},2^k \mathscr{A} + \pa_x\mathscr{A})\|
   _{L^2_{t,x}}\lesssim C^2 \varepsilon2^{-(\frac d2-3)k}d_k, \label{dAdA}\\
 &\| \Pi_k(\delta\mathscr{A},2^k \mathscr{A} + \pa_x\mathscr{A})\|
   _{  2^{-k}L^1_tL^2_x}\lesssim C^2  2^{-(\frac d2-3)k}c_k d_k, \label{dAdA_2}\\
 &\|P_k[(\delta\lambda,\delta\mathscr{A})\cdot(\lambda,\mathscr{A})\cdot \mathscr{A}],\ \delta \mathscr{A}\cdot(\lambda,\mathscr{A})^2 \|
   _{L^2_{t,x} }\lesssim C^2\varepsilon 2^{-(\frac d2-3)k}d_k, \label{dcubic}\\
   &\|\Pi_k(|D|^{-1}\delta \lambda, \lambda \cdot \mathscr{A}), \Pi_k(|D|^{-1}  \lambda, \delta \lambda \cdot \mathscr{A}+\lambda \cdot \delta \mathscr{A})\|_{L^\infty_t L^2_x}\label{highcubic}  \lesssim C^3 \varepsilon^2 2^{-(\frac{d}{2}-3)k}d_k, \\
   &\|\Pi_k(|D|^{-1}\delta \lambda, \lambda \cdot \mathscr{A}), \Pi_k(|D|^{-1}  \lambda, \delta \lambda \cdot \mathscr{A}+\lambda \cdot \delta \mathscr{A})\|_{2^{-2k}L^1_t L^2_x}\label{highcubic_2}  \lesssim C^3 \varepsilon  2^{-(\frac{d}{2}-3)k}c_k d_k, \\
   &\|\Pi_k\big(|D|^{-1}\delta \lambda,  (\pa_x \mathscr{A},(\mathscr{A},\lambda)^2 )\cdot \lambda\big)\|_{L^1_t L^2_x}\label{highquntic}\\
   &+ \|\Pi_k\big(|D|^{-1}  \lambda,   \pa_x \delta\mathscr{A}\cdot \lambda + \pa_x \mathscr{A}\cdot \delta \lambda+ (\delta\mathscr{A},\delta\lambda) \cdot (\mathscr{A},\lambda)\cdot \lambda \big)\|_{L^1_t L^2_x} \notag\\
   &\qquad \lesssim  C^3 \varepsilon  2^{-(\frac{d}{2}-2)k}c_k d_k. \notag
\end{align}
\end{lemma}
\begin{proof}
Each bound is the difference analogue of the corresponding bound in Lemma~\ref{prim_bound}, and is proved in the same way, with the background bounds \eqref{BA1}--\eqref{BA5} supplemented by their difference analogues \eqref{dBA1}--\eqref{dBA5}. The strategy is as follows: if the difference variable is localized at frequency $\le k$, it suffices to spend one derivative, and the case is covered directly by the difference bounds. We therefore omit this case below and consider only the most difficult case, in which the difference variable is localized at frequency $>k$.

Specifically, the bounds \eqref{dlala} and \eqref{dlala2} follow from \eqref{dBA3a}, \eqref{dBA4a}, \eqref{dBA5} and \eqref{ladW}, together with \eqref{BA5} for $k\le j\le k+50$.  For \eqref{dlala2} the same argument applies with \eqref{dBA3b}, \eqref{dBA4b} and \eqref{dlaW} in place of \eqref{dBA3a}, \eqref{dBA4a} and \eqref{ladW}, respectively.  The bound \eqref{dhkA} follows from \eqref{dBA1} and \eqref{dBA2}, in the same way as \eqref{hkA}.  The remaining bounds are proved below.

For \eqref{dlala3}, it suffices to consider the most difficult case, in which the difference variable carries the high frequency, namely the term $L(\lambda_{\le k},\delta\lambda_k)$; the term in which it is localized at low frequency is covered by the general strategy above. Its $2^{-k}L^2_{t,x}$ norm is estimated as in the proof of \eqref{lala2}, with the difference analogues used for the factors involving $\delta\lambda$, while its $L^\infty_tL^2_x$ norm is obtained by taking the $L^\infty_{t,x}$ norm of the low-frequency factor, applying \eqref{dBA1} to the high-frequency factor, and summing over the frequencies.

For \eqref{dlala4}, we observe that $|D|^{-1}\delta\lambda$ behaves like $\lambda$, since the operator $|D|^{-1}$ gains one derivative. Its estimate is therefore analogous to the $L^2_{t,x}$ estimate of \eqref{lala2} and the $L^\infty_tL^2_x$ estimate of \eqref{dlala3}.

For \eqref{dAdA}, we argue as for \eqref{AdA}, with \eqref{dBA1} in place of the corresponding background bound; the estimate reduces to the summation
\begin{equation}
  2^{\frac d2 k}\sum_{j\ge k}C^2 2^{-(d-3)j}c_jd_j \lesssim C^2\varepsilon 2^{-(\frac d2-3)k}d_k,
\end{equation}
which is controlled by $d_k$ even in the endpoint case $d=3$.

For \eqref{dAdA_2}, the estimate is analogous to that of \eqref{dAdA}; since the $L^2_{t,x}$ norm can be used in full, \eqref{dBA2} gives a better bound, and the estimate is controlled by the summation
\begin{equation}
  2^{(\frac d2 + 1) k}\sum_{j\ge k}C^2 2^{-(d-2)j}c_jd_j \lesssim C^2  2^{-(\frac d2-3)k}c_k d_k,
\end{equation}

For \eqref{dcubic}, it suffices to consider the most difficult case, namely the high-high contribution to the first term $(\delta\lambda,\delta\mathscr{A})\cdot(\lambda,\mathscr{A})\cdot \mathscr{A}$ in which the difference variable carries the highest frequency; the second term and the remaining frequency orderings are analogous. In this case the contribution is
\begin{align}
	&\sum_{k_1\ge k_{2,3}\ge k-2}P_k((\delta\lambda,\delta\mathscr{A})_{k_1}\cdot(\lambda,\mathscr{A})_{k_2}\cdot \mathscr{A}_{k_3})\notag\\
	& =\sum_{\substack{k_1\ge k_{2,3}\ge k-2\\|k_1-\max\{k_2,k_3\}|\le 2}} P_k((\delta\lambda,\delta\mathscr{A})_{k_1}\cdot(\lambda,\mathscr{A})_{k_2}\cdot \mathscr{A}_{k_3})
\end{align}
We first consider the case $k_2\ge k_3$, for which we apply \eqref{dBA2}  and \eqref{dlala} to estimate
\begin{align}
    &~\sum_{\substack{|k_1-k_2|\le 2\\ k_1 \ge k_3\ge k-4}} \|P_k((\delta\lambda,\delta\mathscr{A})_{k_1}\cdot(\lambda,\mathscr{A})_{k_2}\cdot \mathscr{A}_{k_3}) \|_{L^2_{t,x} }\notag\\
	\lesssim&~\sum_{\substack{|k_1-k_2|\le 2\\ k_1 \ge k_3\ge k-4}}  2^{\frac{d}{2}k}\|P_k((\delta\lambda,\delta\mathscr{A})_{k_1}\cdot(\lambda,\mathscr{A})_{k_2}\cdot \mathscr{A}_{k_3}) \|_{L^2_{t}L^1_x }\notag\\
	\lesssim &~  \sum_{\substack{|k_1-k_2|\le 2\\ k_1 \ge k_3\ge k-4}}2^{\frac{d}{2}k}\|  (\delta\lambda,\delta\mathscr{A})_{k_1}\cdot \mathscr{A}_{k_3}\|_{L^2_{t,x}} \|(\lambda,\mathscr{A})_{k_2}\|_{L^\infty_tL^2_x } \notag\\
	\lesssim &~ \sum_{\substack{|k_1-k_2|\le 2\\ k_1 \ge k_3\ge k-4}} 2^{\frac{d}{2}k }  \cdot C^2  2^{-\frac{|k_1-k_3|}{2}} 2^{k_1-(\frac{d}{2}-1)k_3}d_{k_1}c_{k_3}\cdot C 2^{-(\frac{d}{2}-1)k_2}c_{k_2}  \notag\\
	\lesssim &~ \sum_{k_1 \ge k_3\ge k-4} C^3 2^{\frac{d}{2}k- (\frac{d}{2}-\frac{3}{2})(k_1+k_3)}c_{k_3} c_{k_1}d_{k_1}\notag\\
	\lesssim &~ C^3 \varepsilon^2 2^{-(\frac{d}{2}-3)k}d_k.
\end{align}
For the complementary case $k_2\le k_3$, we pair $k_1$ with $k_2$ and single out $k_3$, as above: applying \eqref{dlala2} to the pair $(\delta\lambda,\delta\mathscr{A})_{k_1}\cdot(\lambda,\mathscr{A})_{k_2}$ and the energy bound \eqref{BA1} to the remaining factor $\mathscr{A}_{k_3}$ yields
 \begin{align}
	\lesssim &~ \sum_{\substack{|k_1-k_3|\le 2 \\k_1 \ge k_2\ge k-4}}  2^{\frac{d}{2}k } \cdot C^2 2^{-(\frac{d}{2}-2)k_1}c_{k_2}d_{k_1}\cdot C 2^{-(\frac{d}{2}-1)k_3}c_{k_3}\notag\\
	\lesssim &~ \sum_{k_1 \ge k_2\ge k-4} C^3 2^{\frac{d}{2}k-(d-3)k_1}  c^2_{k_1}d_{k_1}\notag\\
	\lesssim &~ C^3 \varepsilon^2 2^{-(\frac{d}{2}-3)k}d_k.
 \end{align}

For \eqref{highcubic}, we first estimate the $L^\infty_tL^2_x$ norm. It suffices to consider the most difficult case, in which the difference variable in $\Pi_k(|D|^{-1}\lambda,\delta\lambda\cdot\mathscr{A})$ carries the highest frequency; the remaining cases are analogous. Arguing as in \eqref{dcubic}, we obtain
\begin{align}
	&~\sum_{\substack{ |k_2-k_3|\le 4\\ k_2\ge k_{1}\ge k-2}} \|\Pi_k(|D|^{-1}  \lambda_{k_1}, \delta \lambda_{k_2} \cdot \mathscr{A}_{k_3})\|_{L^\infty_t L^2_x}\notag\\
	\lesssim &~  2^{\frac{d}{2}k}\sum_{\substack{ |k_2-k_3|\le 4\\ k_2\ge k_{1}\ge k-2}}\|\Pi_k(|D|^{-1}  \lambda_{k_1}, \delta \lambda_{k_2} \cdot \mathscr{A}_{k_3})\|_{L^\infty_t L^1_x}\notag\\
	\lesssim &~  2^{\frac{d}{2}k}\sum_{\substack{ |k_2-k_3|\le 4\\ k_2\ge k_{1}\ge k-2}}\| |D|^{-1}  \lambda_{k_1}\|_{L^\infty_{t,x}}\|\delta \lambda_{k_2}\|_{L^\infty_t L^2_x}\| \mathscr{A}_{k_3}\|_{L^\infty_t L^2_x}\notag\\
	\lesssim &~ 2^{\frac{d}{2}k}\sum_{k_2\ge k_1 \ge k-2 } C c_{k_1}\cdot 2^{-(\frac{d}{2}-2)k_2}d_{k_2}\cdot C 2^{-(\frac{d}{2}-1)k_3}c_{k_3}\notag\\
	\lesssim &~ C^3 \varepsilon^2 2^{-(\frac{d}{2}-3)k} d_{k},
\end{align}

We next estimate \eqref{highcubic_2}.  For the first term $\Pi_k(|D|^{-1}\delta\lambda,\lambda\cdot\mathscr{A})$, the factor $|D|^{-1}\delta\lambda$ behaves like $\lambda$, as observed in \eqref{dlala4}, and the estimate is analogous to that of \eqref{cubic}.  It remains to consider the second term; as above, it suffices to treat the most difficult case, in which the difference variable in $\Pi_k(|D|^{-1}\lambda,\delta\lambda\cdot\mathscr{A})$ carries the highest frequency; the remaining cases are analogous.  We estimate
\begin{align}
	&~2^{2k}\sum_{\substack{|k_2-k_3|\le 4 \\k_2\ge k_{1}\ge k-2}} \|\Pi_k(|D|^{-1}  \lambda_{k_1}, \delta \lambda_{k_2} \cdot \mathscr{A}_{k_3})\|_{L^1_t L^2_x}\notag\\
	\lesssim &~  2^{(\frac{d}{2}+2)k}\sum_{\substack{|k_2-k_3|\le 4 \\k_2\ge k_{1}\ge k-2}}\|\Pi_k(|D|^{-1}  \lambda_{k_1}, \delta \lambda_{k_2} \cdot \mathscr{A}_{k_3})\|_{L^1_{t,x}}\notag\\
	\lesssim &~  2^{(\frac{d}{2}+2)k}\sum_{\substack{|k_2-k_3|\le 4 \\k_2\ge k_{1}\ge k-2}}2^{-k_1}\|  L( \lambda_{k_1}, \delta \lambda_{k_2})\|_{  L^2_{t,x}}\| \mathscr{A}_{k_3}\|_{L^2_{t,x}}\notag\\
	\lesssim &~ 2^{(\frac{d}{2}+2)k}\sum_{k_2\ge k_1 \ge k-2 } C^2 2^{-k_1-(\frac{d}{2}-2)k_2} c_{k_1}d_{k_2}\cdot C 2^{-\frac{d}{2}k_2}c_{k_2}\notag\\
	\lesssim &~ C^3 \varepsilon  2^{-(\frac{d}{2}-3)k} c_k d_{k}.
\end{align}

For \eqref{highquntic} we treat the two terms separately.  For the first term, $|D|^{-1}\delta\lambda$ behaves like $\lambda$ (see \eqref{dlala4}), so its estimate is analogous to that of \eqref{quintic}, with the difference bounds used for the factors involving $\delta\lambda$.  For the second term it suffices to consider the most difficult contribution, namely $\Pi_k(|D|^{-1}\lambda,\pa_x\delta\mathscr{A}\cdot\lambda)$, the contributions involving $\pa_x\mathscr{A}\cdot\delta\lambda$ being analogous and the quartic contribution being treated below.  Since $\pa_x\delta\mathscr{A}$ carries one derivative, it is enough to consider the case in which it is localized at frequency larger than $k'$; otherwise the inverse derivative $|D|^{-1}$ absorbs it.  We estimate
\begin{align}
  &~\sum_{\substack{|k_1-k_2|\le 4\\ k_1\ge k'\ge k}}\|P_k L\big(|D|^{-1}\lambda_{k'},P_{k'}(\pa_x\delta\mathscr{A}_{k_1} \cdot \lambda_{k_2} )\big)\|_{L^1_t L^2_x}\notag\\
  \lesssim &~ \sum_{\substack{|k_1-k_2|\le 4\\ k_1\ge k'\ge k}}2^{\frac{d}{2}k-k'}\|L(\lambda_{k'},\lambda_{k_1})\|_{L^2_{t,x}} \|\pa_x\delta\mathscr{A}_{k_2}\|_{L^2_{t,x}}\notag\\
  \lesssim &~ \sum_{\substack{|k_1-k_2|\le 4\\ k_1\ge k'\ge k}}2^{\frac{d}{2}k-k'} C^2 2^{-(\frac{d}{2}-1)k_2}c_{k'}c_{k_1}\cdot C 2^{-(\frac{d}{2}-2)k_2} d_{k_2}\notag\\
  \lesssim &~  \sum_{\substack{|k_1-k_2|\le 4\\ k_1\ge k'\ge k}}C^32^{\frac{d}{2}k-k'-(d-3)k_2}  c_{k'}c_{k_1}  d_{k_2}\notag\\
  \lesssim &~ C^3 \varepsilon  2^{-(\frac{d}{2}-2)k}c_k d_{k}.
\end{align}
Finally, we treat the quartic contribution, where $k_{\operatorname{med}}$ and $k_{\operatorname{max}}$ denote the median and the maximum of the triple $(k_1,k_2,k_3)$.  We estimate
\begin{align}
  &~\|\Pi_k\big(|D|^{-1}  \lambda,    (\delta\mathscr{A},\delta\lambda) \cdot  (\mathscr{A},\lambda) \cdot \lambda \big)\|_{L^1_t L^2_x}\notag\\
  \lesssim &~ \sum_{k'\ge k}\sum_{\substack{ \min_i k_i \ge k'\\ |k_{\operatorname{med}}-k_{\operatorname{max}}|\le 2 }}2^{\frac{d}{2}k-k'}\sup_{y,y'}\|\lambda_{k'}\cdot (\delta\mathscr{A},\delta \lambda)^{y}_{k_1}\|_{L^2_{t,x}} \|(\mathscr{A},\lambda)_{k_2} \cdot \lambda^{y'}_{k_3} \|_{L^2_{t,x}}
\end{align}
Here it suffices to consider the case in which $(\delta\mathscr{A},\delta\lambda)$ is localized at frequency larger than $k'$; otherwise the inverse derivative absorbs it.  We may then assume without loss of generality that $k_2=k_{\operatorname{med}}$, and we obtain
\begin{align}
    &~\sum_{k'\ge k}\sum_{\substack{k_1\ge k'\\ |k_{\operatorname{med}}-k_{\operatorname{max}}|\le 2 }}2^{\frac{d}{2}k-k'}\sup_{y,y'}\|\lambda_{k'}\cdot (\delta\mathscr{A},\delta \lambda)^{y}_{k_1}\|_{L^2_{t,x}} \|(\mathscr{A},\lambda)_{k_2} \cdot \lambda^{y'}_{k_3} \|_{L^2_{t,x}}\notag\\
    \lesssim &~ \sum_{k'\ge k}\sum_{\substack{k_1\ge k'\\ |k_{\operatorname{med}}-k_{\operatorname{max}}|\le 2 }}2^{\frac{d}{2}k-k'}\cdot C^2 2^{-(\frac{d}{2}-2)k_1}c_{k'}d_{k_1}\cdot C^2 2^{-(\frac{d}{2}-1)k_2}c_{k_2}c_{k_3} \notag\\
    \lesssim &~ \sum_{\substack{k'\ge k \\k_1\ge k',k_3}}C^4 2^{\frac{d}{2}k-k'-(d-3)k_1}c_{k'}c_{k_1}^2  d_{k_1}\notag\\
    \lesssim &~ C^4 \varepsilon^2 2^{-(\frac{d}{2}-2)k}c_k d_k.
\end{align}
This completes the proof.
\end{proof}

As a consequence of Lemmas~\ref{prim_bound} and~\ref{delta_prim_bound}, we obtain the following bounds for nonlinearities of $\mathscr{A}$ and $\delta \mathscr{A}$.
\begin{proposition}\label{prop:K}
	Under the bootstrap assumptions \eqref{BA1}--\eqref{BA5}, we have
	\begin{equation}\label{K_A}
	  \|\mathcal{K}(\mathscr{A}_k)\|_{L^2_{t,x}}\lesssim \langle C \varepsilon  \rangle C \varepsilon 2^{-(\frac{d}{2}-2+\sigma)k}c_k(\sigma).
	\end{equation}

	Moreover, under the additional assumptions \eqref{dBA1}--\eqref{dBA5}, we have
		\begin{align}
	  &\|\delta \mathcal{K}(h_k)\|_{L^2_{t,x}}\lesssim \langle C \varepsilon \rangle C \varepsilon 2^{-(\frac{d}{2}-2)k}d_k,\label{dK_h}\\
	  &\|\delta \mathcal{K}(A_{k,\mathrm{cf}})\|_{L^2_{t,x}}\lesssim \langle C \varepsilon \rangle C \varepsilon 2^{-(\frac{d}{2}-3)k}d_k,\label{dK_A_cf}\\
	  &\|\delta A_{k,\mathrm{df}}\|_{L^\infty_t L^2_x\cap\ 2^{-k}L^2_{t,x}}\lesssim \langle C \varepsilon \rangle C \varepsilon 2^{-(\frac{d}{2}-3)k}d_k.\label{dA_df}
	\end{align}
\end{proposition}
\begin{proof}
	By the definition of $\mathcal{K}(\mathscr{A}_k)$ in \eqref{exp_scrA}, the bound \eqref{K_A} follows from \eqref{hkA}, \eqref{AdA}, \eqref{lala2} and \eqref{cubic}.

	By the definition of $\delta \mathcal{K}(h_k)$ in \eqref{exp_dhk}, the bound \eqref{dK_h} follows from \eqref{dhkA}, \eqref{dlala3}, \eqref{dlala4}, \eqref{dAdA} and \eqref{highcubic}.

	By the definition of $\delta \mathcal{K}(A_{k,\mathrm{cf}})$ in \eqref{para_dA_cf}, the bound \eqref{dK_A_cf} follows from \eqref{dhkA}, \eqref{dAdA} and \eqref{dcubic}.

	Finally, using the relation \eqref{dAdf_lambda}, the bound \eqref{dA_df} follows from \eqref{dlala3}, \eqref{dlala4} and \eqref{highcubic}. This completes the proof.
\end{proof}

To estimate the Schr\"odinger nonlinearities we shall need the following $L^1_{t,x}$ estimates.

 \begin{lemma}\label{flux}
\begin{align}
	 &\| L(\lambda_k, P_k((\lambda,\mathscr{A})^2 \cdot\lambda )) \|_{L^1_{t,x}} \lesssim C^4 \varepsilon^2 (2^{-(\frac {d}{2} -1  + \sigma)k} c_k(\sigma))^2,   \label{Est_lamb_1}\\
	&\|\lambda_k \cdot \Pi_k(\pa_x \mathscr{A},\lambda ),\ \lambda_k \cdot L(2^k \mathscr{A}_{\ge k-50}, \lambda_k)\|_{L^1_{t,x}} \lesssim C^3 \varepsilon (2^{-(\frac {d}{2} -1  + \sigma)k} c_k(\sigma))^2,  \label{Est_lamb_2} \\
	&\|   L( 2^k (A_{\operatorname{df}},W)_{\le k-60}, \lambda_k,\lambda_k) \|_{L^1_{t,x}} \lesssim C^4 \varepsilon^2 (2^{-(\frac {d}{2} -1  + \sigma)k} c_k(\sigma))^2, \label{Est_lamb_3}\\
  &\|  L( 2^k \pa_x h_{\le k-60}, \lambda_k, \lambda_k) \|_{L^1_{t,x}} \lesssim C^4 \varepsilon^2 (2^{-(\frac {d}{2} -1  + \sigma)k} c_k(\sigma))^2, \label{Est_lamb_5}
\end{align}
In particular, in the three-dimensional case we have
\begin{equation}\label{Est_lamb_d3}
  \| L(\pa_x A_{\operatorname{cf},\le k-60}, \lambda_k,\lambda_k) \|_{L^1_{t,x}}\lesssim C^4 \varepsilon^2 (2^{-(\frac {d}{2} -1  + \sigma)k} c_k(\sigma))^2,
\end{equation}
while in the four-dimensional case we have
\begin{equation}\label{Est_lamb_d4}
  \|  L(2^k A_{\operatorname{cf},\le k-60}, \lambda_k,\lambda_k) \|_{L^1_{t,x}}\lesssim C^4 \varepsilon^2 (2^{-(\frac {d}{2} -1  + \sigma)k} c_k(\sigma))^2,
\end{equation}
\end{lemma}
\begin{proof}
	 The bound \eqref{Est_lamb_1} follows by applying \eqref{lala} to the individual frequency pairs and summing over all pairs; the off-diagonal factor $2^{-\frac12|j-k|}$ in \eqref{lala} makes the summation convergent.

	 The bound \eqref{Est_lamb_2} follows by applying \eqref{BA2} to the parabolic variable $\mathscr{A}$ and \eqref{lala} to the pair $(\lambda_k,\lambda_k^y)$.

  For the bound \eqref{Est_lamb_3} we use Lemma~\ref{hidden_1} to handle the term $\lambda \wedge \bar{\lambda}$.
  We introduce the symmetric subset $\mathcal{S}_k\subset\mathcal{S}$ by
\begin{equation}\label{S_kdef}
  \mathcal{S}_k := \{(k_1,k_2)\mid |k_{1,2}-k|\le 60,\ |k_1-k_2|\le 4\},
\end{equation}
Then, taking $\mathcal{S}=\mathcal{S}_k$ in Lemma~\ref{hidden_1}, we obtain
\begin{equation}\label{A_exp_0}
  \begin{aligned}
    &~\lambda\wedge \bar{\lambda} \\
    =&~  \Big(\sum_{k_2\le k_1-4}+\sum_{\substack{|k_1-k_2|\le 4 \\ |k_{1,2}-k|\ge 60 }}\Big)L(\lambda_{k_1},\lambda_{k_2}) \\
   & +  \sum_{(k_1,k_2)\in\mathcal{S}_k} \Bigl[
            \pa_x  L(|D|^{-1}\lambda_{k_1}\cdot  \lambda_{k_2})
           + L\big(|D|^{-1}\lambda_{k_1}, P_{k_2} ( \lambda \cdot\mathscr{A} ) \big)\Bigr],
  \end{aligned}
\end{equation}
Therefore, using the definition of $W$, we can write
 \begin{align}
 &~2^k(A_{\operatorname{df}}, W)_{\le k-60}\notag\\
  = &~2^k  P_{\le k-60}|D|^{-1}(\lambda\wedge \bar{\lambda})  +   2^k \Pi_{\le k-60}(|D|^{-1}\lambda,\lambda)\notag\\
   = &~ 2^k  \sum_{\substack{
  k_{1,2}\le k - 50  }} |D|^{-1}L(\lambda_{k_1},\lambda_{k_2}) +\sum_{\substack{|k_1-k_2|\le 4\\
  k_{1,2}\ge k + 50  }}  |D|^{-1}L(\lambda_{k_1},\lambda_{k_2})\label{A_exp_1} \\
   &  + \sum_{(k_1,k_2)\in\mathcal{S}_k} \Bigl[
              L( \lambda_{k_1}\cdot  \lambda_{k_2})
           + |D|^{-1} L\big( \lambda_{k_1}, P_{k_2} ( \lambda \cdot\mathscr{A} ) \big) \Bigr]\label{A_exp_2}.
\end{align}

For the contribution of the first term in \eqref{A_exp_1} we decompose $\lambda_k = \sum_{\gamma = 1}^d \lambda_{k,\gamma}$; it suffices to estimate
\begin{align}
 &\sum_{\substack{
  k_{1,2}\le k - 50  }}\sup_{y} \big\| \||D|^{-1}L(\lambda_{k_1},\lambda_{k_2})\|_{L^\infty_{\hat{x}_\gamma}}\|\lambda_{k,\gamma}\|^2_{L^2_{\hat{x}_\gamma}} \big\|_{L^1_{t,x_\gamma}}\notag\\
 \lesssim &~  \sum_{ k_{1,2}\le k - 50} 2^{\frac{d-2}{2}(k_1+k_2)}\big\| \| \lambda_{k_1} \|_{L^2_{\hat{x}_\gamma}} \| \lambda_{k_2} \|_{L^2_{\hat{x}_\gamma}}\|\lambda_{k,\gamma}\|^2_{L^2_{\hat{x}_\gamma}} \big\|_{L^1_{t,x_\gamma}}\notag\\
 \lesssim &~  \sum_{ k_{1,2}\le k - 50} C^4 c_{k_1}c_{k_2} (2^{-(\frac{d}{2}-1+\sigma)k}c_k(\sigma))^2 \notag\\
 \lesssim &~ C^4 \varepsilon^2 (2^{-(\frac{d}{2}-1+\sigma)k}c_k(\sigma))^2.
\end{align}
The contribution of the second term in \eqref{A_exp_1} is analogous, except that $\lambda_{k_{1,2}} = \sum_{\gamma = 1}^d\lambda_{k_{1,2},\gamma}$ is decomposed instead.

For the contribution of the first term in \eqref{A_exp_2}, using \eqref{BA5} directly and fixing an increasing $3$-tuple $J$ of indices in $\{1, \ldots,d\}$ as before, we obtain
\begin{align}
   &~ \sum_{|k_{1,2}-k|\le 60} \sup_{y }\| L(\lambda_{k_1},\lambda_{k_2}) \cdot|\lambda_{k}^y|^2\|_{L^1_{t,x}}\notag\\
\lesssim &~ \sum_{|k'-k|\le 60}  \|\lambda_{k'}\|_{L^4_{t,x}}^2 \|\lambda_k\|^2_{L^4_{t,x}}\notag\\
\lesssim &~ \sum_{|k'-k|\le 60} 2^{(d-3)k}  \|\lambda_{k'}\|_{L^4_{t,y_J}L^2_{z_J}}^2 \|\lambda_k\|^2_{L^4_{t,y_J}L^2_{z_J}}\notag\\
\lesssim &~ 2^{(d-3)k}  \cdot  (C 2^{-(\frac{d}{2}-\frac{5}{4})k} c_{k})^2 \cdot (C 2^{-(\frac{d}{2}-\frac{5}{4}+\sigma)k} c_{k}(\sigma))^2 \notag\\
\lesssim &~ C^4 \varepsilon^2 (2^{-(\frac{d}{2}-1+\sigma)k}c_k(\sigma))^2.\label{estia_1}
\end{align}

For the contribution of the second term in \eqref{A_exp_2}, we take the $L^1_t L^\infty_x$ norm and use the $L^1_t L^2_x$ estimate \eqref{cubic}.

    For the bound \eqref{Est_lamb_5} we again decompose $\lambda_k = \sum_{\gamma = 1}^d \lambda_{k,\gamma}$; it suffices to prove
\begin{equation}
  \sup_{y }\|2^k \pa_x h_{\le k-60}\cdot  |\lambda_{k,\gamma}^y|^2  \|_{L^1_{t,x}}\lesssim C^3 \varepsilon (2^{-(\frac{d}{2}-1+\sigma)k}c_k(\sigma))^2.
\end{equation}
To this end we use Lemma~\ref{hidden_curvature}.  We again use the notation $\mathcal{S}_k\subset\mathcal{S}$ defined in \eqref{S_kdef} together with the corresponding truncation $\mathcal{H}^k$ by
\begin{equation}
   \qquad \mathcal{H}^k := \sum_{(k_1,k_2)\in \mathcal{S}_k}\mathcal{H}[k_1,k_2],
\end{equation}
and decompose
\begin{equation}\label{decomp_h}
  h = h-\mathcal{H}^k + \mathcal{H}^k =: \tilde{h} + \mathcal{H}^k,
\end{equation}
The contribution of $\mathcal{H}^k$ can be controlled in much the same way as \eqref{estia_1}.
It remains to estimate the contribution of $\tilde{h}$, using Lemma~\ref{hidden_curvature} with $\mathcal{S}=\mathcal{S}_k$.  The function $\tilde{h}$ satisfies
\begin{equation}\label{tildeh_eq}
  (\pa_t -\pa^\gamma \pa_\gamma)\tilde{h} = \pa_x h\cdot \pa_x h +  \pa_x W + \operatorname{Ric} -(\pa_t -\pa^\gamma\pa_\gamma)  \mathcal{H}^k .
\end{equation}
By the cancellation identity of Lemma~\ref{hidden_curvature},
\begin{equation}
  \begin{aligned}
     &~\pa_x W + \operatorname{Ric} -(\pa_t -\pa^\gamma\pa_\gamma) \mathcal{H}^k  \\
 = &~ \Big(\sum_{k_2\le k_1-4}+\sum_{\substack{|k_1-k_2|\le 4 \\ |k_{1,2}-k|\ge 60 }}\Big)L(\lambda_{k_1},\lambda_{k_2}) \\
 & + \sum_{(k_1,k_2)\in \mathcal{S}_k}L\big(|D|^{-1}\lambda_{k_1}, P_{k_2}(\mathscr{A}\cdot \lambda)\big)  +  L \big( |D|^{-2}\lambda_{k_1}, P_{k_2}[(\mathscr{A}^2,\pa_x\mathscr{A},\lambda^2)\cdot  \lambda] \big)
  \end{aligned}
\end{equation}
For $j\le k-60$, applying $P_j$ to both sides of \eqref{tildeh_eq}, we obtain
\begin{equation}
  \begin{aligned}
      &~(\pa_t -\Delta)\tilde{h}_j \\
      =&~ L(h_{\le j}, 2^{2j}h_j) + L((\mathscr{A},\lambda)_{\le j}, (\mathscr{A},\lambda)_j) + \Pi_j(\mathscr{A},\mathscr{A}) + \sum_{\substack{|k_1-k_2|\le 4 \\ |k_{1,2}-k|\ge 60 }} L(\lambda_{k_1},\lambda_{k_2})\\
      & + \sum_{(k_1,k_2)\in \mathcal{S}_k}L\big(|D|^{-1}\lambda_{k_1}, P_{k_2}(\mathscr{A}\cdot \lambda)\big)  +  L \big( |D|^{-2}\lambda_{k_1}, P_{k_2}[(\mathscr{A}^2,\pa_x\mathscr{A},\lambda^2)\cdot  \lambda] \big),
  \end{aligned}
\end{equation}
By Duhamel's formula,
\begin{align}
  &~2^k \pa_x \tilde{h}_j\notag\\
   = &~ 2^k e^{t\Delta}(\pa_x h_{0,j} - \pa_x\mathcal{H}_j^k|_{t=0}) + \int_0^t 2^{2j+k}e^{s\Delta}P_j L(h_{\le j}, \pa_x h_j)(t-s)\mathrm{d}s \label{flux_A_exp_1}\\
   & + \int_0^t  2^{j+k}e^{s\Delta}P_j\Big[L((\mathscr{A},\lambda)_{\le j}, (\mathscr{A},\lambda)_j) \label{flux_A_exp_2}\\
   &\qquad\qquad\qquad\qquad+ \Pi_j(\mathscr{A},\mathscr{A}) + \sum_{\substack{|k_1-k_2|\le 4 \\ |k_{1,2}-k|\ge 60 }} L(\lambda_{k_1},\lambda_{k_2})\Big](t-s)\mathrm{d}s \label{flux_A_exp_3}\\
   & +  \int_0^t 2^{j}e^{s\Delta}P_j\Big[\sum_{(k_1,k_2)\in \mathcal{S}_k}L\big( \lambda_{k_1}, P_{k_2}(\mathscr{A}\cdot \lambda)\big) \label{flux_A_exp_4}\\
   & \qquad\qquad\qquad\qquad\qquad\ +  L \big( |D|^{-1}\lambda_{k_1}, P_{k_2}[(\mathscr{A}^2,\pa_x\mathscr{A},\lambda^2)\cdot  \lambda] \big)\Big](t-s)\mathrm{d}s.\label{flux_A_exp_5}
\end{align}
We make the following bootstrap assumption: let $C_1\in[1,+\infty)$ be the smallest constant such that
  	\begin{equation} \label{Est_lamb_4}
		\sup_{y  }\sup_{s \in \mathbb{R}_+} \| \boldsymbol{1}_{\{t\ge s\}}|\lambda_{k,\gamma}^y|^2 \cdot  2^k \pa_x \tilde{h}_j(t-s)\, \|_{L^1_{t,x}} \leq C_1  C^4 \varepsilon  ( 2^{-(\frac d2-1 + \sigma)k} c_k(\sigma))^2 c_j.
	\end{equation}

To control the first term in \eqref{flux_A_exp_1} we need bounds for both $\| \pa^2_x \mathcal{H}_j^k|_{t=0}\|_{L^1_{x_\gamma}L^2_{\hat{x}_\gamma}}$ and $\|\pa^2_x h_{0,j}\|_{L^1_{x_\gamma}L^2_{\hat{x}_\gamma}}$.  The former is estimated as follows:
 \begin{align}
  &~\| \pa^2_x\mathcal{H}_j^k|_{t=0}\|_{L^1_{x_\gamma}L^2_{\hat{x}_\gamma}}\notag\\
  \lesssim &~ 2^{\frac{d+3}{2}j-2k}\sum_{|k_{1,2}-k|\le 60}\|L(\lambda_{0,k_1},\lambda_{0,k_2})\|_{L^1_x}\notag\\
  \lesssim &~ 2^{\frac{d+3}{2}j-2k} \cdot  2^{-(d-2)k}\cdot c_{k}^2 \notag\\
  \lesssim &~ \varepsilon c_j 2^{-(\frac{d}{2}-\frac{3}{2})j}.
 \end{align}
For the latter we apply $P_j$ to the elliptic equation for $h_0$ obtained from the harmonic coordinate condition \eqref{eq_g_0}:
\[
(1-h_0^{\alpha\beta}\mathscr{R}_\alpha\mathscr{R}_\beta)\Delta (h_0)_j = [P_j, h_0\partial_x^2]h_0 + P_j(\partial_x h_0\cdot\partial_x h_0).
\]
By the smallness of $h_0$ and the boundedness of the composition of the Riesz transforms with the Littlewood--Paley projection, the operator $(1-h_0^{\alpha\beta}\mathscr{R}_\alpha\mathscr{R}_\beta\circ P_{|\cdot-j|\le2})$ is invertible.  Consequently,
\begin{align}
\|\partial_x^2 (h_0)_j\|_{L^1_{x_\gamma}L^2_{\hat{x}_\gamma}}
&\lesssim \|\Delta (h_0)_j\|_{L^1_{x_\gamma}L^2_{\hat{x}_\gamma}}
\lesssim \|(1-h_0^{\alpha\beta}\mathscr{R}_\alpha\mathscr{R}_\beta\circ P_{|\cdot-j|\le2})\Delta (h_0)_j\|_{L^1_{x_\gamma}L^2_{\hat{x}_\gamma}} \notag\\
&\lesssim \|[P_j, h_0\partial_x^2]h_0 + P_j(\partial_x h_0\cdot\partial_x h_0)\|_{L^1_{x_\gamma}L^2_{\hat{x}_\gamma}}
\lesssim \varepsilon 2^{- \frac{d-3}{2}  j}c_j.
\end{align}

Expanding $\lambda_k = \sum_{\gamma = 1}^d \lambda_{k,\gamma}$ and applying Lemma~\ref{lem:heat_divcurl}, we obtain
\begin{align}
 	&\sup_{y }\sup_{s \in \mathbb{R}_+} \| \boldsymbol{1}_{\{t\ge s\}} |\lambda_{k,\gamma}^y|^2 \cdot  2^k e^{(t-s)\Delta}(\pa_x h_{0,j} - \pa_x \mathcal{H}^k_j|_{t=0})\, \|_{L^1_{t,x}}\notag\\
  \lesssim &~   2^{ \frac{d-3}{2} j+k}\sup_{y  } \sup_{s \in \mathbb{R}_+} \| \boldsymbol{1}_{\{t\ge s\}}  \|   e^{(t-s)\Delta}(\pa^2_x h_{0,j}-\pa^2_x \mathcal{H}^k_j|_{t=0}) \|_{L^2_{\hat{x}_\gamma}}\cdot\|\lambda_{k,\gamma}^y\|_{L^2_{\hat{x}_\gamma}}^2  \, \|_{L^1_{t,x_\gamma}}\notag\\
  \lesssim &~ 2^{ \frac{d-3}{2} j }  \|  ( \pa^2_x h_{0,j}, \pa^2_x \mathcal{H}_j^k|_{t=0} ) \|_{L^1_{x_\gamma}L^2_{\hat{x}_\gamma}} (\|\lambda_k\|^2_{L^\infty_t L^2_x} + \|2^{-k}\mathcal{M}(\lambda_{k,\gamma})\|_{L^1_{t,x}}),
\end{align}
which is acceptable.

The contribution of the second term in \eqref{flux_A_exp_1} is estimated directly from the bootstrap assumption:
\begin{align}
 	&\sup_{y }\sup_{s \in \mathbb{R}_+} \| \boldsymbol{1}_{\{t\ge s\}}|\lambda_{k,\gamma}^y|^2 \cdot  \int_0^{t-s}  2^{2j}e^{s'\Delta}P_j L(h_{\le j}, 2^k \pa_x \tilde{h}_{j})(t-s-s')\mathrm{d}s'\, \|_{L^1_{t,x}}\notag\\
  \lesssim &~  \|h_{\le j}\|_{L^\infty_{t,x}}  \sup_{y }\sup_{\tau \in \mathbb{R}_+}\| \boldsymbol{1}_{\{t\ge \tau\}} |\lambda_{k,\gamma}^y|^2    2^k \pa_x \tilde{h}_{j} (t-\tau)\, \|_{L^1_{t,x}}\notag\\
  \lesssim &~ C_1 \varepsilon \cdot C^4 \varepsilon ( 2^{-(\frac d2-1 + \sigma)k} c_k(\sigma))^2 c_j.
\end{align}

For \eqref{flux_A_exp_2}--\eqref{flux_A_exp_3}, the estimate is analogous to that of \eqref{A_exp_1}; we single out the following term, the others being analogous:
\begin{align}
 &\sum_{k_{1,2}\le k-50}\sup_{y}\sup_{s \in \mathbb{R}_+}2^{k-j}\big\|\boldsymbol{1}_{\{t\ge s\}}\|P_j L(\lambda_{k_1},\lambda_{k_2})(t-s)\|_{L^\infty_{\hat{x}_\gamma}}\|\lambda_{k,\gamma}\|^2_{L^2_{\hat{x}_\gamma}} \big\|_{L^1_{t,x_\gamma}}\notag\\
 \lesssim &~ \sup_{y}\sup_{s \in \mathbb{R}_+}\sum_{k_{1,2}\le k-50}2^{k+ \frac{d-3}{2} j + \frac{d-1}{2}\min(k_1,k_2)}\notag\\
 &\qquad \qquad\qquad \qquad \cdot\big\|\boldsymbol{1}_{\{t\ge s\}}\|\lambda_{k_1}(t-s)\|_{L^2_{\hat{x}_\gamma}} \|\lambda_{k_2}(t-s) \|_{L^2_{\hat{x}_\gamma}}\|\lambda_{k,\gamma}\|^2_{L^2_{\hat{x}_\gamma}} \big\|_{L^1_{t,x_\gamma}}\notag\\
 \lesssim &~ \Big(\sum_{j'\le j }C^4 c_{j'}c_j 2^{-\frac{1}{2}|j-j'|}  + \sum_{j\le k_{1,2}\le k-50} 2^{-\frac{1}{2}|j-\max(k_1,k_2)|}c_{k_1}c_{k_2}\Big) (2^{-(\frac{d}{2}-1+\sigma)k}c_k(\sigma))^2\notag\\
 \lesssim &~ C^4 \varepsilon c_j (2^{-(\frac{d}{2}-1+\sigma)k}c_k(\sigma))^2.
\end{align}
which is acceptable.

The contribution of \eqref{flux_A_exp_4}--\eqref{flux_A_exp_5} is controlled by taking the $L^1_t L^\infty_x$ norm and then using the $L^1_tL^2_x$ estimates \eqref{cubic} and \eqref{quintic}.
In summary, we have proved
  	\begin{equation} \label{Est_lamb_4_1}
		\sup_{y  }\sup_{s \in \mathbb{R}_+} \| \boldsymbol{1}_{\{t\ge s\}}|\lambda_{k,\gamma}^y|^2 \cdot  2^k \pa_x \tilde{h}_j(t-s)\, \|_{L^1_{t,x}} \lesssim \langle C_1 \varepsilon\rangle  C^4 \varepsilon  ( 2^{-(\frac d2-1 + \sigma)k} c_k(\sigma))^2 c_j.
	\end{equation}
  Hence $\langle C_1 \varepsilon\rangle\lesssim C_1$, and since $\varepsilon \ll 1$ the minimality of $C_1$ gives $C_1 \lesssim 1$.

For the bound \eqref{Est_lamb_d3} we again argue by a bootstrap argument: let $C_2\in[1,+\infty)$ be the smallest constant such that
   \begin{equation}\label{Est_lamb_d3_1}
 \sup_{y}\sup_{s \in \mathbb{R}_+} \| \boldsymbol{1}_{\{t\ge s\}}  |\lambda^y_{k,\gamma}|^2 \cdot \pa_x A_{\operatorname{cf},j}(t-s)  \|_{L^1_{t,x}}\lesssim C_2 C^4 \varepsilon^2 (2^{-(\frac {d}{2} -1  + \sigma)k} c_k(\sigma))^2c_j,
\end{equation}
  Using \eqref{para_A_cf}, and recalling that the initial datum $A_{\operatorname{cf}}|_{t=0}$ is determined by the Coulomb condition \eqref{Cou} and the harmonic coordinate condition \eqref{Gam}, we write
  \begin{align}
    &~\pa_x A_{\operatorname{cf},j}\notag \\
    =&~ \int^t_0  2^{2j}e^{s\Delta}P_j L(h_{\le j},\pa_x A_j)(t-s) \mathrm{d}s\label{A_exp_3}\\
 & + \int^t_0 2^j e^{s\Delta}P_j \big[L(\pa_x \mathscr{A}_{\le  j},  \mathscr{A}_j)  + \Pi_j(\mathscr{A}, \pa_x\mathscr{A})+ P_j (   \mathscr{A}\cdot\lambda^2 )\big] (t-s) \label{A_exp_4}\mathrm{d}s,
  \end{align}
The contribution of \eqref{A_exp_3} is controlled by the bootstrap assumption \eqref{Est_lamb_d3_1}, in much the same way as \eqref{flux_A_exp_1}:
\begin{align}
  &~\sup_{y }\sup_{s \in \mathbb{R}_+} \| \boldsymbol{1}_{\{t\ge s\}}|\lambda_{k,\gamma}^y|^2 \cdot  \int_0^{t-s}  2^{2j}e^{s'\Delta}P_j L(h_{\le j},  \pa_x A_{\operatorname{cf},j})(t-s-s')\mathrm{d}s'\, \|_{L^1_{t,x}}\notag\\
  \lesssim &~  \|h_{\le j}\|_{L^\infty_{t,x}}  \sup_{y }\sup_{\tau \in \mathbb{R}_+}\| \boldsymbol{1}_{\{t\ge \tau\}} |\lambda_{k,\gamma}^y|^2  \pa_x A_{j} (t-\tau)\, \|_{L^1_{t,x}}\notag\\
  \lesssim &~  C_2 \varepsilon \cdot  C^4 \varepsilon^2 (2^{-(\frac {d}{2} -1  + \sigma)k} c_k(\sigma))^2c_j, \label{Estia_2}
\end{align}

The contribution of \eqref{A_exp_4} is estimated by taking the $L^1_tL^\infty_x$ norm once again and then using the $L^1_tL^2_x$ estimates \eqref{AdA} and \eqref{cubic}.

Arguing as in the proof of \eqref{Est_lamb_4_1}, we obtain $\langle C_2 \varepsilon\rangle\lesssim C_2$, hence $C_2 \lesssim 1$; this proves \eqref{Est_lamb_d3}.

For the bound \eqref{Est_lamb_d4} we argue similarly, but this time we decompose $\lambda_k = \sum_{\gamma = 1}^d \lambda_{k,\gamma}$ and make the analogous bootstrap assumption: let $C_3\in[1,+\infty)$ be the smallest constant such that

    \begin{equation}\label{Est_lamb_d4_1}
 \sup_{y}\sup_{s \in \mathbb{R}_+} \| \boldsymbol{1}_{\{t\ge s\}}  |\lambda^y_{k,\gamma}|^2 \cdot 2^k A_{\operatorname{cf},j}(t-s)    \|_{L^1_{t,x}}\lesssim C_3 C^4 \varepsilon^2 (2^{-(\frac {d}{2} -1  + \sigma)k} c_k(\sigma))^2c_j,
\end{equation}
where we used that, with the initial datum treated as above,
  \begin{align}
    &~2^k A_{\operatorname{cf},j}\notag \\
    =&~ \int^t_0  2^{2j+k}e^{s\Delta}P_j L(h_{\le j},  A_j)(t-s) \mathrm{d}s\label{A_exp_5}\\
 & + \int^t_0 2^k e^{s\Delta}P_j \big[L(\pa_x \mathscr{A}_{\le  j},  \mathscr{A}_j)  + \Pi_j(\mathscr{A}, \pa_x\mathscr{A})+ P_j (   \mathscr{A}\cdot\lambda^2 )\big] (t-s) \label{A_exp_6}\mathrm{d}s,
  \end{align}
The contribution of \eqref{A_exp_5} is controlled by the bootstrap assumption in the same way as \eqref{Estia_2}, and we shall not repeat the argument.

For the contribution of \eqref{A_exp_6} it suffices to consider the second and third terms, these being the most difficult cases.  The second term is controlled by
\begin{align}
  &~\sum_{j\le j' } 2^{k+j'-2j}\sup_{y}\sup_{s \in \mathbb{R}_+}\big\|\|P_j L(\mathscr{A}_{j'},\mathscr{A}_{j'})(t-s)\|_{L^\infty_{\hat{x}_\gamma}}\|\lambda_{k,\gamma}^y\|^2_{L^2_{\hat{x}_\gamma}}\big\|_{L^1_{t,x_\gamma}}\notag\\
  \lesssim &~ \sum_{j\le j' } 2^{(d-3)j + j'+k}\sup_{y}\sup_{s \in \mathbb{R}_+} \big\|\| \mathscr{A}_{j'} (t-s)\|_{L^2_{\hat{x}_\gamma}}\|\lambda^y_{k,\gamma}\|_{L^2_{\hat{x}_\gamma}}\big\|^2_{L^2_{t,x_\gamma}}\notag\\
  \lesssim &~ \sum_{j\le j' } C^4 2^{-|j-j'|} c_{j'}^2 (2^{-(\frac {d}{2} -1  + \sigma)k} c_k(\sigma))^2\notag\\
  \lesssim &~ C^4 \varepsilon   (2^{-(\frac {d}{2} -1  + \sigma)k} c_k(\sigma))^2c_{j }.
\end{align}
For the third term the argument is analogous: it suffices to regard $\lambda^2$ as $\pa_x \mathscr{A}$, since the two enjoy the same space-time estimates; the only difference is that, when applying the bilinear estimate, one decomposes the high frequency, $\lambda_{k_{\operatorname{max}}} = \sum_{\gamma = 1}^d \lambda_{k_{\operatorname*{max}},\gamma}$.
Finally, arguing as above, we obtain $C_3 \lesssim 1$, which proves \eqref{Est_lamb_d4}.

This completes the proof.
\end{proof}

We shall also need the following estimates for the nonlinearities of the difference variables.
\begin{lemma}\label{flux_delta}
\begin{align}
	 &\| L(2^{-k}\delta\lambda_k, 2^{-k}P_k((\delta\lambda,\delta\mathscr{A})  \cdot(\lambda,\mathscr{A})  \cdot\lambda )) \|_{L^1_{t,x}} \lesssim C^4 \varepsilon^2 (2^{-(\frac {d}{2} -1)k} d_k )^2,   \label{Est_dlamb_1}\\
	&\|2^{-k}\delta\lambda_k \cdot 2^{-k}\delta\Pi_k(\pa_x \mathscr{A},\lambda  )\|_{L^1_{t,x}}   +\|  2^{-k}\delta\lambda_k \cdot 2^{-k}\delta L(2^k \mathscr{A}_{\ge k-60}, \lambda_k)\|_{L^1_{t,x}}\label{Est_dlamb_2}\\
  &\qquad\qquad  \lesssim C^3 \varepsilon (2^{-(\frac {d}{2} -2  )k} d_k )^2,  \notag\\
	&\|   L( 2^k (A_{\operatorname{df}},W)_{\le k-60}, 2^{-k}\delta\lambda_k,2^{-k}\delta\lambda_k) \|_{L^1_{t,x}} \lesssim C^4 \varepsilon^2 (2^{-(\frac {d}{2} -1)k} d_k )^2, \label{Est_dlamb_3}\\
  &\|  L( 2^k \pa_x h_{\le k-60}, 2^{-k}\delta\lambda_k, 2^{-k}\delta\lambda_k) \|_{L^1_{t,x}} \lesssim C^4 \varepsilon^2 (2^{-(\frac {d}{2} -1)k} d_k )^2, \label{Est_dlamb_4}\\
  &\|L(  \delta\mathscr{A}_{\le k-60}, 2^{-k}\delta \lambda_k, \lambda_k)\|_{L^1_{t,x}} \lesssim C^4 \varepsilon^2 (2^{-(\frac {d}{2} -1)k} d_k )^2, \label{Est_dlamb_5}
  \\
     &\|L(2^k \delta h_{\le k}, 2^{-k}\delta \lambda_k, \lambda_k)\|_{L^1_{t,x}} \lesssim C^4 \varepsilon^2 (2^{-(\frac {d}{2} -1)k} d_k )^2. \label{Est_dlamb_6}
\end{align}
In particular, for $d\ge4$ we have
\begin{align}
  & \|  L(2^k A_{\operatorname{cf},\le k-60}, 2^{-k}\delta\lambda_k,2^{-k}\delta\lambda_k) \|_{L^1_{t,x}}\lesssim C^4 \varepsilon^2 (2^{-(\frac {d}{2} -1)k} d_k )^2,\label{Est_dlamb_d4}
\end{align}
while for $d=3$ we have
\begin{align}
  &\| L(\pa_x A_{\operatorname{cf},\le k-60}, 2^{-k}\delta\lambda_k,2^{-k}\delta\lambda_k) \|_{L^1_{t,x}}\lesssim C^4 \varepsilon^2 (2^{-(\frac {d}{2} -1)k} d_k )^2,\label{Est_dlamb_d3}
\end{align}
\end{lemma}

\begin{proof}
The bounds \eqref{Est_dlamb_1}--\eqref{Est_dlamb_6}, \eqref{Est_dlamb_d4} and \eqref{Est_dlamb_d3} are the difference analogues of \eqref{Est_lamb_1}--\eqref{Est_lamb_5}, \eqref{Est_lamb_d4} and \eqref{Est_lamb_d3}, and are proved in the same way, with the difference bounds used for the factors involving the difference variables.  In particular, as observed in the proof of \eqref{dlala4}, the quantity $2^{-k}\delta\lambda_k$ behaves like $\lambda_k$; accordingly, the bounds \eqref{Est_dlamb_3}, \eqref{Est_dlamb_4}, \eqref{Est_dlamb_d3} and \eqref{Est_dlamb_d4} are proved exactly as \eqref{Est_lamb_3}, \eqref{Est_lamb_5}, \eqref{Est_lamb_d3} and \eqref{Est_lamb_d4}, respectively.

For the remaining bounds, it suffices to consider the most difficult case, in which the difference variable is localized at the highest frequency; the other cases are analogous (cf.\ the strategy in the proof of Lemma~\ref{delta_prim_bound}), and we omit the repetition.

The bound \eqref{Est_dlamb_1} follows by applying \eqref{dlala} to $(2^{-k}\delta \lambda_k, (\delta\lambda, \delta \mathscr{A}))$ and the original \eqref{lala} to the remaining pair, and then summing over the frequencies.

The bound \eqref{Est_dlamb_2} is analogous to \eqref{Est_lamb_2}, except that \eqref{dBA2}, \eqref{dBA4a} and \eqref{dBA4b} are used instead.

The bound \eqref{Est_dlamb_5} follows by arguments analogous to those in \eqref{Est_lamb_3} and \eqref{Est_lamb_d3}.  We first consider the contribution of $(\delta A_{\operatorname{df}},\delta W)$; using \eqref{A_exp_0}, we have
 \begin{align}
 &~ (\delta A_{\operatorname{df}},\delta W)_{\le k-60}\notag\\
   = &~  \sum_{\substack{
  k_{1,2}\le k - 50  }} |D|^{-1}L( \delta \lambda_{k_1},\lambda_{k_2}) +\sum_{\substack{|k_1-k_2|\le 4\\
  k_{1,2}\ge k + 50  }}  |D|^{-1}L( \delta \lambda_{k_1},\lambda_{k_2})\label{dA_exp_1} \\
   & +  \sum_{(k_1,k_2)\in\mathcal{S}_k}
            \pa_x  L(|D|^{-1}\delta\lambda_{k_1}\cdot  \lambda_{k_2})
              \label{dA_exp_2}\\
   & +  \sum_{(k_1,k_2)\in\mathcal{S}_k} \Bigl[
             L\big(|D|^{-1}\delta\lambda_{k_1}, P_{k_2} ( \lambda \cdot\mathscr{A} )\big)
            +  L\big(|D|^{-1}\lambda_{k_1}, P_{k_2} ( \delta\lambda \cdot\mathscr{A} +\lambda \cdot\delta\mathscr{A})\big)
            \Bigr]   \label{dA_exp_3}.
\end{align}
The contribution of \eqref{dA_exp_1} is estimated as that of \eqref{A_exp_1}, with \eqref{dBA4a} and \eqref{dBA4b} used instead of \eqref{BA4}; since the difference variable carries the highest frequency, the resulting dyadic sum is controlled through the $\ell^1$-summability of the frequency envelopes, exactly as in \eqref{dAdA}.
The contribution of \eqref{dA_exp_2} is identical to that of the first term in \eqref{A_exp_2}, because $|D|^{-1}\delta\lambda$ behaves like $\lambda$; we refer to the argument above.  Finally, the contribution of \eqref{dA_exp_3} is estimated by taking the $L^1_tL^\infty_x$ norm and then using the $L^1_tL^2_x$ estimate \eqref{highcubic_2}.

For the remaining contribution of $\delta A_{\operatorname{cf}}$, we expand \eqref{para_dA_cf} and obtain
 \begin{align}
  &~\delta A_{j,\mathrm{cf}} \notag\\
  =&~\int_0^t e^{s\Delta} 2^{2j}(L( h_{\le j},  \delta A_j) + L(\delta h_{\le j}, A_j) )(t-s)\mathrm{d}s\label{exp_Acf_1}\\
 &+\int_0^t e^{s\Delta}\Big[L(\pa_x \delta \mathscr{A}_{\le  j}, \mathscr{A}_j) + L(\pa_x \mathscr{A}_{\le  j}, \delta\mathscr{A}_j)+ \Pi_j(\delta \mathscr{A}, \pa_x\mathscr{A})\label{exp_Acf_2}\\
 &\qquad\qquad\  + P_j (   \delta\mathscr{A}\cdot\lambda^2  + \lambda \cdot \delta \lambda \cdot \mathscr{A})\Big](t-s)\mathrm{d}s\label{exp_Acf_3},
 \end{align}

As in the proof of \eqref{Est_lamb_5}, we again argue by a bootstrap argument: let $C_4\in[1,+\infty)$ be the smallest constant such that
   \begin{equation}\label{Est_dlamb_d3_1}
 \sup_{y,y'}\sup_{s \in \mathbb{R}_+} \| \boldsymbol{1}_{\{t\ge s\}}  \lambda^y_{k,\gamma}\cdot 2^{-k}\delta\lambda_{k,\gamma}^{y'} \cdot \delta A_{\operatorname{cf},j}(t-s)  \|_{L^1_{t,x}}\lesssim C_4 C^3 \varepsilon^2 2^{-(d-2)k} c_j   d^2_k.
\end{equation}

The contribution of the first term in \eqref{exp_Acf_1} is controlled by the bootstrap assumption \eqref{Est_dlamb_d3_1}, in the same way as \eqref{Estia_2}, and we shall not repeat the argument;
The contribution of the second term in \eqref{exp_Acf_1} is estimated analogously; repeating the proof of \eqref{Est_lamb_3}, we obtain
 \begin{equation}\label{Estia_3}
   \sup_{y,y'}\sup_{s \in \mathbb{R}_+} \| \boldsymbol{1}_{\{t\ge s\}}  \lambda^y_{k,\gamma}\cdot 2^{-k}\delta\lambda_{k,\gamma}^{y'} \cdot  A_{\operatorname{cf},j}(t-s)  \|_{L^1_{t,x}}\lesssim   C^3 \varepsilon^2 2^{-(d-2)k-j} c_j  d^2_k,
 \end{equation}
Using \eqref{Estia_3}, together with the $L^\infty_{t,x}$ norm of $\delta h_{\le j}$, the extra derivative is absorbed by the inverse derivative $2^{-j}$ on the right-hand side.

The contribution of \eqref{exp_Acf_2} is estimated by taking the $L^1_tL^\infty_x$ norm and then using the $L^1_tL^2_x$ estimates \eqref{dAdA_2} and \eqref{highcubic_2}.

For the contribution of \eqref{exp_Acf_3}, it suffices to consider the most difficult case, in which the difference variable in $(\delta\lambda\cdot\lambda\cdot\mathscr{A})$ carries the highest frequency and exceeds $j$; the reasoning is the same as above.  Furthermore, it suffices to consider the case where $\lambda$ and $\delta\lambda$ have comparable frequencies, since $\mathscr{A}$ enjoys better space-time regularity; it then remains to control the contribution of the following term
\begin{equation}
 \sum_{\substack{|j_1-j_2|\le 2 \\ j_1\ge j_3\ge j}} P_j (\delta \lambda_{j_1}\cdot \lambda_{j_2}\cdot \mathscr{A}_{j_3}).\label{exp_Acf_3_1}
\end{equation}
Let $d=3$, so $d_{j_1}=d_k$. On a smooth interval with finite normalized sources, let $C_6\ge1$ be the least constant in the bounds of Proposition~\ref{prop:source} at $\sigma=0$ in place of $C^4\varepsilon^2$, for both $P_k$ and the conic projections. For $|j_1-j_2|\le2$ and $j\le j_3\le j_1$, Cauchy--Schwarz on the exterior waves, Lemma~\ref{lem:heat_divcurl} and \eqref{energy} give
\begin{align}
 &~\sup_{y,y'}\sup_{s\ge0}\bigg\|\boldsymbol{1}_{\{t\ge s\}}\lambda_{k,\gamma}^y\,2^{-k}\delta\lambda_{k,\gamma}^{y'}
 \int_0^{t-s}e^{s'\Delta} P_j(\delta\lambda_{j_1}\cdot\lambda_{j_2}\cdot\mathscr A_{j_3})(t-s-s')\,\mathrm ds'\bigg\|_{L^1_{t,x}}\notag\\
 \lesssim&~C_6 2^{-k}c_kd_k\,2^{j-k}\int_0^\infty
 \|P_j(\delta\lambda_{j_1}\cdot\lambda_{j_2}\cdot\mathscr A_{j_3})(s')\|_{L^1_{x_\gamma}L^2_{\hat x_\gamma}}\,\mathrm ds'\notag\\
 \lesssim&~C_6 2^{-k}c_kd_k\,2^{2j-k}\|P_j(\delta\lambda_{j_1}\lambda_{j_2}\mathscr A_{j_3})\|_{L^1_{t,x}}\notag\\
 \lesssim&~ C_6 2^{-k}c_kd_k\,2^{2j-k}
 \|\delta\lambda_{j_1}\|_{L^4_{t,x}}\|\lambda_{j_2}\|_{L^4_{t,x}}\|\mathscr A_{j_3}\|_{L^2_{t,x}}\notag\\
 \lesssim&~ C_6 2^{-k}c_kd_k\,2^{2j-k}
 (C2^{3j_1/4}d_{j_1})(C2^{-j_2/4}c_{j_2})(C2^{-3j_3/2}c_{j_3})\notag\\
 \lesssim&~C^3 C_6 2^{-k}c_kd_k^2c_{j_2}c_{j_3}\,2^{2j-k+j_1/2-3j_3/2}.\notag
\end{align}
Here we used \eqref{dBA5}, \eqref{BA5} and \eqref{BA2}, with \eqref{L2W} for $W_{j_3}$. For $j_1\le k+52$,
\[
 2^{2j-k+j_1/2-3j_3/2}\lesssim
 2^{-(k-j)/2-3(j_3-j)/2}\le 2^{-2(j_3-j)/3}.
\]
For $j_1>k+52$, pair the four waves as in \eqref{Est_lamb_1}. By \eqref{dlala2}, \eqref{def_B} and Lemma~\ref{lem:lambda_divcurl}, the same summand is bounded by
\begin{align}
 &~\int_0^\infty e^{-c2^{2j}s'}\,
 \sup_y\|\boldsymbol{1}_{\{t\ge s+s'\}}\lambda_{k,\gamma}^y(t)\,\delta\lambda_{j_1}(t-s-s')\|_{L^2_{t,x}}\notag\\
 &\qquad{}\cdot\sup_{y'}\|\boldsymbol{1}_{\{t\ge s+s'\}}2^{-k}\delta\lambda_{k,\gamma}^{y'}(t)\,\lambda_{j_2}(t-s-s')\|_{L^2_{t,x}}\,
 \|\mathscr A_{j_3}\|_{L^\infty_{t,x}}\,\mathrm ds'\notag\\
 \lesssim&~\int_0^\infty e^{-c2^{2j}s'}\,\mathrm ds'\,
 (C^2 2^{-(j_1-k)/2+j_1/2}c_kd_{j_1})
 (C_6 2^{-(j_2-k)/2-j_2/2}d_kc_{j_2})\,\|\mathscr A_{j_3}\|_{L^\infty_{t,x}}\notag\\
 \lesssim&~2^{-2j}(C^2 2^{k/2}c_kd_{j_1})(C_6 2^{k/2-j_2}d_kc_{j_2})(C2^{j_3}c_{j_3})\notag\\
 \lesssim&~C^3 C_6 2^{-k}c_kd_k^2c_{j_2}c_{j_3}\,2^{2k-j_1+j_3-2j}.\notag
\end{align}
Taking the two bounds to the powers $2/3$ and $1/3$, we obtain
\[
 (2^{2j-k+j_1/2-3j_3/2})^{2/3}(2^{2k-j_1+j_3-2j})^{1/3}
 = 2^{-2(j_3-j)/3}.
\]
Thus the contribution of \eqref{exp_Acf_3_1} to \eqref{Est_dlamb_d3_1} is bounded by
\begin{align}
 &~C^3 C_6 2^{-k}c_kd_k^2\sum_{j_3\ge j}2^{-2(j_3-j)/3}c_{j_3}
 \sum_{j_1\ge j_3}\sum_{|j_2-j_1|\le2}c_{j_2}\notag\\
 \lesssim&~C^3 C_6\varepsilon^2 2^{-k}c_jd_k^2
 \sum_{j_3\ge j}2^{-(2/3-\delta')(j_3-j)}
 \lesssim C^3 C_6\varepsilon^2 2^{-k}c_jd_k^2 .
 \label{eq:flux-cubic-sum}
\end{align}
For $d\ge4$ use the first bound with $L^4_{t,y}L^2_z$, $L^4_{t,y}L^2_z$ and $L^2_{t,y}L^\infty_z$ on the three inputs. Bernstein in the $d-3$ variables $z$ gives
\begin{align}
 &~2^{(d-1)j-k}\|P_j(\delta\lambda_{j_1}\lambda_{j_2}\mathscr A_{j_3})\|_{L^1_{t,x}}\notag\\
 \lesssim&~C^3d_{j_1}c_{j_2}c_{j_3}\,
 2^{(d-1)j-k-(d/2-9/4)j_1-(d/2-5/4)j_2-3j_3/2}\notag\\
 \lesssim&~C^3 d_kc_{j_2}c_j\,
 2^{-(1-\delta')(k-j)-(d-7/2-\delta')(j_1-j)-(3/2-\delta')(j_3-j)}.\notag
\end{align}
Multiplying by $C_6 2^{-(d-2)k}c_kd_k$ and summing gives $C^3C_6\varepsilon^2 2^{-(d-2)k}c_jd_k^2$.

For the bound \eqref{Est_dlamb_6}, as before, we split off the high-high part of $h$ and write $h=\tilde h+\mathcal{H}$, where $\mathcal{H}$ is the full high-high contribution to $\operatorname{Im}(\bar\psi\lambda_{\alpha\beta})$, i.e.\ the sum of $\mathcal{H}[k_1,k_2]$ over all pairs with $|k_1-k_2|\le4$ (see Lemma~\ref{hidden_curvature}). Since all high-high pairs are collected in $\mathcal{H}$, the nolinearities in equation for $\tilde h$ contains no high-high interaction of the form $\lambda\cdot \lambda$, which is essential for the estimates below. Applying $\delta$ gives $\delta h=\delta\tilde h+\delta\mathcal{H}$, where
\begin{equation}
  \delta \mathcal{H} =  \sum_{|k_1 -k_2|\le 4}
             2^{-k} L(|D|^{-1}\delta \lambda_{k_1}\cdot  \lambda_{k_2})
\end{equation}
The contribution of $\delta\mathcal{H}$ is estimated as in \eqref{estia_1}, noting only that $|D|^{-1}\delta\lambda$ behaves like $\lambda$.

It remains to estimate $\delta\tilde h$. We use the cancellation identity of Lemma~\ref{hidden_curvature} with $\mathcal{S}=\{(k_1,k_2):|k_1-k_2|\le4\}$:
\begin{equation}
  \begin{aligned}
     &~\pa_x W + \operatorname{Ric} -(\pa_t -\pa^\gamma\pa_\gamma) \mathcal{H}  \\
 = &~  \sum_{k_2\le k_1-4} L(\lambda_{k_1},\lambda_{k_2}) \\
 & + \sum_{|k_1-k_2|\le 4}L\big(|D|^{-1}\lambda_{k_1}, P_{k_2}(\mathscr{A}\cdot \lambda)\big)  +  L \big( |D|^{-2}\lambda_{k_1}, P_{k_2}[(\mathscr{A}^2,\pa_x\mathscr{A},\lambda^2)\cdot  \lambda] \big)
  \end{aligned}
\end{equation}

Applying $\delta$ to the equation for $\tilde h$ and expanding the result by Duhamel's formula gives
 \begin{align}
  &~2^k \delta \tilde{h}_j\notag\\
   = &~ 2^k e^{t\Delta}(\delta h_{0,j} - \delta\mathcal{H}_j|_{t=0}) + \int_0^t 2^{2j}e^{s\Delta} [L(h_{\le j}, \delta h_j)+L(\delta h_{\le j},  h_j)] (t-s)\mathrm{d}s \label{flux_dA_exp_1}\\
   & + \int_0^t  2^{ k}e^{s\Delta}P_j\Big[L((\delta\mathscr{A},\delta\lambda)_{\le j}, (\mathscr{A},\lambda)_j) +L((\mathscr{A},\lambda)_{\le j}, (\delta\mathscr{A},\delta\lambda)_j)\label{flux_dA_exp_2}\\
   &\qquad\qquad\qquad \quad+ \Pi_j(\delta\mathscr{A},\mathscr{A})  \Big](t-s)\mathrm{d}s \label{flux_dA_exp_3}\\
   & +  \int_0^t  e^{s\Delta}P_j\Big[\sum_{|k_1-k_2|\le 4}L\big( \delta\lambda_{k_1}, P_{k_2}(\mathscr{A}\cdot \lambda)\big) +L\big( \lambda_{k_1}, P_{k_2}(\delta\mathscr{A}\cdot \lambda+\mathscr{A}\cdot \delta\lambda)\big) \label{flux_dA_exp_4}\\
   & \qquad +  L \big( |D|^{-1}\delta \lambda_{k_1}, P_{k_2}[(\mathscr{A}^2,\pa_x\mathscr{A},\lambda^2)\cdot  \lambda] \big) \label{flux_dA_exp_5}\\
   & \qquad  + L \big( |D|^{-1}  \lambda_{k_1}, P_{k_2}[ \pa_x \delta\mathscr{A}\cdot \lambda + \pa_x \mathscr{A}\cdot \delta \lambda+ (\delta\mathscr{A},\delta\lambda) \cdot (\mathscr{A},\lambda)\cdot \lambda] \big)    \Big](t-s)\mathrm{d}s\label{flux_dA_exp_6}
\end{align}
 We again argue by a bootstrap argument: let $C_5\in[1,+\infty)$ be the smallest constant such that
  	\begin{equation} \label{Est_dlamb_6_1}
		\sup_{y,y'  }\sup_{s \in \mathbb{R}_+} \| \boldsymbol{1}_{\{t\ge s\}} \lambda_{k,\gamma}^y \cdot 2^{-k}\delta \lambda_k^{y'}  \cdot   2^k \pa_x \tilde{h}_j(t-s)\, \|_{L^1_{t,x}} \lesssim   C_5    C^4 \varepsilon 2^{-(d-2)k}\cdot(c_j d_k^2 + d_j c_k d_k).
    \end{equation}
The contribution of \eqref{flux_dA_exp_1} is estimated in much the same way as that of \eqref{flux_A_exp_1}: the first term is handled by Lemma~\ref{lem:heat_divcurl}, and the last two directly from the bootstrap assumption \eqref{Est_dlamb_6_1}.

The contribution of \eqref{flux_dA_exp_2} is estimated exactly as before: since $\delta$ acts on frequencies $\le j$, it may be regarded as an ordinary derivative, and the same argument applies.

The contribution of \eqref{flux_dA_exp_3} is handled as the cubic term \eqref{exp_Acf_3_1} above:
\begin{align}
 	&~\sup_{y,y'  }\sup_{s \in \mathbb{R}_+} \| \boldsymbol{1}_{\{t\ge s\}} \lambda_{k,\gamma}^y \cdot 2^{-k}\delta \lambda_{k,\gamma}^{y'}  \cdot  \int_0^{t-s}  2^k e^{s'\Delta}\Pi_j(\delta \mathscr{A},\mathscr{A})(t-s-s')\mathrm{d}s'\, \|_{L^1_{t,x}}\notag\\
    \lesssim	&~\sup_{y,y'  }\sup_{s \in \mathbb{R}_+} \int_{\mathbb{R}_+}\| \boldsymbol{1}_{\{t\ge s+s'\}} \lambda_{k,\gamma}^y \cdot 2^{-k}\delta \lambda_{k,\gamma}^{y'}  \cdot     2^k e^{(t-s-s')\Delta} \Pi_j(\delta \mathscr{A},\mathscr{A})(s')\|_{L^1_{t,x}}\mathrm{d}s'\notag\\
 \lesssim &~ \sum_{j'\ge j}\sup_{y,y'}\sup_{\tau \in \mathbb{R}_+}  C^2 2^{-(d-2)k}c_k d_k \cdot 2^{(d-1)j} \|\delta \mathscr{A}_{j'}(s')\cdot\mathscr{A}_{j'}(s')\|_{L^1_{s',x}}\notag\\
 \lesssim &~ \sum_{j'\ge j}\sup_{y,y'}\sup_{\tau \in \mathbb{R}_+}  C^2 2^{-(d-2)k}c_k d_k \cdot C^2 2^{-(d-1)(j'-j)}c_{j'}d_{j'}\notag\\
 \lesssim &~   C^4 c_j d_j c_k d_k \lesssim C^4 \varepsilon c_j d_k^2 . 
\end{align}
which is acceptable.

The contributions of \eqref{flux_dA_exp_4} and \eqref{flux_dA_exp_6} are estimated by taking the $L^1_tL^\infty_x$ norm and then using the $L^1_tL^2_x$ estimates \eqref{highcubic_2} and \eqref{highquntic}.

Finally, combining the estimates above and arguing as in the proof of \eqref{Est_lamb_4_1}, we obtain $C_5\lesssim 1$.  Summing the bootstrap bound in $j$ then yields the desired estimate \eqref{Est_dlamb_6}.
\end{proof}

Using the results above, the source terms of the Schrödinger variable can be controlled.
\begin{proposition}\label{prop:source}
    Under the bootstrap assumptions \eqref{BA1}--\eqref{BA5} and \eqref{dBA1}--\eqref{dBA5}, for $\sigma \ge 0$ we have
    \begin{align}
    &\|\mathcal{E}(\varphi_k),\ 2^{-k}\mathcal{M}(\varphi_k),\  2^{-k}\widetilde{\mathcal{M}}(\varphi_k)\|_{L^1_{t,x}} \lesssim  C^4\varepsilon^2   2^{- (d-2+2\sigma)k}c_k(\sigma)^2, \\
    &\|\mathcal{E}(2^{-k}\delta\varphi_k),\ 2^{-k}\mathcal{M}(2^{-k}\delta\varphi_k),\  2^{-k}\widetilde{\mathcal{M}}(2^{-k}\delta\varphi_k)\|_{L^1_{t,x}} \lesssim  C^4\varepsilon^2   2^{- (d-2 )k}d_k^2,
    \end{align}
\end{proposition}
\begin{proof}
The first estimate is obtained by combining the $L^1_{t,x}$ bounds of Lemma~\ref{flux} for $\varphi_k$ with the bootstrap assumptions and the definitions of $\mathcal{E}$, $\mathcal{M}$ and $\widetilde{\mathcal{M}}$; the second is its difference analogue and follows from Lemma~\ref{flux_delta}. Both lemmas precede the present proposition and are proved without invoking it, so the logical order is acyclic.
\end{proof}

\section{Proof of the main theorem}\label{sec:proof}
\subsection{Closing the bootstrap arguments}
We now close the bootstrap arguments.  Recall from Corollary~\ref{cor:boot} that it suffices to establish the theorem under the bootstrap assumptions \eqref{BA1}--\eqref{BA5} and \eqref{dBA1}--\eqref{dBA5}.  Proposition~\ref{prop:source} provides the key nonlinear estimates for the source terms of the Schrödinger variable, while Proposition~\ref{prop:K}  controls the source terms of the parabolic variables.  Using these estimates, we immediately improve the constants in the bootstrap bounds.

\begin{proposition}\label{pro:closeBA}
Under the bootstrap assumptions \eqref{BA1}--\eqref{BA5} and \eqref{dBA1}--\eqref{dBA5}, the constant $C$ can be improved to $\langle C^4\varepsilon\rangle$ in \eqref{BA1}--\eqref{BA5} and \eqref{dBA1}--\eqref{dBA5}.
\end{proposition}
\begin{proof}
For \eqref{BA1} and \eqref{BA2}, the estimates follow from the mass balance law \eqref{energy} for the Schrödinger variable and the analogous energy estimates for $(\pa_x h,A)$, the source terms being controlled by Proposition~\ref{prop:source} and Proposition~\ref{prop:K}.  For \eqref{BA3} and \eqref{BA4}, the bounds follow from the div-curl lemmas Lemma~\ref{lem:Ah_divcurl} and Lemma~\ref{lem:lambda_divcurl}.  For \eqref{BA5}, the required bound is a consequence of the interaction Morawetz estimate in Proposition~\ref{pro:Morawetz}.  The difference bounds \eqref{dBA1}--\eqref{dBA5} are established in the same way, with Lemma~\ref{delta_prim_bound} and Lemma~\ref{flux_delta} in place of Lemma~\ref{prim_bound} and Lemma~\ref{flux}, while \eqref{dBA5} again follows from Proposition~\ref{pro:Morawetz}, applied to the difference variable $\delta\varphi_k$.
\end{proof}

Since $\langle C^4\varepsilon\rangle\le C/2$ (with $C=\varepsilon^{-1/16}$), the improvement of the constant is strict, and a standard continuity argument closes the bootstrap: the bounds \eqref{BA1}--\eqref{BA5} provide the global a priori control of the solution, while \eqref{dBA1}--\eqref{dBA5} give \eqref{dependence} for $d\ge4$ and the analogous bound in $\dot B^{\frac d2-2}_{2,\infty}$ for $d=3$, hence uniqueness and continuous dependence in all dimensions $d\ge3$.

At this point the bootstrap argument has yielded the a priori bounds for smooth solutions of the differentiated system \eqref{eq_lam}, \eqref{eq_g}--\eqref{eq_A} in the modified heat gauge. It remains to explain how such solutions are obtained from the geometric initial data.

For smooth initial data, the short-time existence of a smooth solution of the SMCF in the temporal gauge $V=0$ is due to Huang and Tataru \cite{HuT25}. To pass to the modified heat gauge, we perform a change of coordinates $x\mapsto y(t,x)$ and a $\mathrm{U}(1)$ rotation of the normal frame; the dynamical variables transform as
\[
\tilde g_{\alpha\beta}=\partial_{y_\alpha}x^{\alpha'}\partial_{y_\beta}x^{\beta'}g_{\alpha'\beta'},\qquad
\tilde\lambda_{\alpha\beta}=\partial_{y_\alpha}x^{\alpha'}\partial_{y_\beta}x^{\beta'}e^{i\theta}\lambda_{\alpha'\beta'},
\]
\[
\tilde A_\alpha=\partial_{y_\alpha}x^{\alpha'}A_{\alpha'}-\partial_{y_\alpha}\theta,\qquad
\tilde B=B-\partial_t\theta.
\]
The velocity of the new coordinates is $\tilde V^\gamma=-\partial_t y^\gamma$; hence the modified heat coordinate condition $\tilde V^\gamma=\tilde g^{\alpha\beta}\tilde\Gamma^\gamma_{\alpha\beta}+\tilde W^\gamma$ is equivalent to
\[
\partial_t y^\gamma-\Delta_{\tilde g}y^\gamma=-\tilde W^\gamma,
\]
where we used $\Delta_{\tilde g}y^\gamma=-\tilde g^{\alpha\beta}\tilde\Gamma^\gamma_{\alpha\beta}$, and the initial data $y(0,x)=x+\phi(x)$ are chosen so that the harmonic coordinate condition holds at $t=0$. Similarly, since $\tilde B=B-\partial_t\theta$ and $\tilde A_\alpha=\partial_{y_\alpha}x^{\alpha'}A_{\alpha'}-\partial_{y_\alpha}\theta$, the modified heat gauge condition $\tilde B=\nabla^\alpha\tilde A_\alpha+\tilde{\mathcal{F}}$ is equivalent to
\[
\partial_t\theta-\Delta_{\tilde g}\theta=B-\nabla^\alpha A_\alpha-\tilde{\mathcal{F}},
\]
with $\theta|_{t=0}=b$ chosen so that the Coulomb condition holds at $t=0$. These linear parabolic equations with smooth coefficients are solvable on a short time interval; by construction, the resulting diffeomorphism $y(t,\cdot)$ and rotation $\theta(t,\cdot)$ satisfy
\[
\tilde V^\gamma=\tilde g^{\alpha\beta}\tilde\Gamma^\gamma_{\alpha\beta}+\tilde W^\gamma,\qquad
\tilde B=\nabla^\alpha\tilde A_\alpha+\tilde{\mathcal{F}},
\]
so the transformed solution is in the modified heat coordinates and heat gauge and evolves according to the differentiated system of Section~2.

For initial data merely in the critical Besov space, we apply this construction to a standard smooth regularization of the data (cf.\ Remark~\ref{rem:initial_data}) and use the uniform a priori bounds obtained above to pass to the limit. This yields the global solution in the critical regularity and completes the proof of Theorem~\ref{thm:boot}.
 \bibliographystyle{plain}
\bibliography{ref}
\end{document}